\documentclass{article}

\usepackage{amsmath}
\usepackage{amsthm} 
\usepackage{natbib}
\usepackage[colorlinks,citecolor=blue,urlcolor=blue,filecolor=blue,backref=page,hypertexnames=false]{hyperref}
\usepackage{graphicx}

\usepackage{float}
\usepackage{subcaption}
\usepackage{bm}
\usepackage{amsfonts}
\usepackage{enumerate}
\usepackage[capitalize,nameinlink,noabbrev]{cleveref}
\usepackage{autonum}
\usepackage{thmtools}
\usepackage{thm-restate}
\usepackage[dvipsnames]{xcolor}
\usepackage{extarrows}

\numberwithin{equation}{section}

\theoremstyle{plain}
\newtheorem{thm}{Theorem}
\numberwithin{thm}{section}

\newtheorem{lem}[thm]{Lemma}
\newtheorem{cor}[thm]{Corollary}
\newtheorem{prop}{Proposition}
\newtheorem{propA}[thm]{Proposition}

\newtheorem{rem}[thm]{Remark}

\theoremstyle{definition}

\newtheorem*{dfn*}{Definition}

\newcommand{\Asection}[1]{Section \ref*{#1} (\nameref*{#1})}

\newcounter{Bsec}\stepcounter{Bsec}

\newcommand{\mI}{I}
\newcommand{\mQ}[1]{Q_{\mathbf{x}}(#1)}
\newcommand{\bx}{\mathbf{x}}
\newcommand{\by}{\mathbf{y}}
\newcommand{\CRP}{\textrm{CRP}}
\newcommand{\bX}{\mathbf{X}}
\newcommand{\ceil}[1]{\lceil #1 \rceil}
\newcommand{\ca}{\ceil{\alpha}}
\newcommand{\mr}{m^{(r)}}
\newcommand{\Mr}{M^{(r)}}
\newcommand{\eig}[1]{\underline{\nu}_{#1}}
\newcommand{\Eig}[1]{\overline{\nu}_{#1}}
\newcommand{\Tau}{\textnormal{T}}
\newcommand{\sm}{\setminus}
\newcommand{\1}{{\bf 1}}

\newcommand{\eps}{\ensuremath{\varepsilon}}

\newcommand{\tr}{\textnormal{tr}}
\newcommand{\conv}{\textnormal{conv}}
\newcommand{\diag}[1]{\textnormal{diag}(#1)}
\newcommand{\dd}[1]{\frac{\partial}{\partial #1}}
\renewcommand{\d}[1]{\textnormal{d}#1}
\newcommand{\diam}{\textnormal{diam}\hspace{0.05cm}}
\newcommand{\dist}{\textnormal{dist}}

\newcommand{\norm}[1]{\|#1\|}
\newcommand{\znorm}[1]{\zz\|#1\zz\|}
\newcommand{\xnorm}[1]{\xx\|#1\xx\|}

\newcommand{\Var}{\textnormal{Var}}
\newcommand{\zz}{\big}
\newcommand{\xx}{\Big}
\newcommand{\re}[1]{\frac{1}{#1}}
\newcommand{\inv}[1]{#1^{-1}}

\newcommand{\E}{\mathbb{E}\,}

\newcommand{\cond}{\,|\,}

\newcommand{\iid}{\stackrel{\textrm{iid}}{\sim}}

\newcommand{\Beta}{\textnormal{Beta}}
\newcommand{\Normal}{\mathcal{N}}
\newcommand{\Exp}{\textnormal{Exp}}

\newcommand{\DP}{\textnormal{DP}}
\newcommand{\Dir}{\textnormal{Dir}}

\newcommand{\Unif}{\textnormal{Unif}}

\newcommand{\R}{\ensuremath{\mathbb{R}}}

\newcommand{\N}{\ensuremath{\mathbb{N}}}
\renewcommand{\P}{\mathbb{P}}

\newcommand{\bV}{\ensuremath{\mathbf{V}}}

\newcommand{\cA}{\ensuremath{\mathcal{A}}}
\newcommand{\cB}{\ensuremath{\mathcal{B}}}
\newcommand{\cC}{\ensuremath{\mathcal{C}}}
\newcommand{\cD}{\ensuremath{\mathcal{D}}}
\newcommand{\cE}{\ensuremath{\mathcal{E}}}
\newcommand{\cF}{\ensuremath{\mathcal{F}}}
\newcommand{\cG}{\ensuremath{\mathcal{G}}}

\newcommand{\cI}{\ensuremath{\mathcal{I}}}
\newcommand{\cJ}{\ensuremath{\mathcal{J}}}
\newcommand{\cK}{\ensuremath{\mathcal{K}}}

\newcommand{\cM}{\ensuremath{\mathcal{M}}}

\newcommand{\cX}{\ensuremath{\mathcal{X}}}

\title{Analysis of the maximal posterior partition in the Dirichlet Process
Gaussian Mixture Model}
\author{{\L}ukasz Rajkowski,\\
Faculty of Mathematics, Informatics and Mechanics\\
University of Warsaw}

\begin{document}
\maketitle
\begin{abstract}
Mixture models are a natural choice in many applications, but it can be
difficult to place an a priori upper bound on the number of components. To
circumvent this, investigators are turning increasingly to Dirichlet process
mixture models (DPMMs). It is therefore important to develop an
understanding of the strengths and
weaknesses of this approach. This work considers the MAP (maximum a posteriori)
clustering for the Gaussian DPMM (where the cluster means have Gaussian
distribution and, for each cluster, the observations within the cluster have
Gaussian distribution). Some desirable properties of the MAP partition are
proved: `almost disjointness' of the convex hulls of clusters (they may have at
most one point in common) and (with natural assumptions) the comparability of sizes of those clusters that
intersect any fixed ball with the number of observations (as the latter goes to
infinity). Consequently, the number of such clusters remains bounded.
Furthermore, if the data arises from independent identically distributed
sampling from a given distribution with bounded support then the asymptotic MAP
partition of the observation space maximises a function which has a
straightforward expression, which depends only on the within-group covariance
parameter. 
As the operator norm of this covariance parameter decreases, the number of
clusters in the MAP partition becomes
arbitrarily large, which may lead to the overestimation of the number of mixture components.
\end{abstract}



\section{Introduction}\label{sec:intro}
\subsection{Motivation and new contributions}\label{subsec:motiv}
Clustering is a central task in statistical data analysis. A Bayesian approach
is to model data as coming from a~random
mixture of distributions and derive the posterior distribution on the space of
possible divisions into clusters.
When there is not a natural a priori upper bound on the number of clusters, an
increasingly popular technique to use is Dirichlet Process Mixture
Models (DPMMs). It is therefore important to develop an understanding of the strengths and
weaknesses of this approach.

\cite{bib:miller} restrict attention to the number of clusters produced by such
a procedure and are somewhat critical of the method. Their main result implies that in a very
general setting, the Bayesian posterior estimate of the {\em number} of clusters is not
consistent, in the sense that for any $t \in \{1,2,\ldots\}$ almost surely 
\[ \limsup_{n \rightarrow \infty} \mathbb{P}(T_n = t|X_1,\ldots, X_n) < 1,\] 
where $X_1, X_2,\ldots$ is an i.i.d. sample from a mixture with $t$ components
and $T_n$ denotes the number of clusters to which the data are assigned. Here
$\P$ is the probability in the probability space on which $X_1,X_2,\ldots$ are defined.

The Miller and Harrison inconsistency result relates to estimation of the
{\em number} of clusters, not the classification itself. While they do not pursue this, 
they do provide examples of the structure estimators, such as the \emph{MAP}
(maximal a posteriori) partition, which maximises the
posterior probability and the \emph{mean partition}, introduced in
\cite{bib:huelsenbeck2007inference}, which minimises the sum of the
squared distance between the mean partition and all partitions sampled by the
MCMC algorithm which they run, where the {\em distance} is the minimum number of
individuals that have to be deleted from both partitions to make them the same. 

This article presents developments that concern the 
properties of the MAP estimator
in a {\em Gaussian} mixture model, where the cluster {\em centres} are generated according
to a Gaussian distribution and, conditioned on the cluster centre, the
observations within a cluster are generated by Gaussian distribution. The
clusters are generated according to a Dirichlet Process. Analysing the MAP
partition seems to be a natural
choice. It is listed, for example, in \cite{bib:fritsch} as an established
method.
Of course, the set of
possible candidates for the maximiser has to be restricted, since the space
of partitions is too large for an exhaustive search. For example, \cite{bib:dahl2006} suggests choosing the
MAP estimator from a sample from the posterior. He notes, however, a potential
problem of this approach; there may be only a small difference in the posterior
probability between two significantly different partitions. This may indicate
that the classifier is giving the wrong answer as a consequence of mis-specification
of the within-cluster covariance parameter. We investigate such instability in
our examples.

The conclusions of our analysis may be summarised as follows: 
\begin{enumerate}
 \item The convex hulls of the clusters are pairwise `almost disjoint' (they
 may have at most one point in common, which must be a data point).
 \item The clusters are of reasonable size; if $\big(\frac{1}{n}\sum_{j=1}^n
 \|x_j\|^2\big)_{n=1}^\infty$ (the sequence of means of squared Euclidean norms)
 is bounded, then $\liminf_{n \rightarrow \infty} \frac{m_{n}^{[r]}}{n} > 0$ for any
 $r>0$, where
 $m_{n}^{[r]}$ denotes the number of observations in the smallest cluster
 (in the MAP partition of the first $n$ observations) with
 non-empty intersection with $B(\bm{0},r)$ (the ball of radius $r$, centred at
 the origin).
 \item This implies that for any $r>0$ the
 number of clusters in the $n$-th MAP partition required to cover
 observations inside $B(\bm{0},r)$ remains bounded as $n \rightarrow \infty$.
 \item When the data sequence comes from an i.i.d. sampling with bounded support there is an elegant formula to describe the
 limit of the MAP clustering; it is the partition of the observation space that maximises
 the function $\Delta$ given by \cref{eq:DeltaDef}.
In general, though it is a hard problem to find the global maximiser
for this expression. 
Furthermore, the only parameter that this function depends on is the within-group covariance parameter.
\item The {\em negative} finding of the paper is that the clustering is very
sensitive to the specification of the within-cluster variance and model
mis-specification can lead to very misleading clustering. For example, if the
data is i.i.d. from an input distribution which is uniform over a ball of radius
$r$ in $\mathbb{R}^2$ and the within-cluster variance parameter is $\sigma^2 I$,
then for small $\sigma$, the classifier partitions the ball into several,
seemingly arbitrary, convex sets. This classifier therefore has to be treated
with caution.
\end{enumerate}
\subsection{Organisation of the article} \label{subsec:overview}

We now present a brief overview of the structure of the
paper. In \cref{sec:mainres} we give key definitions
and provide complete mathematical
statements of the main results together with intuitive explanations.
\cref{sec:exmp} presents examples which illustrate the results obtained in
the article. These examples show the MAP clustering obtained in various situations where
the data comes from i.i.d. sampling. They indicate that this procedure may
fail to produce reasonable output. The examples are supported by
numerical simulations, which are described in Supplement B.
\cref{sec:basic} contains a detailed presentation of the asymptotic proposition together
with some related developments. In \cref{sec:discuss} we state the open problems
and plans for future work.

\section{Main results}\label{sec:mainres}
\subsection{The Model}\label{subsec:prelim}
This section presents definitions of fundamental notions of our considerations
together with some of their basic properties and relevant formulas. We show how they can be used to construct a~statistical model in which
we expect the data to be generated from different sources of randomness, 
without an a priori upper bound on the number of these sources a~priori.
We start with the definition of the Dirichlet Process, formally introduced in
\cite{bib:ferguson1973bayesian}.

\begin{dfn*}
Let $\Omega$ be a space and $\cF$ a $\sigma$-field of its subsets. Let $\alpha>0$ and
$G_0$ be a probability measure on $(\Omega,\cF)$. The \emph{Dirichlet Process} on
$\Omega$ with parameters $\alpha$ and $G_0$ is a stochastic process
$(G(A))_{A\in\cF}$ such that for every finite partition
$\{A_1,\ldots,A_p\}\subseteq \cF$ of
$\Omega$ the random vector $(G(A_1),\ldots,G(A_p))$ has Dirichlet distribution
with parameters $\alpha G_0(A_1),\ldots, \alpha G_0 (A_p)$. In this case we
write $G\sim \DP(\alpha, G_0)$.
\end{dfn*}

As considered in \cite{bib:antoniak1974mixtures}, the Dirichlet Process can be
used to construct a mixture model in which the number of clusters is not known a
priori. The details are given in the following definition.
\begin{dfn*}
Let $(\Theta,\cF)$ be the parameter space and $(\cX,\cB)$ the observation
space. Let $\alpha>0$ and $G_0$ be a
probability measure on $(\cX,\cF)$. Let $\{F_{\theta}\}_{\theta\in\Theta}$ be a
family of probability distributions on $(\cX,\cB)$. The \emph{Dirichlet Process
mixture model} is defined by the following scheme for generating
a random sample from the space $(\cX,\cF)$
\begin{equation}\label{eq:DPMgen}
\begin{split}
G&\sim \DP(\alpha,G_0)\\
\bm{\theta}=( \theta_1,\ldots,\theta_n )\cond G&\iid G\\
x_i\cond \bm{\theta},G&\sim F_{\theta_i}\qquad \textrm{independently for $i\leq n$}.
\end{split}
\end{equation}
\end{dfn*}

In \cite{bib:blackwell} it is shown that the first two stages of \eqref{eq:DPMgen} may
be replaced by the following recursive procedure: 
\begin{equation}\label{eq:PUSgen}
\theta_1\sim G_0,
\quad
\theta_i\cond \theta_1,\ldots,\theta_{i-1}\sim 
\frac{\alpha}{\alpha+i-1}G_0 + \sum_{j=1}^{i-1} \re{\alpha+i-1} \delta_{\theta_j},
\end{equation}
where $\delta_\theta$ is the probability measure that assigns probability
1 to the singleton $\{\theta\}$. Of course, this procedure can be used to
generate sequences of arbitrary length; the distribution of the resulting infinite sequence
$(\theta_i)_{i=1}^\infty$ produced in this way is called the \textit{Hoppe urn scheme}. Note that a
realisation of this scheme defines a partition of $\N$ by a natural equivalence
relation $(i\sim j)\equiv (\theta_i=\theta_j)$. Restriction of this partition to
sets ${[n]}$ for $n\in\N$ is called the \textit{Chinese Restaurant Process}.

\begin{dfn*}
The \emph{Chinese Restaurant Process} with parameter $\alpha$ is a~sequence of
random partitions $(\cJ_n)_{n\in\N}$, where $\cJ_n$ is a~partition of
$[n]=\{1,2,\ldots,n\}$, that satisfies
\begin{equation} \label{eq:crp}
\cJ_{n+1}\cond \cJ_{n}=\{J_1,\ldots,J_k\}\sim\left\{
\begin{array}{cl}
\{J_1,\ldots,J_{i}\cup\{n+1\},\ldots,J_k\} & \textrm{with probability
$\frac{|J_i|}{n+\alpha}$}\\
\{J_1,\ldots,J_k,\{n+1\}\} & \textrm{with probability $\frac{\alpha}{n+\alpha}$ }
\end{array}\right..
\end{equation}
We write $\cJ_n\sim \CRP(\alpha)_n$.
\end{dfn*}

The Dirichlet Process mixture model 
for $n$ observations is therefore equivalent to
\begin{equation}
\begin{array}{rcll}
\cJ&\sim&\CRP(\alpha)_n&\\
\bm{\theta}=(\theta_J)_{J\in \cJ}\cond \cJ&\iid&G_0&\\
\bx_J=(x_j)_{j\in J}\cond \cJ,\bm{\theta}&\iid&F_\theta&\textrm{for $J\in \cJ$}.
\end{array}
\end{equation}
We will refer to this formulation as the \emph{CRP-based model}.
In this paper we focus our attention on the Gaussian case, in which
$\Theta=\R^d$, $\cX=\R^d$, $\cF$ and $\cB$ are $\sigma$-fields of Borel sets,
$G_0=\Normal(\mu,\Tau)$ and
$F_\theta=\Normal(\theta,\Sigma)$ for
$\theta\in\Theta$, where $\mu\in\R^d$ and $\Tau,\Sigma\in\R^{d,d}$ are the parameters
of the model. This will be called the \textit{CRP-based Gaussian model}.
We also limit ourselves to the case where $\mu=0$, however it may be
easily seen that this is not a~real restriction; the sampling from the
zero-mean Gaussian model and transposing the output by the vector
$\mu$ is equivalent to sampling from the Gaussian model with mean $\mu$.
Therefore all the clustering properties of the model can be investigated with
the assumption that $\mu=0$.

\begin{rem}\label{prop:normalForm}
The
conditional probability of partition $\cJ$ in the zero-mean Gaussian model,
given the observation vector $\mathbf{x}=(x_j)_{j=1}^n$, is proportional to
\begin{equation}\label{eq:posterior}
C^{|\cJ|}\prod_{J\in
\cJ}\frac{|J|!}{|J|^{(d+2)/2}\det R_{|J|}}\cdot
\exp\Big\{\frac{1}{2}\sum_{J\in \cJ}|J|\cdot \znorm{ R_{|J|}^{-1}R^2 \overline{\bx_J}}^2\Big\}
=:\mQ{\cJ}
\end{equation}
where $C=\alpha/\sqrt{\det T}$, 
$ R=\Sigma^{-1/2}$, 
$ R_{m}=(\Sigma^{-1}+T^{-1}/m)^{1/2}$ for $m\in\N$, $\norm{\cdot}$
is the standard Euclidean norm in $\R^d$ and $\overline{\bx_J} = \re{|J|}\sum_{j\in J} x_j$.
\end{rem}
\begin{proof}
See Supplement A.
\end{proof}

Having established the model we are now able to use it for
inference about the data structure. A~natural choice is to 
choose the partition that maximises the posterior probability given by
\eqref{eq:posterior}. This leads to the notion of the MAP partition. 

\begin{dfn*}
The \emph{maximal a posteriori} (MAP) partition of
$[n]$ with observed $\mathbf{x}~=~(x_i)_{i=1}^n$ is any partition of $[n]$ that maximises
$\mQ{\cdot}$ (equivalently, the posterior probability). We denote a~maximiser by
$\hat{\cJ}(\textbf{x})$ (note: a priori this may not be unique).
\end{dfn*}

\subsection{Results}\label{subsec:results}

The first result is \cref{prop:MAPconv} which states that
the MAP partition divides the data into clusters whose convex hulls are
disjoint, with the possible exception of one datum. 
\begin{prop}\label{prop:MAPconv} 
For every $n\in\N$ if $J_1,J_2\in \hat{\cJ}(x_1,\ldots,x_n)$, $J_1\neq J_2$
and $A_k$ is the convex hull of the set $\{x_i\colon i\in J_k\}$ for $k=1,2$ then 
$A_1\cap A_2$ is an empty set or a singleton $\{x_i\}$ for some $i\leq n$.
\end{prop}
\begin{proof}
See Supplement A.
\end{proof}

\begin{figure}[H]
\begin{subfigure}{.32\textwidth}
\centering
\includegraphics[width=\textwidth]{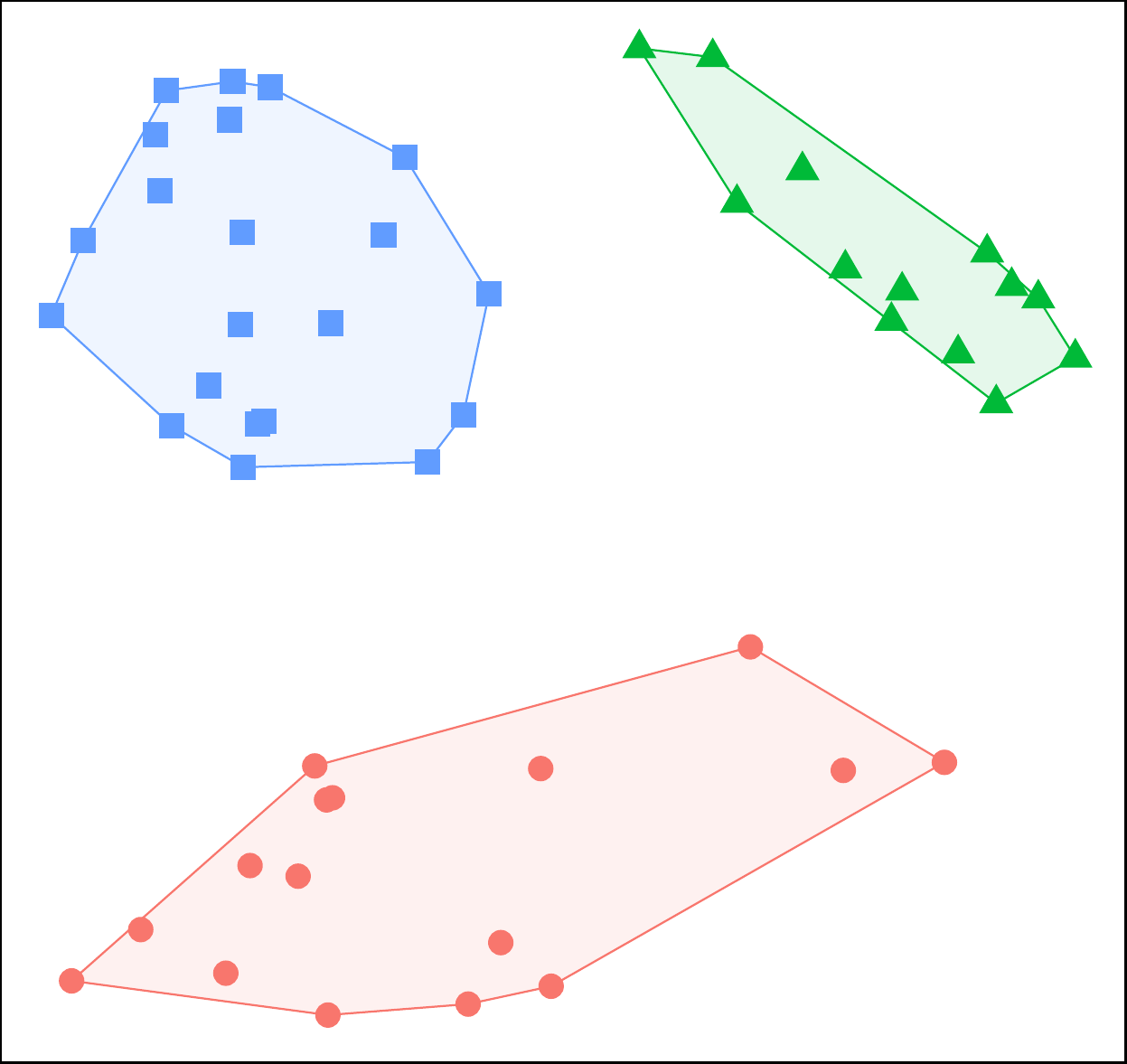}
\caption{This is the desired partition which is also convex.}
\end{subfigure}\hfill
\begin{subfigure}{.32\textwidth}
\centering
\includegraphics[width=\textwidth]{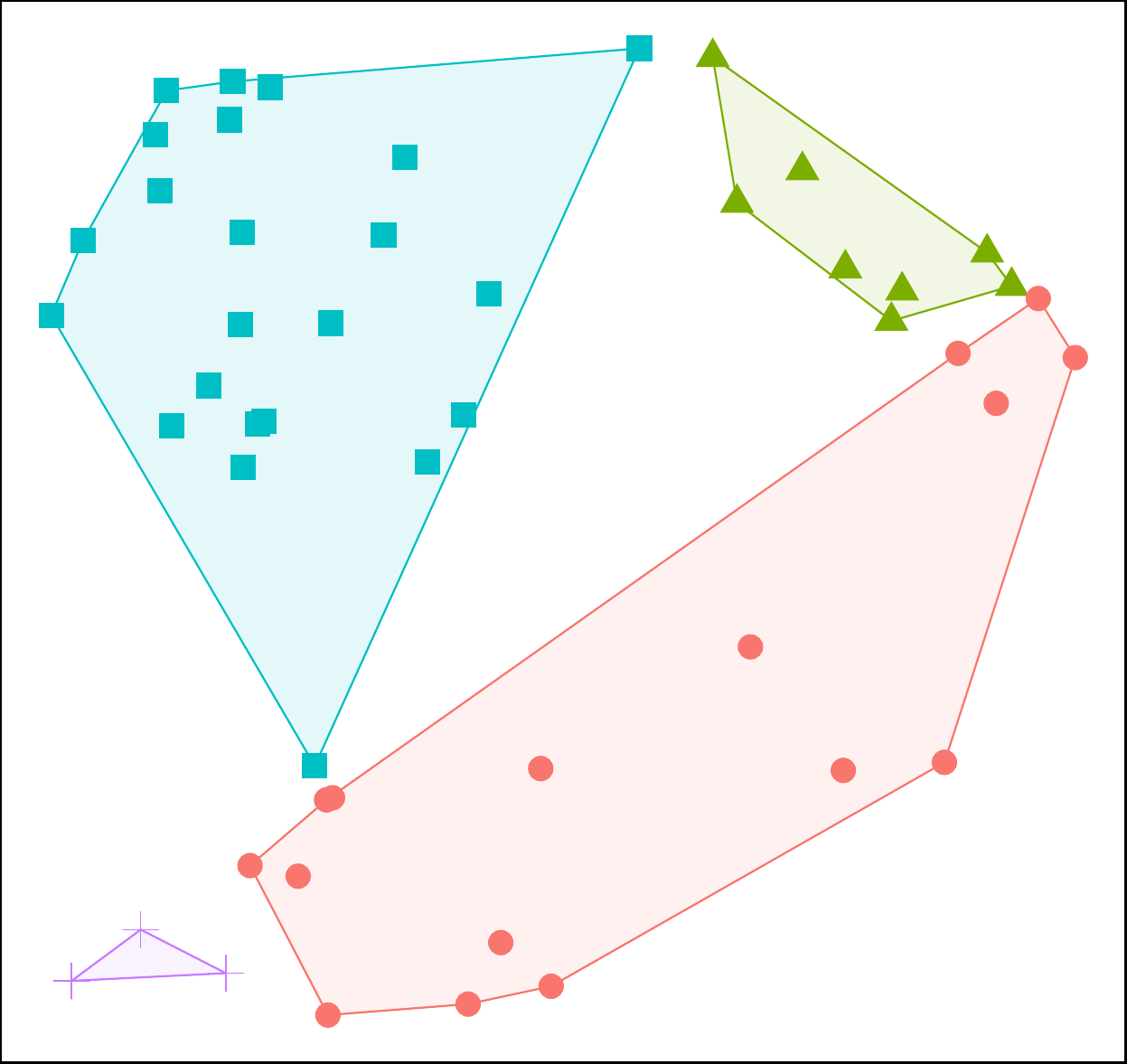}
\caption{This is a convex partition which is not ideal.}
\end{subfigure}\hfill
\begin{subfigure}{.32\textwidth}
\centering
\includegraphics[width=\textwidth]{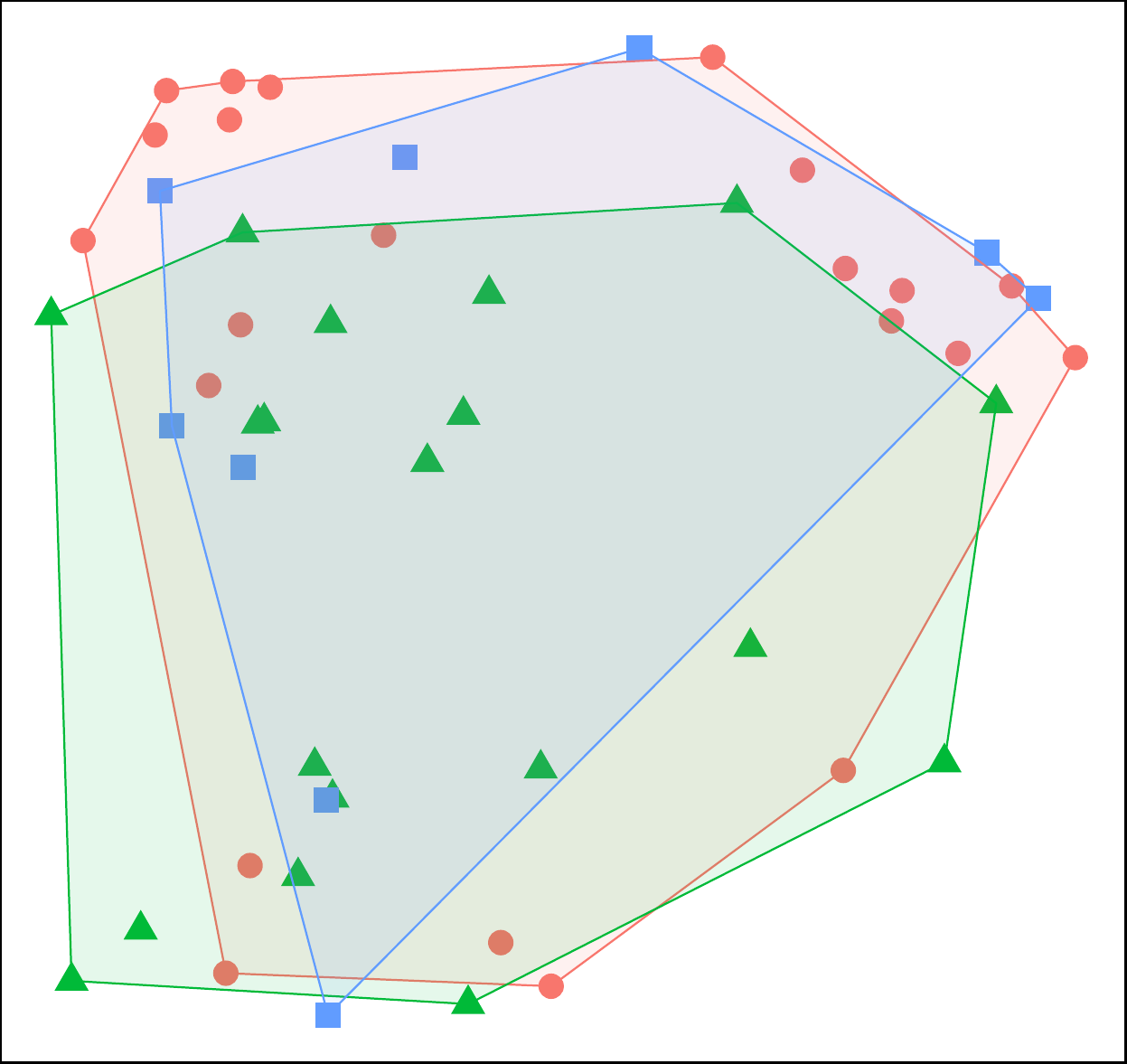}
\caption{This partition is not convex and it is clearly a bad one.}
\end{subfigure}
\caption{
Illustration of the convexity property of a partition of the data.
Clusters are indicated by the shape and colour of the points.
}
\label{fig:convex}
\end{figure}
We say that a partition satisfying the property described by \cref{prop:MAPconv}
is a \emph{convex partition}. As \cref{fig:convex} indicates, this is a rather desirable feature of a
clustering mechanism.

\smallskip
The next development give information about the size and number of the clusters. \cref{prop:xmininter} states that when the sequence
of sample `second moments' is bounded then the size of the smallest
cluster in the MAP partition among those that intersect a ball of given radius is
comparable with the sample size.

\begin{prop}\label{prop:xmininter}
If $\sup_n \re{n}\sum_{i=1}^n \norm{x_n}^2<\infty$ then 
$$\liminf_{n\to\infty}\min\{|J|\colon J\in\hat{\cJ}(x_1,\ldots,x_n), \exists_{j\in
J}\norm{x_j}<r\}/n>0$$ for every $r>0$.
\end{prop}
\begin{proof}
See Supplement A.
\end{proof}

\begin{figure}[H]
\begin{subfigure}{.32\textwidth}
\centering
\includegraphics[width=\textwidth]{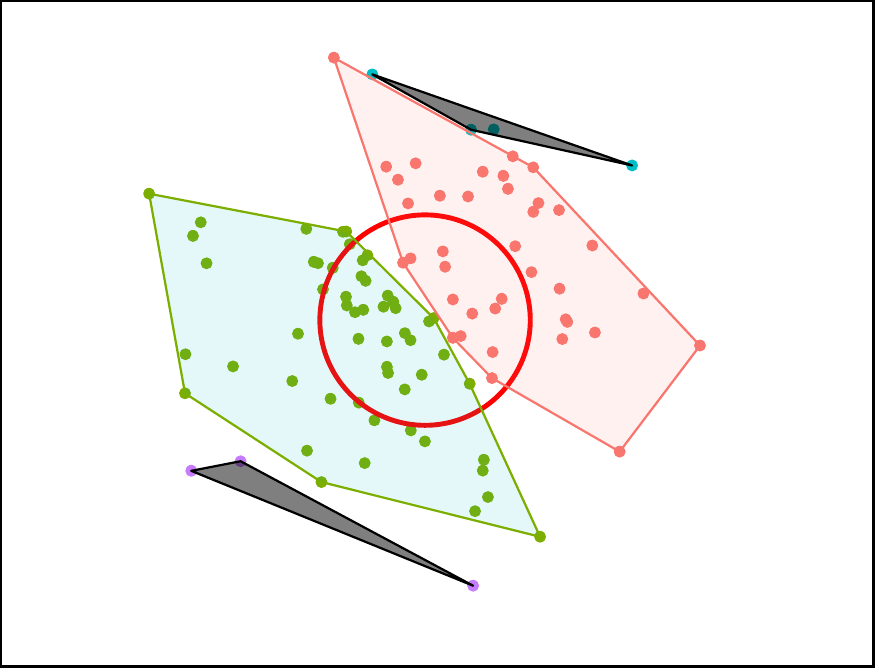}
\caption{$n=100$}
\end{subfigure}\hfill
\begin{subfigure}{.32\textwidth}
\centering
\includegraphics[width=\textwidth]{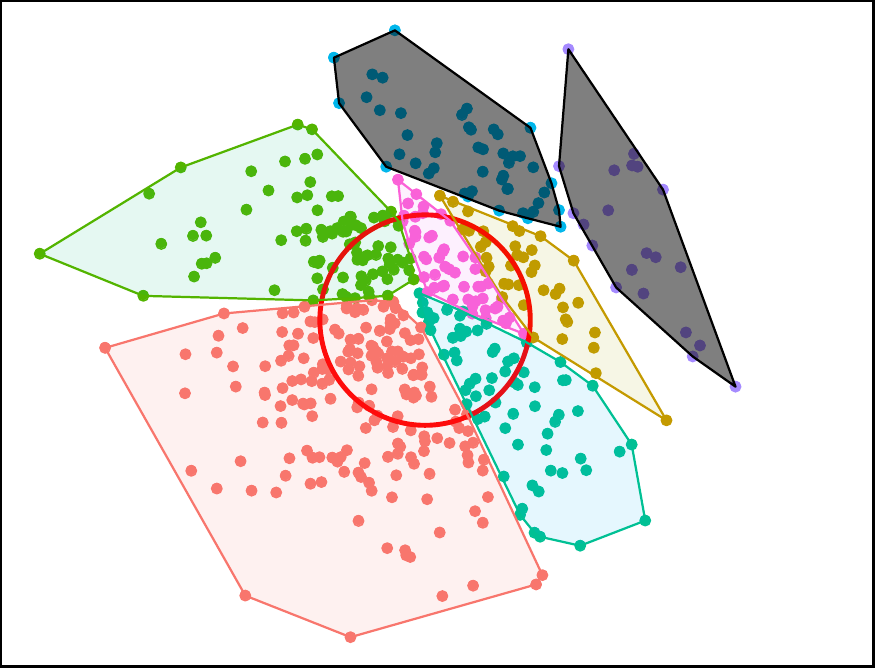}
\caption{$n=500$}
\end{subfigure}\hfill
\begin{subfigure}{.32\textwidth}
\centering
\includegraphics[width=\textwidth]{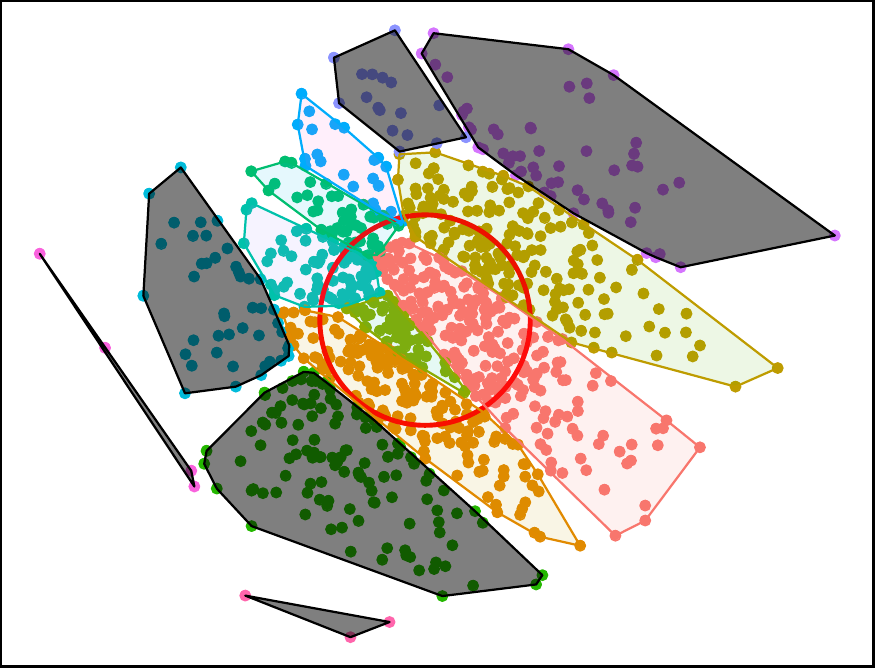}
\caption{$n=1000$}
\end{subfigure}
\caption{
Illustration of \cref{prop:xmininter} and \cref{cor:newbound}. The red
circle is arbitrarily fixed and the clusters it intersects are coloured. The number of
observations in each coloured cluster is proportional to $n$ and the number of
these clusters remains bounded as $n\to\infty$.
}
\label{fig:mapsize}
\end{figure}
The assumption $\sup_n \re{n}\sum_{i=1}^n \norm{x_n}^2<\infty$ allows the data
sequence to be unbounded but it does ensure that it does not grow too quickly. It is
easy to see that an assumption of this kind is necessary, otherwise it would be
possible for each new observation to be large enough to create a new singleton
cluster.

A~simple consequence of \cref{prop:xmininter} is that under these assumptions the number of
components in the MAP partition that intersect a given ball is almost surely bounded. 

\begin{cor}\label{cor:newbound}
If $\big(\re{n}\sum_{i=1}^n \norm{x_i}^2\big)_{n=1}^\infty$ is bounded then
for every $r>0$ the number of clusters that intersect $B(\bm{0},r)$ is bounded, i.e.
$$\limsup_{n\to\infty} |\{J\in\hat{\cJ}(x_1,\ldots,x_n)\colon \exists_{j\in
J}\norm{x_j}<r\}|<\infty.$$
\end{cor}
\begin{proof}
The proof follows easily from the fact that the size
of the smallest cluster that intersects $B(\bm{0},r)$ is bounded from above by the
number of observations divided by the number of clusters intersecting
the ball.
\end{proof}

In order to formulate the central result of the paper we need to introduce
several notions. Let $P$ be a probability distribution on $\R^d$ and $X$ a
random variable with distribution $P$. Let $\Delta$
be the function on the space of finite
families of measurable sets defined by the following formula
\begin{equation}\label{eq:DeltaDef}
\Delta(\cG)=\re{2}\sum_{G\in\cG}P(G)\znorm{R\E (X\cond X\in G)}^2+
\sum_{G\in\cG}P(G)\ln P(G).
\end{equation}
where $R^2$ is the inverse of the within-cluster covariance matrix $\Sigma$ and 
$\E(X\cond X\in G)$ is the expected value of $X$ conditioned on $X\in G$.

We consider the \emph{symmetric distance metric} over $P$-measurable
sets, which is defined by
$d_P(A,B)=P\big( (A\sm B)\cup(B\sm A)\big)$. This can be easily extended to a metric $\overline{d_P}$
over finite families of measurable subsets of $\R^d$ (details are given in
\cref{sec:convMAP}). Also we say that a family of
measurable sets $\cA$ is a \emph{$P$-partition} if $P(\bigcup_{A\in\cA}A)=1$ and
$P(A\cap B)=0$ for all $A,B\in \cA$, $A\neq B$. 
Let $\bm{M}_\Delta$ denote the set of finite $P$-partitions that
maximise the function $\Delta$.

Consider $X_1,X_2,\ldots\iid P$ and let $\hat{\cA}_n$ be the family of the convex hulls of clusters of observations in
$\hat{\cJ}(X_1,\ldots,X_n)$.

\begin{prop}\label{thm:distto0} Assume that $P$ has bounded support and is
continuous with respect to Lebesgue measure. Then $\bm{M}_\Delta\neq\emptyset$
and almost surely $\inf_{\cM\in\bm{M}_\Delta}
\overline{d_P}(\hat{\cA}_n,\cM)\to 0$.
\end{prop}
\begin{proof}
The proof follows from \cref{thm:MAPmax}. See Supplement A for details.
\end{proof}

The function $\Delta$ does not depend on the concentration parameter $\alpha$
or the between-groups covariance parameter. It therefore follows, somewhat surprisingly, that in
the limit the shape of the MAP partition does not depend on these two parameters.

It can be shown that as the norm of the within group covariance matrix tends to 0, the variance of the
conditional expected value gains larger importance in maximising the function $\Delta$ in formula \eqref{eq:DeltaDef} and this variance increases as the
number of clusters increases.
Therefore by manipulating the within group covariance parameter, when the input
distribution is bounded it is possible
to obtain an arbitrarily large (but fixed) number of clusters in the MAP
partition as
$n\to \infty$, as \cref{prop:KMAPlarge} states.
This is also an indication of the inconsistency of the 
procedure used since it implies that when the input comes from a finite mixture of
distributions with bounded support, then setting the $\Sigma$ parameter too small
leads to an overestimation of the number of clusters.

\begin{prop}\label{prop:KMAPlarge}
Assume that $P$ has bounded support and is
continuous with respect to Lebesgue measure. Then for every $K\in\N$ there
exists an $\eps>0$ such that if $\norm{\Sigma}<\eps$ then $|\hat{\cJ}_n|>K$ for
sufficiently large $n$.
\end{prop}
\begin{proof}
See Supplement A.
\end{proof}

It is worth pointing out that \cref{prop:MAPconv} and \cref{prop:xmininter} hold also for
\emph{finite} Gaussian mixture models with Dirichlet prior on the probabilities of
belonging to a given cluster. \cref{thm:distto0} also remains true with
$\bm{M}_\Delta$ replaced by $\bm{M}_\Delta^K$ -- the set of $P$-partitions with
at most $K$ clusters that maximise the function $\Delta$, where $K$ is the
number of clusters assumed by the model. The details are left for Supplement A.

\section{Examples}\label{sec:exmp}

This section presents some examples which illustrate the main propositions of the article.
In \cref{subsec:segment} we compute the convex partition that maximises
$\Delta$ when $P$ is a~uniform distribution on the interval $[-1,1]$. 
\cref{subsec:exp} gives an example of a distribution with well-defined moments,
for which the maximiser of $\Delta$ necessarily has infinitely many
clusters, although 
for any $r < \infty$, the number of clusters that intersect a ball of radius $r$
is finite. This
example illustrates the content of \cref{thm:manyclust}, where it is shown that
with appropriate choice of model parameters, if the input distribution is
exponential then the number of clusters in the sequence of MAP partitions
becomes arbitrarily large.
 \cref{subsec:mnorm}
 investigates Gaussian mixture models; the MAP partition does not properly identify the two clusters
 when the mixture distribution is bi-modal.
Finally, in \cref{subsec:disc} we consider the uniform distribution on the unit disc in
$\R^2$. The partition maximising the function $\Delta$ cannot be obtained by
analytical methods, but it may be approximated. The results approximate the
optimal partition of the unit disc and illustrate the convexity of
\cref{prop:MAPconv}.
All examples are substantiated with computer simulations, presented in the main
text or in Supplement B.

\subsection{Uniform distribution on an interval}\label{subsec:segment}
We find the convex partition that maximises $\Delta$ if $P$ is a~uniform
distribution on $[-1,1]$. Firstly we find an optimal partition with fixed number
of clusters $K$. Since it is
convex, it is defined by the lengths of $K$ consecutive
subintervals of $[-1,1]$. Let those be $2p_1,\ldots,2p_n$. Computations in Supplement A show that
with $K$ fixed the optimal division is $p_1=p_2=\ldots=p_K=1/K$. Using this, it
is computed that
the optimal number of clusters is $K=\lfloor
R/\sqrt{3}\rfloor$ or $K=\lceil R/\sqrt{3}\rceil$, where $\lfloor x
\rfloor$ and $\lceil x \rceil$ are the largest integer not greater than
$x$ and the smallest integer not less than $x$, respectively. 
It is worth noting that the variance of the data within a segment of length
$2R/\sqrt{3}$ is equal to $R$, so in this case the MAP clustering splits the data in a way
that adjusts the empirical within-group covariance to the model assumptions.

It should be underlined that in this example, if $\Sigma$ is small, the MAP
partition has more than one cluster. 
The clustering is therefore misleading, since in this case there is exactly one population (which is uniform $[-1,1]$). The number of clusters in the MAP
partition becomes arbitrarily large as $\Sigma$ goes to $0$, as \cref{prop:KMAPlarge}
states.

This would suggest that, in general, a sensible choice of $\Sigma$ should be
made a priori. The sample variance would give an upper bound on $\Sigma$
(since the data variance is the sum of between-group and within-group
variances), but there is no natural lower bound for this parameter.
In this example the partitioning mechanism itself is clearly far from
satisfactory when it produces more than two clusters; the divisions seem very
arbitrary.

\subsection{Exponential distribution}\label{subsec:exp}
When the input distribution is exponential with parameter 1, then for a relevant
choice of model parameters (e.g. $\alpha=\Tau=1$, $\Sigma=4$) there is
no finite partition that maximises $\Delta$;
the value of the function $\Delta$ for a given convex partition can be increased by
taking any interval of length larger than 3 and dividing it into two equally
probable parts. See Supplement A for the proof.

Since the exponential distribution does not have bounded support, our
considerations regarding the relation between the function $\Delta$ and the MAP
clustering cannot be applied directly. However, by using similar methods we can establish
that for exponential input the MAP procedure creates an arbitrarily
large number of clusters. This is stated in \cref{thm:manyclust}, whose proof is
presented in Supplement A.

\subsection{Mixture of two normals}\label{subsec:mnorm}
Let the input distribution be a mixture of two normals
($P=\re{2}(\nu_{-1.01}+\nu_{1.01})$), where $\nu_m$ is the normal distribution
with mean $m$ and variance 1). It can be proved that this distribution is
bi-modal (however slightly; see Supplement A). Choose the model parameters
consistent with the input distribution, i.e.
$d=\alpha=\Sigma=\Tau=1.$
It can be computed numerically that
$\Delta(\{(-\infty,0],(0,\infty)\})\approx -0.0046<0=\Delta(\{\R\})$. 
An intuitive partition of the data into positive and negative is induced by the
partition $\{(-\infty,0],(0,\infty)\}$ and hence, by
\cref{cor:betterness}, for sufficiently large data input
the posterior score for the two clusters partition is smaller than the posterior
score for a single cluster. This may be taken as an indication of inconsistency of
the MAP estimator in this setting.

\subsection{Uniform distribution on a disc}\label{subsec:disc}
This gives an example of non-uniqueness of the optimal partition, since the family of optimal partitions is clearly invariant under rotation
around $(0,0)$. Let $P$ be uniform distribution on $B(\bm{0},1)$. It can be easily seen that $\Delta(B(\bm{0},1))=0$. Let
$R$ be the identity matrix and let $B_1^+$ ($B_1^-$) be a~subset of $B(\bm{0},1)$ with non-negative
(negative) first coordinate. Then
$
\Delta(\{B_r^+,B_r^-\})=2r^2/9-\ln 2.
$
Therefore, for sufficiently large
$r$, a partition of $B(\bm{0},1)$ into halves is better than a single cluster, hence the
optimal convex partition $\cE$ is not a~single cluster. Since a single cluster
is the only convex partition of $B(\bm{0},1)$ that is rotationally invariant
about the origin, it follows that the optimal partition is not unique.

The simulation in this case also give a nice illustration of the
convexity of the MAP partition, proved in \cref{prop:MAPconv} and show the
arbitrary nature of the partitioning when $r$ is large. 

\begin{figure}[H]
\includegraphics[width=\textwidth]{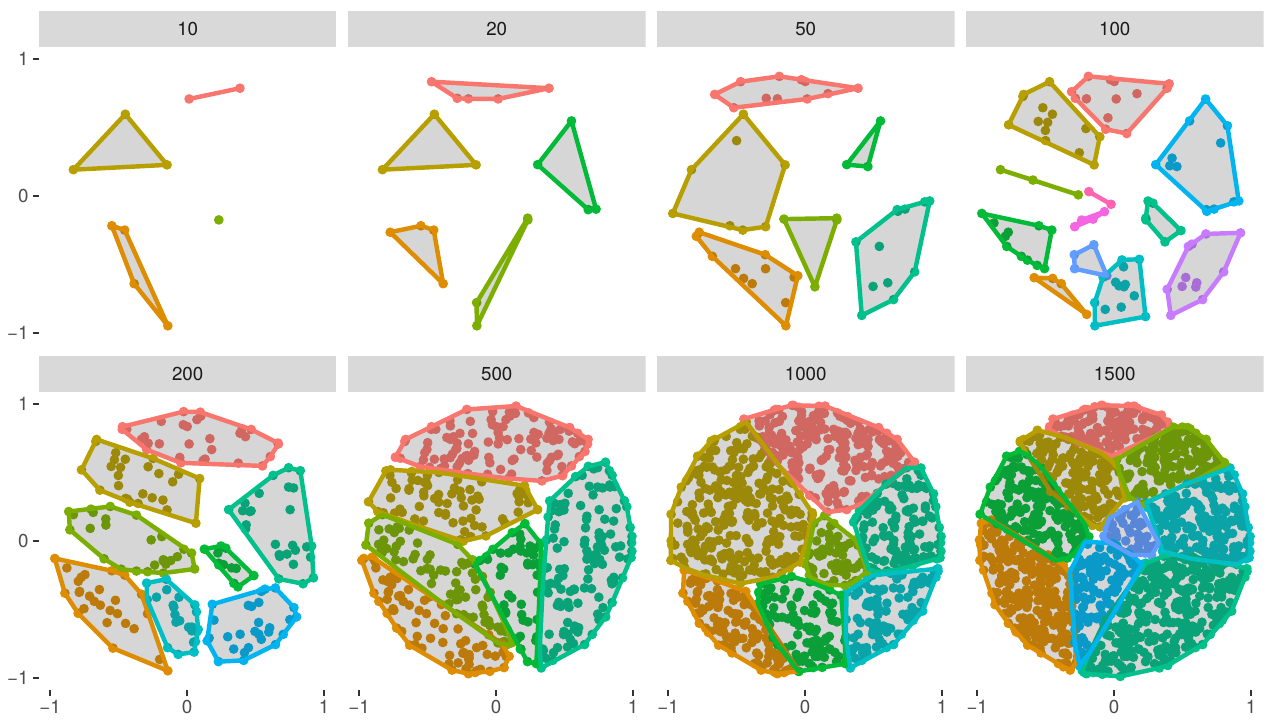}
\caption{
Clustering in the MAP partition of the first $k=100,500,1000,1500,2000$
observations (in columns) in the i.i.d. sample from the uniform distribution of 
disc $B(\bm{0},1)$. Different clusters are denoted by different colours.
}
\label{fig:circle}
\end{figure}

\nopagebreak
\subsection{The MAP clustering properties}
\begin{figure}
\centering
\includegraphics[width=.99\textwidth]{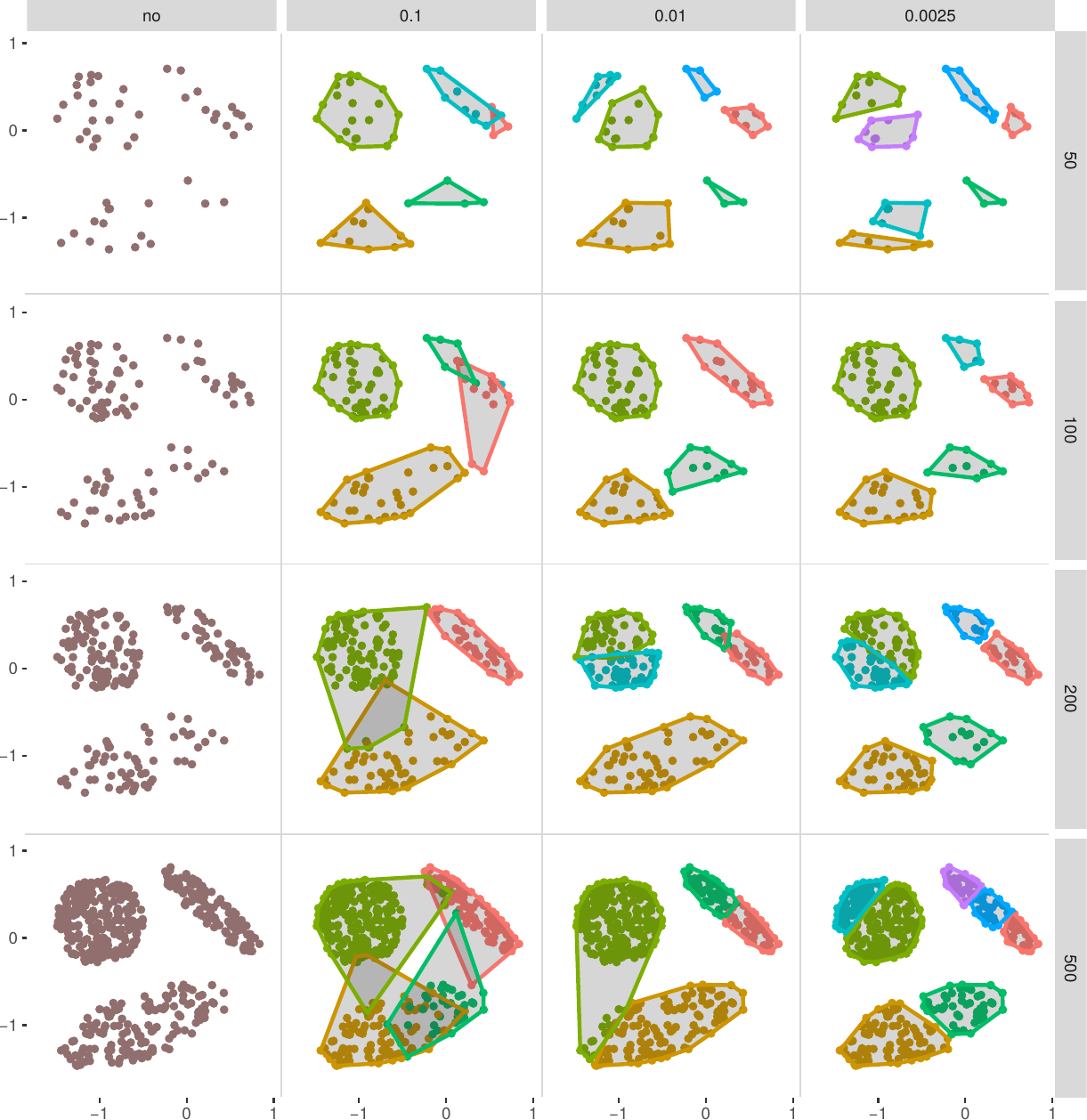}
\caption{
Clustering in the MAP partition of the first $k=50, 100, 200,
500$ observations (in columns) in the i.i.d. sample from the mixture of three uniform
distributions on a disjoint ellipses. The MAP clustering was constructed for $\alpha=1$,
$\Tau=\mI$ and $\Sigma=\sigma^2\mI$ where $\sigma^2\in \{1,.1,.01,.0025\}$ (in rows).
Different clusters are denoted by different colours, the convex hulls of the clusters
are also marked.
It is clear that some of the partitions presented are not convex, particularly
for large $\sigma^2$. 
This is due to the fact that the method is less than perfect. As $\sigma^2$
increases, the likelihood component of the formula for the posterior is less
significant and hence partitions with the same prior (where clusters are of the
same size) have similar posterior score. Therefore, with high probability,
sampling from the posterior will not choose the MAP partition, or even a partition that
reasonably resembles the MAP clustering.  
We mentioned this instability in \cref{subsec:motiv}.
}
\label{fig:sevclust}
\end{figure}

This short simulation study presents the performance of the MAP estimator when
the input distribution is a mixture of uniform distributions on three pairwise
disjoint ellipses. The output is shown on \cref{fig:sevclust}. It shows that the
MAP clustering detects the mixture components or at least the clusters it creates are the
sub-groups of the true mixture components (all depending on the within-group
covariance parameter $\Sigma$). It also provides a
nice illustration for two properties of the MAP partition: firstly the convexity property
(\cref{prop:MAPconv}) and secondly -- the fact that when the within-group
covariance parameter is decreasing, the number of cluster in the MAP partition
grows, as stated in \cref{prop:KMAPlarge}.
\section{Detailed presentation of \cref*{thm:distto0}}\label{sec:basic}

\subsection{Classification of Randomly Generated Data}\label{sec:random}

Let $P$ be a~probability distribution on $(\R^d,\cB)$ and let
$(X_n)_{n=1}^\infty$ be a~sequence of independent copies of a~random variable
$X$ with distribution $P$. Then $\hat{\cJ}_n=\hat{\cJ}(X_1,\ldots,X_n)$ goes
a random partition of $[n]$. 
Note that if $\E\norm{X}^4<\infty$ (here and subsequently, $\E$ denotes the
expected value) then by the strong law of large numbers almost surely $\re{n}\sum_{i=1}^n
\norm{X_i}^2\to \E\norm{X}^2<\infty$ and therefore the assumptions of
\cref{prop:xmininter} are satisfied almost surely. Useful corollaries of this observation are listed below.

\begin{cor}\label{cor:min}
If $\E\norm{X}^4<\infty$ then for every $r>0$ 
\vspace*{-.3cm}
\begin{enumerate}[(a)]
 \setlength{\itemsep}{1pt}
 \setlength{\parskip}{0pt}
 \setlength{\parsep}{0pt}
\item $\liminf_{n\to\infty}\min\{|J|\colon J\in\hat{\cJ}_n, \exists_{j\in
J}\norm{X_j}<r\}/n>0$ almost surely.
\item the number of clusters in $\hat{\cJ}_n$ that intersect $B(\bm{0},r)$ is bounded.
\end{enumerate}
\end{cor}

An easy consequence of \cref{cor:min} is
\begin{cor}\label{cor:max}
If the support of $P$ is bounded then 
\vspace*{-.3cm}
\begin{enumerate}[(a)]
 \setlength{\itemsep}{1pt}
 \setlength{\parskip}{0pt}
 \setlength{\parsep}{0pt}
	\item $\liminf_{n\to\infty}\min\{|J|\colon J\in\hat{\cJ}_n\}/n>0$ almost surely.
	\item $|\hat{\cJ}_n|$ is almost surely bounded.
\end{enumerate}
\end{cor}
\begin{proof}
If the support of $P$ is bounded then $\E \norm{X^4}<\infty$. Therefore we can
use \cref{cor:min} where we take $r$ sufficiently large so that $B(\bm{0},r)$ contains
the support of $P$.
\end{proof}

The assumptions of \cref{cor:max} cannot be relaxed to those of \cref{cor:min}. 
It turns out that there exists a~probability distribution $P$ with a countable number of
atoms sufficiently far apart, whose probabilities are chosen so that $\E
\norm{X}^4<\infty$ and almost surely the most recent observation
creates a~singleton in the sequence of MAP partitions infinitely often, i.e. there exists a sequence
$(n_k)_{k=1}^\infty$ such that $\{x_{n_k}\}\in \hat{\cJ}_{n_k}$. This violates part (a) of \cref{cor:max}. On the
other hand, for appropriate parameter choice, sampling from the exponential
distribution leads to the number of clusters in the MAP partition tending to infinity,
which contradicts part (b) of \cref{cor:max}. Proofs of these facts are left for
Supplement A. These facts are now formally stated in the following two theorems:

\begin{thm}\label{thm:nosmallm}
If $d=1$ and $\alpha=\Tau=\Sigma=1$ then for $P=\sum_{m=0}^\infty
q(1-q)^{m}\delta_{18^m}$, where $q=(2\cdot 18)^{-1}$, 
almost surely $\liminf_{n\to\infty}m(\hat{\cJ}_n)=1$. 
\end{thm}

\begin{thm}\label{thm:manyclust}
If $P=\Exp(1)$ and the CRP model parameters are
$\alpha=\Tau=1$, $\Sigma<( 32\ln 2 )^{-1}$ then 
the number of clusters in the sequence of MAP partitions almost surely goes to infinity.
\end{thm}

\subsection{The Induced Partition}\label{sec:induced}

Instead of searching for the MAP clustering, one may choose a simpler (and more
arbitrary) way to partition the data. The idea is to choose a
partition of the observation space in advance and then divide the sample
assigning each datum to the element of this partition which contains it. We call this decision rule an \textit{induced partition}.
In this section we give a~formal definition and investigate how it behaves when
the input is identically distributed and how it relates to the formula for the
posterior probability given by~\eqref{eq:posterior}.

\begin{dfn*}
Let $\cA$ be a~fixed partition of $\R^d$. For $n\in\N$ and $A\in \cA$ let $J^A_n=\{i\leq
n\colon X_i\in A~\}$ and define a~random partition of $[n]$ by 
$
\cJ^{\cA}_n=\big\{J^A_n\neq\emptyset\colon A\in\cA\big\}.
$
We say that this partition of $[n]$ is \emph{induced by $\cA$}. 
\end{dfn*}

In the following part of the text, for two sequences $(a_n)_{n=1}^\infty$ and
$(b_n)_{n=1}^\infty$ of nonzero real numbers, we use the notation
$
a_n\approx b_n$ to denote $\lim_{n\to\infty} a_n/b_n=1.
$

\begin{lem}\label{lem:induced}
Let $\cA$ be a~finite $P$-partition of $\R^d$ consisting of Borel sets with
positive $P$ measure. Then almost surely
$
\sqrt[n]{Q_{\bX_{1:n}}(\cJ^\cA)}\approx 
\frac{n}{e}\exp\left\{\Delta(\cA)\right\},
$
where $\Delta$ is the function defined by \eqref{eq:DeltaDef}.\end{lem}
\begin{proof}
See Supplement A.
\end{proof}

\begin{cor}\label{cor:betterness}
If $\cA$, $\cB$ are two finite $P$-partitions of $\R$ such that
$\Delta(\cA)>\Delta(\cB)$
then almost surely 
$
Q_{\bX_{1:n}}(\cJ_n^{\cA})>Q_{\bX_{1:n}}(\cJ_n^{\cB})\ \
\textrm{for sufficiently large $n$.}
$
\end{cor}
\begin{proof}
The proof is straightforward and therefore omitted.

\end{proof}

\cref{cor:betterness} implies that if we look for the optimal, finite induced
partition, it will be a partition of the data induced by the finite partition of the
observation space that maximises the function $\Delta$. This formulation suggests a
strong relationship between the MAP partition and the finite maximisers of $\Delta$, which will
be investigated further in \cref{sec:convMAP}, in the case where $P$ has
bounded support. 
The case where $P$ does not have bounded support is beyond the scope of this work, for reasons presented
in \cref{sec:discuss}. This is a goal for future research.

At the end of this Section, let us provide an interpretation of the function
$\Delta$. Let $\cA$ be a finite partition and $Z_\cA=\E\left(X|\1_A(X)\colon
A\in\cA\right)$ be the conditional expected value of
$X$ given the indicators $\1_A(X)$ for $A\in\cA$.
Then $Z_\cA$ is a discrete random variable which is equal to $\E(X\cond X\in A)$
with probability $P(A)$. This implies that 
$
\Delta(\cA) = \re{2}\E \norm{RZ_\cA}^2 - H(Z_\cA),
$
where the function $H$ assigns to a random variable its entropy. Moreover
$$\E \norm{RZ_\cA}^2= \tr\big(\bV ( RZ_\cA )\big) + \norm{\E RZ_\cA}^2
= \tr\big(R\bV ( Z_\cA )R^t\big) + \norm{R\E Z_\cA}^2$$
in which $\tr(\cdot)$ is the trace function and $\bV(\cdot)$ is the covariance
matrix of a given random vector.
Since $\E Z_A=\E X$ we obtain that
\begin{equation}\label{eq:DeltaForm}
\Delta(\cA) = \re{2}\tr\big(R\bV ( Z_\cA )R^t\big) - H(Z_\cA)+\re{2}\norm{R\E
X}^2.
\end{equation}
Equation \eqref{eq:DeltaForm} justifies the following description of the
function $\Delta$: up to a constant, it may be treated as a difference between the
variance and the entropy of the conditional expected value of a linearly
transformed, $P$-distributed
random variable given its affiliation to one of the sets in the partition.

\subsection{Convergence of the MAP partitions}\label{sec:convMAP}
\cref{cor:betterness} gives us a~convenient characterisation of the partitions
of $\R^d$ that in the limit induce the best possible partitions of sets
$[n]$. At this stage however we do not know yet if the best induced partitions
relate to overall best partitions, namely the MAP partitions. 
A natural question is if the
behaviour of the MAP partition resembles the induced classification introduced in
\cref{sec:induced}, as the sample size goes to infinity, and under what
conditions. 
This section presents partial answers in this regard; it should
be stressed however that all the developments presented here are limited to the case
when the input distribution has bounded support. The reasons for such
limitation
are briefly described in \cref{sec:discuss}. 

As we already know that clusters
in the MAP partition create disjoint convex sets, the analysis of the approximate
behaviour of these partitions would be easier if a~form of `uniform law of large
numbers' with respect to the family of convex sets were true. More precisely
if we let $P_n=\re{n}\sum_{i=1}^n \delta_{X_i}$ we need the following to hold:
\begin{equation}\label{eq:gcc}
\lim_{n\to\infty}\sup_{C\textrm{ convex}} \big|P_n(C)-P(C)\big|=0
\quad\textrm{almost surely.}\tag{$\ast$}
\end{equation}
In other words we require that the class of convex sets is a~\textit{Glivenko-Cantelli
class} with respect to $P$. A~convenient condition for this to hold is
given in \cite{bib:cantelli}, Example 14:

\begin{lem}\label{lem:gliwcan}
If for each convex set $C$ the boundary $\partial C$ can be covered by countably many
hyperplanes plus a~set of $P$-measure zero, then~\eqref{eq:gcc} holds for
$P$.
\end{lem}

In particular, it can easily be seen that the assumptions of \cref{lem:gliwcan}
are satisfied if $P$ has a~density with respect to Lebesgue measure $\lambda_d$ on $\R^d$
(since in this case the Lebesgue measure $\lambda_d$ of the boundary of any
convex set is 0, and hence is also $P$ measure 0).

We can now formulate a functional relation between the posterior probability of
the MAP partition and the value of the function $\Delta$ on the family of convex
hulls of the sets in the MAP partition.
\begin{lem}\label{lem:approxAn}
Assume that $P$ has bounded support and satisfies~\eqref{eq:gcc}. Let
$\hat{\cA}_n$ be the family of the convex hulls of the clusters in the MAP
partition,
i.e.
$\hat{\cA}_n=\big\{\conv\{\bX_j\colon j\in J\}\colon J\in\hat{\cJ}\big\}$. 
Then almost surely
$$\sqrt[n]{Q_{\bX_{1:n}}(\hat{\cJ}_n)}\approx
\frac{n}{e}\exp\{\Delta(\hat{\cA}_n)\}.$$
\end{lem}
\begin{proof}
See Supplement A.
\end{proof}

Now we investigate the convergence of the sequence $\hat{\cA}_n$ defined in
\cref{lem:approxAn}. In order to do so we need a~topology on relevant subspaces of $2^{\R}$.
We begin by recalling two standard metrics used in this context. 

\begin{dfn*}
Let $\cD$ be a~class of closed subsets of $\R^d$. Then the function
$\varrho_H\colon \cD^2\to \R$ defined by
$$
\varrho_H(A,B)=\inf\{\eps>0\colon A\subseteq (B)_\eps, B\subseteq (A)_\eps\},
$$
where $(X)_\eps=\{x\in\R^d\colon \dist(x,X)<\eps\}$, is a~metric on $\cD$. It is called
the \textit{Hausdorff distance}. The fact that it is a metric follows from 1.2.1
in \cite{bib:moszynska}.
\end{dfn*}

\begin{dfn*}
Let $\cM$ be a~$\sigma$-field on $\R^d$ and $\mu$ be a~measure on $(\R^d,\cM)$. Then the function
$d_\mu\colon \cM^2\to\R$ defined by
$
d_{\mu}(A,B)=\mu\big((A\sm B)\cup (B\sm A)\big)
$
is a~pseudometric on $\cM$,
which by definition means that it is symmetric,
nonnegative and satisfies the triangle inequality. It is called the \textit{symmetric difference metric}. The fact that it is a pseudometric is explained in
the beginning of Section 13, Chapter III of \cite{bib:doob}.
Note that since $d_\mu(A,B)=0$ does not imply $A=B$, formally $d_\mu$ is not a~metric
on $\cM$. Although for our consideration the difference of measure 0 is of no
importance, we keep on using the proper \textit{pseudometric} term in this context.
\end{dfn*}

The two following theorems are crucial for establishing
the limits of maximisers. \cref{thm:dHcomp} is Theorem 3.2.14 in
\cite{bib:moszynska}; it ensures the existence of $d_H$-converging
subsequence in every bounded sequence of convex sets.
\cref{thm:dPcomp} is a straightforward consequence of Theorem 12.7 in \cite{bib:Valentine64}
(in the latter $P$ is taken to be the Lebesgue measure).
It states that when $P$ has a density with respect to the Lebesgue measure then the Hausdorff metric restricted to 
$\cK$ is stronger than the symmetric difference metric.

\begin{thm}\label{thm:dHcomp}
The space $(\cK,\varrho_H)$ is finitely compact (i.e. every
bounded sequence has a~convergent subsequence).
\end{thm}
\begin{thm}\label{thm:dPcomp}
If $P$ is continuous with respect to the Lebesgue measure then convergence in
$\varrho_H$ implies convergence in $d_P$ in the space $\cK$.
\end{thm}

Note that the Hausdorff and symmetric difference metrics are defined on
sets. However we are interested in MAP partitions, which are \emph{families} of
sets. Therefore it is convenient to extend the definitions of these metrics to
families of sets, as presented below. \cref{thm:famconv} ensures that the
desirable properties of compactness are preserved by such extension.

\begin{dfn*}\label{dfn*:metrExt}
Let $d$ be a~pseudometric on the family of sets $\cF$. For $K\in\N$ we define
$F_K(\cF)$ to be the space of finite subfamilies of $\cF$ that have at most $K$
elements. Moreover
$\cA=\{A^{(1)},\ldots,A^{(k)}\}\in F_K(\cF)$ and
$\cB=\{B^{(1)},\ldots,B^{(l)}\}\in F_K(\cF)$ we define
\begin{equation}
\bar{d}(\cA,\cB)=\min_{\sigma\in\Sigma_K} \max_{i\leq K}
d(A^{(i)},B^{(\sigma(i))}),
\end{equation}
where $\Sigma_K$ is the set of all permutations of $[K]$ and we assume
$A^{(i)}=\emptyset$ and $B^{(j)}=\emptyset$ for $i>k$ or $j>l$ respectively.
\end{dfn*}
 
\begin{rem}\label{thm:famconv}
If $(\cF,d)$ is a~pseudometric space then $(F_K(\cF),\bar{d})$ is also a
pseudometric space. 
Moreover, if $(\cF,d)$ is finitely compact then $(F_K(\cF),\bar{d})$ is also finitely
compact.
\end{rem}
\begin{proof}
The proof is straightforward. See Supplement A for details.
\end{proof}

Now assume that $P$ has bounded support. Then by
\cref{thm:dHcomp} and~\cref{thm:famconv} it follows that
$(\hat{\cA}_n)_{n=1}^\infty$ has convergent subsequences
which have a limit under $\overline{d_H}$ (note that as the support of $P$ is
bounded, sets
$\hat{\cA}$ are also bounded in the $d_H$ metric). Let us denote the (random) set of
their limits by $\bm{E}$. Note that
by the properties of $d_H$ distance each family in $\bm{E}$ consists of convex,
closed sets. If we assume that $P$ is continuous with respect to the Lebesgue
measure 
then it follows from \cref{lem:approxAn} together with \cref{thm:dPcomp} that
$\bm{E}$ consists of
finite $P$-partitions that maximise the function $\Delta$.

\begin{thm}\label{thm:MAPmax}
Assume that $P$ has bounded support and is continuous with respect to Lebesgue
measure. Then every partition in $\mathbf{E}$ is a~finite $P$-partition that maximises $\Delta$.
\end{thm}
\begin{proof}
See Supplement A.
\end{proof}

Now \cref{thm:distto0} is a straightforward, topological consequence of
\cref{thm:MAPmax}. This is shown in Supplement A.

\section{Discussion}\label{sec:discuss}
It should be clearly stated that the scope of the paper is limited in two
ways. Firstly, only the Gaussian model is considered. It is natural
to ask if the methods used here can be applied for other combinations of
base measure and component distributions. The author is sceptical in this
regard. The proofs of the key \cref{prop:MAPconv} and
\cref{prop:xmininter} rely strongly on the formula \eqref{eq:posterior}. It is
difficult to find a computationally feasible choice of the base and component
measures so that the resulting formula for the posterior probabilities has
similar properties.

Secondly, the limiting results contained in \cref{sec:convMAP} are proved in the
case where the support of the input distribution is bounded. In this case the
model is clearly misspecified. A significant effort was put in order to extend
the results from \cref{sec:convMAP} at least to the case where $P$ is Gaussian.
Unfortunately, there are some technical hurdles which the author was not able
to overcome, which we now outline. The first result in which the boundedness of
the input distribution is used is \cref{lem:approxAn} -- here we use both parts
of \cref{cor:max} which, as shown by \cref{thm:nosmallm} and
\cref{thm:manyclust}, cannot be easily generalised. A natural approach is to fix
large $r>0$ and use \cref{cor:min} -- then the product of those factors in
\eqref{eq:posterior} which come from the clusters that intersect $B(\bm{0},r)$ may be
well approximated using \cref{lem:gliwcan}, since by
\cref{cor:min} there are finitely many clusters intersecting $B(\bm{0},r)$ and the
number of observations in the cluster is comparable with $n$ for each cluster. Unfortunately in this way there is
no control over the impact of the clusters outside $B(\bm{0},r)$ as there are no
lower bounds on their size and upper bounds on their number. However the author
believes that these obstacles are possible to overcome and this remains subject
for the future work.

It should be also underlined the setting of our analysis was not the usual one
for the consistency analysis. 
Indeed, in our formulation of the CRP model our parameter space
is the space of partitions of $[n]$, which is changing with $n$. To
perform a classical consistency analysis we need the parameter space to be fixed
regardless of the number of observations. On the other hand, if we consider the DPMM
formulation, in which the parameter space is the space of all possible realisations
of the Dirichlet Process (i.e. the space of discrete measures on $\R^d$ with
infinitely many atoms) then again our input should come from an
\emph{infinite} mixture of normals, which was not the case in our examples.

However some of our results
from \cref{sec:basic} 
can be applied when the input sequence is a realisation of the DPMM. Indeed, the
convexity result of \cref{prop:MAPconv} does not have any assumptions on the data
sequence. As for \cref{prop:xmininter}, it requires the sequence of mean squared
norms to be bounded. It is easy to prove (see Supplement A) that for
a realisation of the DPMM this assumption holds almost surely and hence for
every $r>0$ the clusters intersecting $B(\bm{0},r)$ in the sequence of the MAP
partitions constructed on the sample from DPMM are of size comparable with
the number of observation and their number is bounded. However, some fundamental
questions remain unanswered in this case (e.g. does the number of clusters in
the MAP partition tend to infinity in this case?) and they are open for further investigation.

Note that the machinery presented can be used for a~different task. The
$P$-partitions that maximise the function $\Delta$ seem to be interesting objects in
their own right. Note that for dimension greater than 1 it seems to be extremely difficult to derive
the maximisers simply by analytical means.
\cref{thm:famconv} and~\cref{thm:distto0} give us a~convenient tool to examine
those maximisers as they may be approximated by performing sampling from the
posterior. This cannot be done faithfully as the normalizing constant in the
formula~\eqref{eq:posterior} cannot be computed explicitly, however there are
standard MCMC techniques that can be applied there (e.g. \cite{bib:neal}).
Further examination of the maximisers of the function $\Delta$ 
is left for future research.

\medskip
\small \textbf{Acknowledgement} The author wishes to express his thanks to an anonymous Referee and the
Associate Editor from the Bayesian Analysis for their helpful comments
and many thoughtful suggestions, which made the text far more readable and better than it
was. Also, the author thanks dr John Noble for helpful discussions.

\newpage
\renewcommand{\thesection}{A}
\begin{center}
\Large\textbf{Supplement A}\\
\small Proofs
\end{center}

\section*{Proof of \cref*{prop:normalForm}}
Is is easy to see that in \cref*{eq:DPMgen}
when the $F_\theta$'s have densities $f_\theta$ with respect to some
$\sigma$-finite measure $\nu$, sampling of $\theta$ may be omitted by taking the marginal distribution of
$\bx_J$ under $\cJ$. Then the model takes the form
\begin{equation}
\begin{array}{rcll}
\cJ&\sim&\CRP(\alpha)_n&\\
\bx_J\cond \cJ&\iid&F^{G_0}_J&\textrm{for $J\in \cJ$}
\end{array}
\end{equation}
where for $F^{G_0}_J$ is a~probability distribution on $\cX^{|J|}$ with the
density
\begin{equation}
f^{G_0}_J(\bx_J)=\int_{\Theta}\prod_{j\in J}f_\theta(x_j)\d{G_0(\theta)}
\end{equation}
with respect to product measure $\nu^{|J|}$. We now compute the exact formula
for $f^{G_0}_J$ when $G_0=\Normal(\bm{0},T)$ and $F_\theta=\Normal(\theta,\Sigma)$.

\smallskip
In order to simplify computations it is
convenient to use the notation $[A]=A^t A$, where $A$ is a~matrix of any
dimensionality. This is ambiguous as it is the same as the notation
introduced in the main text, $[n]=\{1,2,\ldots,n\}$, where $n$ is a~natural
number. However we don't use the latter in the following proof. 

Let $U=T^{-1/2}$. Since the densities of $\theta$ and $\bx_J\cond\theta$ are given by
\begin{equation}
\theta\sim
\frac{\det U}{\sqrt{2\pi}^d}
\exp\left\{-\re{2}[  U\theta ]\right\}
\quad\textrm{and}\quad
\bx_J\cond\theta\sim
\left(\frac{\det R}{\sqrt{2\pi}^d}\right)^{|J|} 
\exp\left\{-\re{2}\sum_{j\in J} [  R(x_j-\theta) ]\right\}
\end{equation}
we obtain that
\begin{equation}\label{eq:normderiv1}
f_J^{G_0}(\bx_J)=
\frac{\det U}{\sqrt{2\pi}^d}
\left(\frac{\det R}{\sqrt{2\pi}^d}\right)^{|J|} 
\int_\Theta
\exp\left\{-\re{2}[  U\theta ]
-\re{2}\sum_{j\in J}[  R(x_j-\theta) ]\right\}\d{\theta}.
\end{equation} 
Let $H_J$ be a positive definite symmetric matrix such that $[H_J]=|J|\cdot
[R]+[U]$, then
\begin{equation}
\begin{split}
[  U\theta ]+\sum_{j\in J}[ R(x_j-\theta) ]
&=[  U\theta ]+ 
\Big(\sum_{j\in J}[  R x_j ]-2\theta^t [  R ]\sum_{j\in J} x_j+|J| [  R\theta ]\Big)=\\
&= [  H_J\theta ] -2\theta^t [  R ]\sum_{j\in J} x_j+
\sum_{j\in J} [  R x_j ]=\\
&=\Big[ H_J \big(\theta-[ H_J]^{-1}[  R ]\sum_{j\in J} x_j\big)\Big]
-\Big[ H_J^{-1} [  R ]\sum_{j\in J} x_j\Big] +\sum_{j\in J} [  R x_j ]
\end{split}
\end{equation}
and hence
\begin{equation}\label{eq:normderiv2}
\begin{split}
\int_\Theta
\exp\left\{-\re{2}[  U\theta ]
-\re{2}\sum_{j\in J}[  R(x_j-\theta) ]\right\}\d{\theta}=
\frac{\sqrt{2\pi}^d}{\det H_J} \exp\left\{\re{2}\Big(
\Big[ H_J^{-1}[  R ]\sum_{j\in J} x_j\Big]
- \sum_{j\in J} [  R x_j ]\Big)\right\}.
\end{split}
\end{equation}
By joining equalities \eqref{eq:normderiv1} with \eqref{eq:normderiv2} and
substituting $H_J=\sqrt{|J|}\cdot R_J$ we obtain
that
\begin{equation}
f_J^{G_0}(\bx_J)=
\Big(\frac{\det R}{\sqrt{2\pi}^d}\Big)^{|J|}\frac{\det U}{\sqrt{|J|}^d\det
R_{|J|}}
\exp\Big\{\frac{1}{2}\Big(|J|\cdot\znorm{ R_{|J|}^{-1} R^2 \overline{\bx_J}}^2
-\sum_{j\in J} \znorm{ Rx_j}^2 \Big) \Big\}.
\end{equation}
Therefore
\begin{equation}\label{eq:bayes1}
\begin{split}
\bm{x}\cond\cJ&\sim
\prod_{J\in\cJ} f^{G_0}_J (\bm{x}_J)=\\
&=\Big(\frac{\det R}{\sqrt{2\pi}^d}\Big)^{n}
\exp\Big\{-\sum_{j\leq n} \znorm{ Rx_j}^2\Big\}(\det U)^{|\cJ|}
\prod_{J\in\cJ}
\frac{1}{\sqrt{|J|}^d\det R_{|J|}}
\exp\Big\{\frac{1}{2}\Big(|J|\cdot\znorm{ R_{|J|}^{-1} R^2 \overline{\bx_J}}^2
\Big) \Big\}\propto\\
&\propto(\det U)^{|\cJ|}
\prod_{J\in\cJ}
\frac{1}{\sqrt{|J|}^d\det R_{|J|}}
\exp\Big\{\frac{1}{2}\Big(|J|\cdot\znorm{ R_{|J|}^{-1} R^2 \overline{\bx_J}}^2 \Big) \Big\}.
\end{split}
\end{equation}
It is easy to see that 
the probability weights in $\CRP(\alpha)_n$ are given by the formula
\begin{equation}\label{eq:prior}
\P(\cJ_n=\cJ)=\frac{\alpha^{|\cJ|}}{\alpha^{(n)}} \prod_{J\in
\cJ}(|J|-1)!\propto
\alpha^{|\cJ|} \prod_{J\in \cJ}(|J|-1)!,
\end{equation}
\nopagebreak
where $a^{(k)}=a(a+1)\ldots (a+k-1)$ and $|\cJ|$ is the number of sets in the
partition $\cJ$. The proof of \cref*{prop:normalForm} follows from
\eqref{eq:bayes1}, \eqref{eq:prior} and the Bayes formula.

\section*{Proof of \cref*{prop:MAPconv}}
Take any $I,J\in \hat{\cJ}_n$. Consider all partitions of $[n]$ that are
obtained by replacing sets $I,J$ with $\tilde{I}$ and $\tilde{J}$ that satisfy
$|\tilde{I}|=|I|$, $|\tilde{J}|=|J|$ and $\tilde{I}\cup\tilde{J}=I\cup J$. Note
that by such operation we do not alter either the number of clusters or the
size of the clusters and therefore the posterior probability of such partitions is an
increasing function of 
\begin{equation}\label{eq:quadratic}
\begin{split}
|I|\cdot \znorm{ R_{|I|}^{-1}R^2
\overline{x_{\tilde{I}}}}^2+
|J|\cdot \znorm{ R_{|J|}^{-1}R^2
\overline{x_{\tilde{J}}}}^2&=
\re{|I|}\cdot \znorm{ R_{|I|}^{-1}R^2
\sum_{i\in\tilde{I}}x_i}^2+
\re{|J|}\cdot \znorm{ R_{|J|}^{-1}R^2
\Big(S-\sum_{i\in \tilde{I}}x_i\Big)}^2
\end{split}
\end{equation}
where $S=\sum_{k\in \tilde{I}\cup\tilde{J}}x_k=\sum_{k\in I\cup J}x_k$. It may be 
seen quite easily that~\eqref{eq:quadratic} defines a~strictly convex quadratic function with respect to
$\sum_{i\in\tilde{I}}x_i$. We investigate its value over a~finite number of possible
replacements. Therefore it achieves its maximal value at the vertices of convex hull of all
possible values of $\sum_{i\in\tilde{I}}x_i$. Since $\hat{\cJ}_n$ is the MAP
partition it follows that
$\sum_{i\in I}x_i$ maximises~\eqref{eq:quadratic}. 

Suppose that $\conv\{x_i\colon i\in I\}$ and
$\conv\{x_j\colon j\in J\}$ have a point in common, which is not $x_i$ for any
$i\leq n$. Then there exist
two equal convex combinations of points in $\{x_i\colon i\in I\}$
and $\{x_j\colon j\in J\}$, at least one of which is non-trivial, i.e.
\begin{equation}
\sum_{i\in I}\lambda_i^I x_i=\sum_{j\in J}\lambda_j^J x_j,\ \ 
\sum_{i\in I}\lambda_i^I =\sum_{j\in J}\lambda_j^J =1,\ \ 
\lambda_i^I ,\lambda_j^J \in [0,1).\ \ 
\end{equation}
(a convex combination is non-trivial if at least two of `lambdas' are non-zero). From this we can deduce that
\begin{equation}\label{eq:convexNew}
\begin{split}
\sum_{i\in I} x_i &= 
\sum_{i'\in I}\Big(\lambda^I_{i'}\sum_{i\in I} x_i\Big)=
\sum_{i\in I} \lambda^I_i x_i 
+\sum_{\substack{(i,i')\in I^2\\ i\neq i'}}\lambda^I_{i'}x_i=
\sum_{j\in J} \lambda^J_j x_j 
+\sum_{\substack{(i,i')\in I^2\\ i\neq i'}}\lambda^I_{i'}x_i=\\
&=\sum_{j\in J} \lambda^J_j \Big(\sum_{i'\in I}\lambda^I_{i'}x_j\Big) 
+\sum_{\substack{(i,i')\in I^2\\ i\neq i'}}\lambda^I_{i'}\Big(\sum_{j\in
J}\lambda^J_j x_i\Big)=
\sum_{(i',j)\in I\times J} \lambda^I_{i'}\lambda^J_{j}\Big(\sum_{i\in
(I\sm\{i'\})\cup\{j\}}x_i\Big).
\end{split}
\end{equation}
Moreover $\lambda_{i'}^I\lambda^J_j\in [0,1)$ for $i'\in I$, $j\in J$ and
$
\sum_{(i',j)\in I\times J} \lambda^I_{i'}\lambda^J_{j}=
\sum_{i'\in I} \lambda^I_{i'} \cdot \sum_{j\in J} \lambda^J_{j}=1,
$
so~\eqref{eq:convexNew} gives a~representation of $\sum_{i \in I}x_i$ as a
non-trivial (since at least one of the two combinations was non-trivial) convex combination of $\sum_{i\in \tilde{I}} x_i$. This is a
contradiction and the proof follows.

\section*{Proof of \cref*{prop:xmininter}}
For the reader's convenience the proof is split into three parts. 
In Subsection \nameref{subsec:prelimlem} we list some facts important for further
analysis. \nameref{subsec:impMAP} presents lemmas regarding the MAP, which are
further used in Subsection \nameref{subsec:proofxmininter}, where the proof of
one of the main results of the paper is presented.
\subsection*{Preliminary lemmas}\label{subsec:prelimlem}
\begin{rem}\label{rem:posprop}
Symmetric, positive definite matrices have the following
properties
\vspace*{-.3cm}
\begin{enumerate}[(a)]
\setlength{\itemsep}{1pt}
\setlength{\parskip}{0pt}
\setlength{\parsep}{0pt}
\item the sum of symmetric positive definite matrices is symmetric positive definite.
\item the inverse of symmetric positive definite matrix is symmetric positive definite.
\item for each symmetric positive matrix $A$ there exist an uniquely defined
symmetric positive matrix $B$ such that $A=B^t B$. We use the notation $B=A^{1/2}$.
\item if $A,B$ are symmetric positive definite matrices and also $A-B$ is
symmetric positive definite then $B^{-1}-A^{-1}$ is symmetric positive definite.
\item if $A,B$ are positive definite then $\det(A+B)\geq\det A$.
\end{enumerate}
\end{rem}

\begin{proof}
Let $A,B\in\R^{d,d}$.
\begin{enumerate}[(a)]
\item If $A,B$ are symmetric then $A+B$ is also symmetric. If $A,B$ are positive
definite then for every $x\in\R^d\sm\{\vec{0}\}$ we have $x^t(A+B)x=x^tAx+x^tBx> 0$ and
hence $A+B$ is also positive definite.
\item If $A$ is symmetric then $A^{-1}$ is also symmetric. If $A$ is positive
definite then by Theorem 7.1 from \cite{bib:zhang} it may be expressed as
$U^*\diag{\lambda_1,\ldots,\lambda_d}U$ where $U$ is unitary matrix and
$U^*$ its conjugate transpose and $\lambda_1,\ldots,\lambda_d>0$. Therefore
$A^{-1}=U^*\diag{\lambda_1^{-1},\ldots,\lambda_d^{-1}}U$ and again by using
Theorem 7.1 we obtain that $\inv{A}$ is positive definite.
\item Since if $A$ is a symmetric matrix then $A^t A=A^2$ and this point is an easy
consequence of Theorem 7.4 in \cite{bib:zhang}.
\item Let $P$ be a symmetric matrix that satisfy $P^2=B$. Positive definiteness
of $A-B$ is equivalent to $x'Ax>x'Bx$ for all $x\in\R^d$. By substituting
$y=Px$ this is equivalent to $y'\inv{P}A\inv{P}y>y'y$ for all $y\in\R^d$. Note
that $\inv{P}A\inv{P}$ is positive definite (as a product of positive definite
matrices) and hence it can be expressed as $U^*\Lambda U$. But then the latest
condition can be expressed as $z'z>z'U^*\inv{\Lambda} Uz$ for all $z\in\R^d$
which in the same way is equivalent to the positive definiteness of
$B^{-1}-A^{-1}$.
\item This is clearly equivalent to $\det(\mI+B\inv{A})\geq \det(\mI)=1$. As
$B\inv{A}$ is positive definite then for every eigenvalue $\lambda$ of
$\mI+B\inv{A}$ we have $\lambda=v'(\mI+B\inv{A})v>1$, where $v$ is its
eigenvector of norm 1. Therefore the determinant of $\mI+B\inv{A}$ is also
greater than 1.
\end{enumerate}
\end{proof}

\begin{rem}\label{prop:Dmprop}
Let  $R_m$ be defined as in the statement of \cref*{prop:normalForm}, then
\vspace*{-.3cm}
\begin{enumerate}[(a)]
\setlength{\itemsep}{1pt}
\setlength{\parskip}{0pt}
\setlength{\parsep}{0pt}
\item $\det R_m\to\det R$
\item if $y_m\to y$ then $R_m y_m\to Ry$
\end{enumerate}
\end{rem}
\begin{proof}
The proof is straightforward and therefore omitted.
\end{proof}

\begin{lem}\label{lem:facprod}
Let $n_1,\ldots,n_k\in\N$, $n_1\leq n_2\leq \ldots\leq n_k$ and $n=\sum_{i=1}^k n_i=an_k+r$,
where $a\in\N$, $r<n_k$. Then
$\prod_{i=1}^{k}n_i! \leq (n_k!)^{a}n_k(n_k-1)\ldots (n_k-r+1)$.
\end{lem}
\begin{proof}
We prove by induction on $n_k$ that the
sequence
$$[n_k]_{a,r}=(\underbrace{1,\ldots,n_k,1,\ldots,n_k,\ldots,1,\ldots,n_k}_{a},n_k-r+1,n_k-r+2,\ldots,n_k)$$
may be ordered so that it is term-wise not less than
$\mathbf{c}=( 1,\ldots,n_{1},1,\ldots,n_2,\ldots,1,\ldots,n_{k} )$. Clearly
the existence of such ordering establishes the lemma. For $n_k=1$ this is self
evident. For $n_k>1$ we apply `greedy' approach -- put all $a+1$ (or
$a$ in case $n_k|n$) $n_k$'s in
the places of $n_{k},n_{k-1},\ldots,n_{k-a+1}$. The fact that $n_k\geq
n_{k-1}\geq\ldots\geq n_1$ ensures that it is possible and all of
$n_k-1,n_{k-1}-1,\ldots,n_{k-a+1}-1,n_{k-a},\ldots,n_1$ are less or equal to
$n_k-1$.
Therefore we may apply inductive assumptions to these numbers thus finishing the
proof of the lemma.
\end{proof}

\begin{lem} \label{lem:sqrtn}
For every $\eps>0$ there exist $K\in\N$ such that if $n_1,\ldots,
n_k\leq n/K$, where $n=\sum_{i=1}^k n_i$, then $\sqrt[n]{\prod_{i=1}^k
n_i!/n!}<\eps$.
\end{lem}
\begin{proof}
Assume that $n_1\leq\ldots\leq n_k\leq n/K$ and let $n=an_k+r$, where
$0\leq r<n_k$. By \cref{lem:facprod} we get that
\begin{equation}
\frac{\prod_{i=1}^k n_i!}{n!}\leq \frac{(n_k!)^a (n_k-r+1)\ldots n_k}{n!}\leq
\frac{1}{1^{ n_k }}\cdot\frac{1}{2^{n_k}}\cdot\ldots\cdot\frac{1}{a^{n_k}}\cdot
\frac{1}{(a+1)^r}\leq \frac{1}{(a!)^{n_k}}.
\end{equation}
Therefore
\begin{equation}
\sqrt[n]{\frac{\prod_{i=1}^k n_i!}{n!}}\leq \re{\sqrt[\frac{n}{n_k}]{a!}}=
\re{\sqrt[\frac{n}{n_k}]{\lfloor \frac{n}{n_k} \rfloor!}}.
\end{equation}
For $K$ large enough this might be arbitrarily small, so the proof follows.
\end{proof}

\subsection*{Important properties of the MAP partition}\label{subsec:impMAP}
Let us fix a sequence $(x_n)_{n=1}^\infty$ in $\R^d$ and let
$\hat{\cJ}_n=\hat{\cJ}(x_1,\ldots,x_n)$. In order to facilitate the analysis, we introduce the following notation: let
$m_n=\min_{J\in\hat{\cJ}_n}|J|$ and $M_n=\max_{J\in\hat{\cJ}_n}|J|$ be the minimum and the
maximum cluster size in the partition $\hat{\cJ}_n$. 
Moreover for $r>0$ let 
\begin{equation}
\mr_n=\min\{|J|\colon J\in \hat{\cJ}_n,\norm{\overline{x_J}}<r\},\quad
\Mr_n=\max\{|J|\colon J\in \hat{\cJ}_n,\norm{\overline{x_J}}<r\}
\end{equation}
be the minimal and the maximal cluster size in the partition $\hat{\cJ}_n$ among the
clusters whose center of mass lies in $B(\bm{0},r)$.
Finally let
\begin{equation}
m^{[r]}_n=\min\{|J|\colon J\in \hat{\cJ}_n,\exists_{j\in J}\norm{x_j}<r\},\quad
M^{[r]}_n=\max\{|J|\colon J\in \hat{\cJ}_n,\exists_{j\in J}\norm{x_j}<r\}
\end{equation}
be the minimal and the maximal cluster size in the partition $\hat{\cJ}_n$ among the
clusters that intersect the ball $B(\bm{0},r)$.

Let $J^m_n,J^M_n\in\hat{\cJ}_n$ satisfy $|J^m_n|=m_n$ and $|J^M_n|=M_n$.
We define $J^{m,(r)}_n$, $J^{M,(r)}_n$, $J^{m,[r]}_n$ and $J^{M,[r]}_n$
similarly (e.g. $J^{m,(r)}_n\in \hat{\cJ}_n$ satisfies
$\znorm{\overline{\bm{x}_{J^{m,(r)}_n}}}<r$ and $|J^{m,(r)}_n|=m^{(r)}_n$).
\begin{propA}
\label{prop:xmax}
If $\big(\re{n}\sum_{i=1}^n \norm{x_i}^2\big)_{n=1}^\infty$ is bounded then
$\liminf_{n\to\infty}M_n/n>0.$
\end{propA}

\begin{proof}
Suppose that $\liminf M_n/n=0$. Then there exists an increasing
sequence $(n_k)_{k\in\N}$ such that $M_{n_k}/n_k<1/k$ for every
$k\in\N$.  
We now prove that
$$\lim_{k\to\infty}\sqrt[n_k]{\mQ{\hat{\cJ}_{n_k}}/\mQ{\{[n_k]\}}}=0,$$ hence obtaining
a contradiction with the definition of the MAP partition. 
By~\eqref{eq:posterior}
\begin{equation}\label{eq:Qquotient}
\begin{split}
\sqrt[n_k]{\mQ{\hat{\cJ}_{n_k}}/\mQ{[{n_k}]}}&=
\sqrt[{n_k}]{C^{|\hat{\cJ}_{n_k}|}/C}\cdot
\sqrt[{n_k}]{\prod_{J\in\hat{\cJ}_{n_k}}|J|!/{n_k}!}\cdot
\sqrt[{n_k}]{\frac{{n_k}^{(d+2)/2}\det
R_{n_k}}{\prod_{J\in\hat{\cJ}_{n_k}}|J|^{(d+2)/2}\det  R_{|J|}}}\cdot\\
&\ \ \cdot\exp\Big\{
\re{2n_k}\big(\sum_{J\in\hat{\cJ}_{n_k}}|J|\znorm{ R_{|J|}^{-1}R^2 \overline{\bx_J}}^2- 
n_k\znorm{ R_{n_k}^{-1}R^2 \overline{x_{[n_k]}}}^2\big)
\Big\}.
\end{split}
\end{equation} 

Firstly note that
\begin{equation}\label{eq:maxx1}
\limsup_{k\to\infty}\sqrt[{n_k}]{C^{|\hat{\cJ}_{n_k}|}/C}=
\limsup_{k\to\infty}C^{( |\hat{\cJ}_{n_k}|-1 )/{n_k}}\leq \max\{1,C\}.
\end{equation}

By \cref{lem:sqrtn}, it follows that, under the assumptions,
\begin{equation}\label{eq:maxx2}
\lim_{k\to\infty}\sqrt[{n_k}]{\prod_{J\in\hat{\cJ}_{n_k}}|J|!/{n_k}!}=0.
\end{equation}

From \cref{prop:Dmprop} 
\begin{equation}\label{eq:maxx3}
\limsup_{k\to\infty}\sqrt[{n_k}]{\frac{{n_k}^{(d+2)/2}\det
 R_{n_k}}{\prod_{J\in\hat{\cJ}_{n_k}}|J|^{(d+2)/2}\det  R_{|J|}}}\leq 
\frac{\limsup_{k\to\infty}\sqrt[{n_k}]{{n_k}^{(d+2)/2}\det  R_{{n_k}}}}
{\liminf_{k\to\infty} \sqrt[{n_k}]{\det  R^{|\hat{\cJ}_{n_k}|}}}\leq
\frac{1}{\min\{1,\det R\}}.
\end{equation}

Recall the inequality between linear and quadratic means which states
that for every sequence $y_1,\ldots,y_l$ of real numbers we have
\begin{equation}\label{eq:amqm}
\Big|\frac{\sum_{i=1}^l y_i}{l}\Big|\leq \sqrt{\frac{\sum_{i=1}^l y_i^2}{l}}
\Longleftrightarrow
l\cdot \Big(\frac{\sum_{i=1}^l y_i}{l}\Big)^2\leq \sum_{i=1}^l y_i^2.
\end{equation}
If we apply~\eqref{eq:amqm} to every coordinate of vectors
$\bm{y}_1,\ldots,\bm{y}_d\in \R^d$ and sum up obtained inequalities we
obtain that
\begin{equation}\label{eq:amqmmult}
l\cdot \xnorm{\frac{\sum_{i=1}^l \bm{y}_i}{l}}^2\leq \sum_{i=1}^l
\norm{\bm{y}_i}^2.
\end{equation}
Therefore, by linearity of multiplication by matrix
\begin{equation}
\begin{split}
\sum_{J\in\hat{\cJ}_n}|J|\znorm{ R_{|J|}^{-1}R^2 \overline{\bx_J}}^2&\leq
\sum_{J\in\hat{\cJ}_n}\sum_{j\in J}  \znorm{ R_{|J|}^{-1}R^2 x_j}^2
\end{split}
\end{equation}

and hence, using \cref{rem:posprop}, we have
\begin{equation}
\begin{split}
\sum_{J\in\hat{\cJ}_n}|J|\znorm{ R_{|J|}^{-1}R^2 \overline{\bx_J}}^2&\leq
\sum_{J\in\hat{\cJ}_n}\sum_{j\in J}  \znorm{ R_{|J|}^{-1}R^2 x_j}^2\leq
\sum_{J\in\hat{\cJ}_n}\sum_{j\in J}  \znorm{ R^{-1}R^2 x_j}^2\leq
\norm{R}_2^2\sum_{i\in [n]}\norm{x_i}^2,
\end{split}
\end{equation}
where $\norm{\cdot}_2$ is a~matrix norm induced by $\norm{\cdot}$ (i.e.
$\norm{A}_2=\sup_{\norm{x}=1}\norm{Ax}$). From this and assumptions of the
Proposition we can easily deduce that
\begin{equation}\label{eq:maxx4}
\re{n_k}\big(\sum_{J\in\hat{\cJ}_{n_k}}|J|\znorm{ R_{|J|}^{-1}R^2 \overline{\bx_J}}^2- 
n_k\znorm{ R_{n_k}^{-1}R^2 \overline{x_{[n_k]}}}^2\big)
\quad \textrm{is bounded from above.}
\end{equation}

Gathering~\eqref{eq:Qquotient},~\eqref{eq:maxx1},~\eqref{eq:maxx2},~\eqref{eq:maxx3} and
\eqref{eq:maxx4} together, we obtain that
$$\limsup\sqrt[n_k]{\mQ{\hat{\cJ}_{n_k}}/\mQ{[n_k]}}=0.$$ 
Hence there exists a
sufficiently large $n$ that satisfies $\P(\hat{\cJ}_n\cond \textbf{x})<\P(\{[n]\}\cond \textbf{x})$. This is a~contradiction.
\end{proof}

\begin{cor}\label{cor:centerMax}
If $\big(\re{n}\sum_{i=1}^n \norm{x_i}^2\big)_{n=1}^\infty$ is bounded then
there exist $r_0>0$ such that $\norm{\overline{x_{J^M_{n}}}}\leq r_0$ for all
$n>0$.
\end{cor}
\begin{proof}
By \cref{prop:xmax} we know that $\gamma:=\liminf_{n\to\infty}
M_n/n>0$, so there exists $N>0$ such that $M_n/n>\gamma/2$ for $n>N$. Suppose
that there exists a sequence $(n_k)_{k=1}^\infty$ such that
$\norm{\overline{x_{J^M_{n_k}}}}\geq k$ for all $k\in\N$. Note
that for $n_k>N$
\begin{equation}
\re{n_k}\sum_{i=1}^{n_k} \norm{x_i}\geq
\re{n_k}\sum_{i\in J^M_{n_k}} \norm{x_i}\geq
\re{n_k}\xnorm{\sum_{i\in J^M_{n_k}} x_i}=
\frac{M_{n_k}}{n_k}\znorm{\overline{x_{J^M_{n_k}}}}\geq
\gamma/2\cdot k,
\end{equation}
which, together with the inequality between the arithmetic and quadratic mean,
contradicts the assumption that the sequence $\big(\re{n}\sum_{i=1}^n
\norm{x_i}^2\big)_{n=1}^\infty$ is bounded. The proof of the Lemma follows from
the contradiction.
\end{proof}
\begin{propA}
\label{lem:xminbound}
If $\big(\re{n}\sum_{i=1}^n \norm{x_i}^2\big)_{n=1}^\infty$ is bounded then
$\liminf_{n\to\infty}m^{(r)}_n/n>0$ for every $r>0$.
\end{propA}
\begin{proof}
Firstly note that it is enough to prove the statement of \cref{lem:xminbound}
for all $r>r_0$ for some given $r_0>0$ -- indeed, $\mr$ is decreasing with
$r$. We take $r_0$ from the statement of \cref{cor:centerMax}. 

Fix $r>r_0$. Note that
$J^{M,(r)}_n=J^M_n$ and hence $\lim\inf_{n\to\infty}\Mr_n/n>0$. Now we
prove that $\lim\inf_{n\to\infty}\mr_n/n>0$. Suppose the contrary. We show that for sufficiently large $n$, the posterior probability of
$\hat{\cJ_n}$ increases if we create one cluster out of $J^{m,(r)}_n$ and
$J^{M,(r)}_n$. Let $\tilde{\cJ}_n$ be a~partition of $[n]$ obtained from $\hat{\cJ}_n$ by
joining $J_n^{m,(r)}$ with $J_n^{M,(r)}$, i.e.
\begin{equation}
\tilde{\cJ}_n=\hat{\cJ}_n \sm \{J_n^{m,(r)},J_n^{M,(r)}\}\cup\{J_n^{m,(r)}\cup
J_n^{M,(r)}\}.
\end{equation}
In order to simplify the notation, we write $m,M$ instead of
$m^{(r)}_n,M^{(r)}_n$ respectively and remember that they are both functions of
$n$. Similarly let us write $\overline{\bx_m},\overline{\bx_M}$ and
$\overline{\bx_{m\cup M}}$ instead of
$\overline{\bx_{J^{m,(r)}_n}},\overline{\bx_{J^{M,(r)}_n}}$ and
$\overline{\bx_{J^{m,(r)}_n\cup J^{M,(r)}_n}}$.
When taking a~quotient
$\P(\hat{\cJ}_n\cond\textbf{x})/\P(\tilde{\cJ}_n\cond\textbf{x})$ most
factors in~\eqref{eq:posterior} cancel out, giving
\begin{equation}\label{eq:xxx1}
\frac{\P(\hat{\cJ}_n\cond\textbf{x})}{\P(\tilde{\cJ}_n\cond\textbf{x})}
= C\frac{m!M!}{(m+M)!}\left( \frac{m+M}{mM} \right)^{(d+2)/2}\frac{\det
R_{m+M}}{\det R_m\cdot \det R_M}\cdot
\exp\left\{D_n\right\}^{1/2},
\end{equation}
where 
\begin{equation}
D_n=m\znorm{R_m^{-1}R^2\overline{\bx_{m}}}^2+
M\znorm{R_M^{-1}R^2\overline{\bx_{M}}}^2-(m+M)\znorm{R_{m+M}^{-1}R^2\overline{\bx_{m\cup
M}}}^2.
\end{equation}
It is now straightforward to verify that
\begin{equation}\label{eq:maxhelp1}
\begin{split}
m\znorm{R\overline{\bx_{m}}}^2+
M\znorm{R\overline{\bx_{M}}}^2-(m+M)\znorm{R\overline{\bx_{m\cup
M}}}^2&=
\frac{mM}{m+M}\znorm{R(\overline{\bx_{m}}-\overline{\bx_{M}})}^2\leq\\
&\hspace{-3cm}\leq m\znorm{R(\overline{\bx_{m}}-\overline{\bx_{M}})}^2\leq
m\norm{R}_2^2\cdot 4r^2.
\end{split}
\end{equation}
Moreover
\begin{equation}
\begin{split}
(m+M)\mI-(m+M) R(R_{m+M}^{-1} )^2 R&=
(m+M)(\mI-(R_{m+M}^{-1}R)^t R_{m+M}^{-1}R)R=\\
&\hspace{-3cm}= (m+M)R\Big(\mI-\big(\mI+(UR^{-1})^t UR^{-1}/(m+M)\big)^{-1}\Big)R=\\
&\hspace{-3cm}= R(UR^{-1})^t UR^{-1}\big(\mI+(UR^{-1})^t UR^{-1}/(m+M)\big)^{-1}R
\end{split}
\end{equation}
and therefore
\begin{equation}\label{eq:maxhelp2}
\limsup_{n\to\infty}\big((m+M)\znorm{R\overline{\bx_{m\cup
M}}}^2-(m+M)\znorm{R_{m+M}^{-1}R^2\overline{\bx_{m\cup
M}}}^2\big)\leq \norm{U}_2^2 r^2.
\end{equation}
By \cref{rem:posprop}, together with~\eqref{eq:maxhelp1} and
\eqref{eq:maxhelp2},
\begin{equation}\label{eq:xsum}
\limsup_{n\to\infty} D_n
\leq m\norm{R}_2^2\cdot 4r^2+\norm{U}_2^2 r^2.
\end{equation}
Stirling formula, which is valid for every $n\in\N$ (cf. \cite{bib:feller}),
states that
\begin{equation}\label{eq:stirling}
\sqrt{2\pi n}(n/e)^n e^{\re{12n+1}}<n!<\sqrt{2\pi n}(n/e)^n e^{\re{12n}}.
\end{equation}
This gives:
\begin{equation}\label{eq:xxx2}
\frac{m!M!}{(m+M)!}\leq 
\sqrt{2\pi}\left(\frac{mM}{m+M}\right)^{1/2}\frac{m^m M^M
}{(m+M)^{(m+M)}}e\leq
\sqrt{2\pi}e\left(\frac{mM}{m+M}\right)^{1/2}\left(\frac{m}{M}\right)^m.
\end{equation}
Now by applying~\eqref{eq:xsum} and~\eqref{eq:xxx2} to~\eqref{eq:xxx1} we obtain that for some
constants $C'$ and $C''$
\begin{equation}\label{eq:finalmM}
\liminf_{n\to\infty}\frac{\P(\hat{\cJ}_n\cond\textbf{x})}{\P(\tilde{\cJ}_n\cond\textbf{x})}\leq
\liminf_{n\to\infty}C'\left(\frac{m+M}{mM}\right)^{(d+1)/2}\left(\frac{mC''}{M}\right)^m=
0,
\end{equation}
as $\liminf_{n\to\infty} m/M\to 0$.
Hence there exist $n$ such that the posterior probability of $\hat{\cJ}_n$ is smaller
than the posterior probability of $\tilde{\cJ}_n$. This contradicts the
definition of $\hat{\cJ}_n$ and finishes the proof of the Lemma.
\end{proof}

\subsection*{Proof of \cref*{prop:xmininter}}\label{subsec:proofxmininter}
Assume that  $\big(\re{n}\sum_{i=1}^n \norm{x_i}^2\big)_{n=1}^\infty$ is
bounded. We want to prove that
$\liminf_{n\to\infty}m_n^{[r]}/n>0$ for every $r>0$.

\smallskip
Take $r_0$ from the statement of \cref{cor:centerMax}. Note that, as in proof of
\cref{lem:xminbound} it is enough to prove the statement of
\cref{prop:xmininter} for $r>r_0$. 

Fix $r>r_0$. Suppose that $\liminf_{n\to\infty}m_n^{[r]}/n=0$ and let
$(n_k)_{k=1}^\infty$ be a sequence such that
$\lim_{k\to\infty}m_{n_k}^{[r]}/n_k=0$. By \cref{lem:xminbound} we
obtain that $\lim_{k\to\infty} \znorm{\overline{x_{J^{m,[r]}_{n_k}}}}=\infty$
(otherwise we would obtain a contradiction). Let 
\begin{equation}
I_n^a = \{j\in J^{m,[r]}_n\colon \norm{x_j}\leq r\},\quad
I_n^b = \{j\in J^{m,[r]}_n\colon \norm{x_j}>r\}.
\end{equation}
Consider a partition $\check{\cJ}_n$ obtained from $\hat{\cJ}_n$ by taking
$I_n^a$ from $J^{m,[r]}_n$ and adding it to $J^M_n$, i.e.
\begin{equation}
\check{\cJ}_n=
\hat{\cJ}_n\sm \{J^{m,[r]}_n,J^M_n\}\cup
\{J^{m,[r]}_n\sm I^a_n, J^M_n\cup I^a_n\}.
\end{equation}
When taking a~quotient
$\P(\hat{\cJ}_n\cond\textbf{x})/\P(\check{\cJ}_n\cond\textbf{x})$ most
factors in~\eqref{eq:posterior} cancel out, giving
\begin{equation}\label{eq:interquot}
\frac{\P(\hat{\cJ}_{n_k}\cond\textbf{x})}{\P(\check{\cJ}_{n_k}\cond\textbf{x})}
= \frac{(a+b)!M!}{b!(a+M)!}\left( \frac{b(a+M)}{(a+b)M} \right)^{(d+2)/2}\frac{\det
R_b\cdot \det
R_{a+M}}{\det R_{a+b}\cdot \det R_M}\cdot
\exp\left\{\check{D}_{n_k}\right\}^{1/2},
\end{equation}
where $M=|J^M_{n_k}|$, $a=|I^a_{n_k}|$, $b=|I^b_{n_k}|$ (in order to simplify
the notation we skip the index $n_k$) and
\begin{equation}
\begin{split}
\check{D}_{n_k}&=(a+b)\znorm{R_{a+b}^{-1}R^2\overline{\bx_{a\cup b}}}^2+
M\znorm{R_M^{-1}R^2\overline{\bx_{M}}}^2-\\
&\quad - b\znorm{R_{b}^{-1}R^2\overline{\bx_{b}}}^2 -
(a+M)\znorm{R_{a+M}^{-1}R^2\overline{\bx_{a\cup M}}}^2.
\end{split}
\end{equation}
in which $\overline{\bx_{a\cup b}}=\overline{\bx_{I^a_{n_k}\cup I^b_{n_k}}}$ and
we define $\overline{\bx_b},\overline{\bx_M},\overline{\bx_{a\cup M}}$
similarly. Note that
\begin{equation}\label{eq:interto0}
\frac{(a+b)!M!}{b!(a+M)!}=\frac{(b+1)^{(a)}}{(M+1)^{(a)}}<
\frac{b+1}{M+1}\stackrel{k\to\infty}{\longrightarrow} 0,
\end{equation}
since $\lim_{k\to\infty}(a+b)/n_k=\lim_{k\to\infty}
m_{n_k}^{[r]}/n_k= 0$ and $\liminf_n M/n>0$. For the similar reason
\begin{equation}\label{eq:intersmall}
\frac{b(a+M)}{(a+b)M}< \frac{a+M}{M}\stackrel{k\to\infty}{\longrightarrow} 1.
\end{equation}
Now let us investigate $\check{D}_n$. The notation is easier after a linear
substitution $y_i=R^2 x_i$ (so that $\overline{\by_I}=R^2\overline{\bx_I}$),
hence obtaining 
\begin{equation}
\check{D}_{n_k}=(a+b)\znorm{R_{a+b}^{-1}\overline{\by_{a\cup b}}}^2+
M\znorm{R_M^{-1}\overline{\by_{M}}}^2-
b\znorm{R_{b}^{-1}\overline{\by_{b}}}^2 -
(a+M)\znorm{R_{a+M}^{-1}\overline{\by_{a\cup M}}}^2.
\end{equation}
Note that 
\begin{equation}\label{eq:Xinterexp1}
\begin{split}
(a+b)\znorm{R_{a+b}^{-1}\overline{\by_{a+b}}}^2 - 
b\znorm{R_{b}^{-1}\overline{\by_{b}}}^2&=
(a+b)\xnorm{R_{a+b}^{-1}\Big(\frac{a}{a+b} \overline{\by_{a}} +\frac{b}{a+b}
\overline{\by_{b}}\Big)}^2 - 
b\znorm{R_{b}^{-1}\overline{\by_{b}}}^2=\\
&=\re{a+b}\znorm{R_{a+b}^{-1}(a \overline{\by_{a}} +b \overline{\by_{b}})}^2 - 
b\znorm{R_{b}^{-1}\overline{\by_{b}}}^2\leq\\
&=\re{a+b}\big( a^2\znorm{R_{a+b}^{-1} \overline{\by_{a}}}^2 +
2ab\znorm{R_{a+b}^{-1} \overline{\by_{a}}}\znorm{R_{a+b}^{-1}\overline{\by_{b}}}
+b^2\znorm{R_{a+b}^{-1} \overline{\by_{b}}}^2 -\\ 
&\quad -b(a+b)\znorm{R_{b}^{-1}\overline{\by_{b}}}^2 \big)=\\
&=\re{a+b}\big( a^2\znorm{R_{a+b}^{-1} \overline{\by_{a}}}^2 +
2ab\znorm{R_{a+b}^{-1} \overline{\by_{a}}}\znorm{R_{a+b}^{-1}\overline{\by_{b}}}
+\overline{\by_{b}}T_1 \overline{\by_{b}}\big)
\end{split}
\end{equation}
where $T_1=b^2R_{a+b}^{-1} - b(a+b)R_b^{-1}$. 
For two positive/negative matrices $M_1,M_2$ we write $M_1\geq M_2$ when
$M_1-M_2$ is positive definite. Then
\begin{equation}\label{eq:intermatrix1}
\begin{split}
T_1&=b^2\big(R^2+U^2/(a+b)\big)^{-1} - b(a+b)(R^2+U^2/b)^{-1} =\\
&=b^2(a+b)\big((a+b)R^2+U^2\big)^{-1} - b^2(a+b)(bR^2+U^2)^{-1} =\\
&=b^2(a+b)\Big( \big((a+b)R^2+U^2\big)^{-1} - (bR^2+U^2)^{-1}\Big) =\\
&=-ab^2(a+b) \big((a+b)R^2+U^2\big)^{-1}R^2(bR^2+U^2)^{-1} =\\
&=-ab \big(R^2+U^2/(a+b)\big)^{-1}R^2(R^2+U^2/b)^{-1} \leq\\
&=-ab (R^2+U^2)^{-1}R^2(R^2+U^2)^{-1}=:-abT_2^2
\end{split}
\end{equation}
Using \eqref{eq:Xinterexp1} and \eqref{eq:intermatrix1} we have that
\begin{equation}\label{eq:interexp1}
\begin{split}
(a+b)\znorm{R_{a+b}^{-1}\overline{\by_{a+b}}}^2 - 
b\znorm{R_{b}^{-1}\overline{\by_{b}}}^2&\leq\re{a+b}\big( a^2\znorm{R_{a+b}^{-1} \overline{\by_{a}}}^2 +
2ab\znorm{R_{a+b}^{-1} \overline{\by_{a}}}\znorm{R_{a+b}^{-1}\overline{\by_{b}}}
-ab\znorm{ T_2\overline{\by_{b}} }^2\big) =\\ 
&=a\Big( \frac{a}{a+b}\znorm{R_{a+b}^{-1} \overline{\by_{a}}}^2 +
2\frac{b}{a+b}\big( \znorm{R_{a+b}^{-1} \overline{\by_{a}}}\znorm{R_{a+b}^{-1}\overline{\by_{b}}}
-\znorm{ T_2\overline{\by_{b}} }^2 \big)\Big) \leq\\ 
&\leq a\Big( \frac{a}{a+b}\znorm{R^{-1} \overline{\by_{a}}}^2 +
2\frac{b}{a+b}\big( \znorm{R^{-1} \overline{\by_{a}}}\znorm{R^{-1}\overline{\by_{b}}}
-\znorm{ T_2\overline{\by_{b}} }^2 \big)\Big) \leq\\ 
&\leq a\Big( \eig{R}^{-2}\frac{a}{a+b}\znorm{\overline{\by_{a}}}^2 +
2\frac{b}{a+b}\znorm{\overline{\by_{b}}}\big( \eig{R}^{-2}\znorm{ \overline{\by_{a}}}
-\eig{T_2}^2 \znorm{ \overline{\by_{b}} } \big)\Big)
\end{split}
\end{equation}
where by $\eig{A}$ we denote the minimal eigenvalue of the square matrix
$A$. 
Similarly we note that
\begin{equation}\label{eq:Xinterexp2}
\begin{split}
M\znorm{R_M^{-1}\overline{\by_{M}}}^2-
(a+M)\znorm{R_{a+M}^{-1}\overline{\by_{a\cup M}}}^2&\leq
M\znorm{R^{-1}\overline{\by_{M}}}^2-
(a+M)\znorm{R_{1}^{-1}\overline{\by_{a\cup M}}}^2=\\
&\hspace{-6cm}= M\znorm{R_M^{-1}\overline{\by_{M}}}^2-
\re{a+M}
\big( a^2\znorm{R_{a+M}^{-1}\overline{\by_{a}}}^2 +
2aMR_{a+M}^{-1}\overline{\by_{a}}\cdot R_{a+M}^{-1}\overline{\by_{M}}+ 
M^2\znorm{R_{a+M}^{-1}\overline{\by_{M}}}^2 
\big)=\\
&\hspace{-6cm}= \re{a+M}
\big((a+M)M\znorm{R_M^{-1}\overline{\by_{M}}}^2-
 a^2\znorm{R_{a+M}^{-1}\overline{\by_{a}}}^2 -
2aMR_{a+M}^{-1}\overline{\by_{a}}\cdot R_{a+M}^{-1}\overline{\by_{M}}-
M^2\znorm{R_{a+M}^{-1}\overline{\by_{M}}}^2 
\big)\leq\\
&\hspace{-6cm}\leq\re{a+M}
\Big(M\big( (a+M)\znorm{R_M^{-1}\overline{\by_{M}}}^2-
M\znorm{R_{a+M}^{-1}\overline{\by_{M}}}^2 \big) -
2aMR_{a+M}^{-1}\overline{\by_{a}}\cdot R_{a+M}^{-1}\overline{\by_{M}}
\Big).
\end{split}
\end{equation}
We can write
\begin{equation}\label{eq:intermaterix2}
\begin{split}
(a+M)(R^2+U^2/M)^{-1}-M\big(R^2+U^2/(a+M)\big)^{-1}&=\\
&\hspace{-4cm}=M(a+M)(MR^2+U^2)^{-1}-M(a+M)\big((a+M)R^2+U^2\big)^{-1}=\\
&\hspace{-4cm}=M(a+M)\Big((MR^2+U^2)^{-1}-\big((a+M)R^2+U^2\big)^{-1}\Big)=\\
&\hspace{-4cm}=aM(a+M)(MR^2+U^2)^{-1}R^2\big((a+M)R^2+U^2\big)^{-1}=\\
&\hspace{-4cm}=a(R^2+U^2/M)^{-1}R^2\big(R^2+U^2/(a+M)\big)^{-1}\leq\\
&\hspace{-4cm}\leq a(R^2)^{-1}R^2(R^2)^{-1}=aR^2.\\
\end{split}
\end{equation}
and hence, by \eqref{eq:Xinterexp2} and \eqref{eq:intermaterix2}
\begin{equation}\label{eq:interexp2}
\begin{split}
M\znorm{R_M^{-1}\overline{\by_{M}}}^2-
(a+M)\znorm{R_{a+M}^{-1}\overline{\by_{a\cup M}}}^2&\leq\re{a+M}
\big(aM \znorm{R\overline{\by_{M}}}^2  -
2aMR_{a+M}^{-1}\overline{\by_{a}}\cdot R_{a+M}^{-1}\overline{\by_{M}}
\big)=\\
&=a\frac{M}{a+M}
\big( \znorm{R\overline{\by_{M}}}^2  -
2R_{a+M}^{-1}\overline{\by_{a}}\cdot R_{a+M}^{-1}\overline{\by_{M}}
\big)\leq\\
&\leq a\frac{M}{a+M}
\big( \znorm{R\overline{\by_{M}}}^2  +
2\znorm{R^{-1}\overline{\by_{a}}}\cdot \znorm{R^{-1}\overline{\by_{M}}}
\big)
\end{split}
\end{equation}
Joining \eqref{eq:interexp1} and \eqref{eq:interexp2} we get that
\begin{equation}\label{eq:interDineq}
\begin{split}
\check{D}_{n_k}&\leq 
a\Big(\eig{R}^{-2}\frac{a}{a+b}\znorm{\overline{\by_{a}}}^2 +
2\frac{b}{a+b}\znorm{\overline{\by_{b}}}\big( \eig{R}^{-2}\znorm{ \overline{\by_{a}}}
-\eig{T_2}^2 \znorm{ \overline{\by_{b}} } \big)+\\
&\quad+\frac{M}{a+M}
\big( \znorm{R\overline{\by_{M}}}^2  +
2\znorm{R^{-1}\overline{\by_{a}}}\cdot \znorm{R^{-1}\overline{\by_{M}}}
\big)\Big)
\end{split}
\end{equation}
Note that since 
$\norm{\overline{\by_{a\cup b}}}
\stackrel{k\to\infty}{\longrightarrow}\infty$ and $\znorm{\overline{\by_a}}\leq
\Eig{R^2}\cdot\znorm{\overline{\bx_a}}\leq
\Eig{R^2}r$ ($\Eig{A}$ is the largest eigenvalue of the square matrix $A$) we have that
\begin{equation}\label{eq:interdown}
\frac{b}{a+b}\znorm{\overline{\by_{b}}}\geq
\znorm{\overline{\by_{a\cup b}}}-
\frac{a}{a+b}\znorm{\overline{\by_{a}}}\stackrel{k\to\infty}{\longrightarrow}\infty.
\end{equation}
Moreover $\znorm{\overline{\by_M}}\leq \Eig{R^2}\cdot\znorm{\overline{\bx_M}}\leq
\Eig{R^2}r$ and therefore by \eqref{eq:interDineq} and \eqref{eq:interdown} we have
\begin{equation}\label{eq:interDtoDown}
\lim_{k\to\infty} D_{n_k}=-\infty
\end{equation}
By taking \eqref{eq:interquot} and using \eqref{eq:interto0},
\eqref{eq:intersmall} and \eqref{eq:interDtoDown}	we obtain that
$\lim_{k\to\infty}\frac{\P(\hat{\cJ}_{n_k}\cond\textbf{x})}{\P(\check{\cJ}_{n_k}\cond\textbf{x})}=0$;
from this contradiction the proof follows.

\section*{Proof of \cref*{thm:distto0}}
Let $\cK_r$ be the space of all closed and convex subsets of $B(\bm{0},r)$. 
Note that $\bm{M}_\Delta$ is closed in $(F_K(\cK_r),\overline{d_P})$ as an intersection of the set of maximisers of $\bm{M}_\Delta$ in
$(F_K(\cK_r),\overline{d_P})$ and the subspace of $P$-partitions, both of them being
closed subspaces of $(F_K(\cK_r),\overline{d_P})$. By \cref*{thm:MAPmax} we know that
$\bm{E}\subseteq \bm{M}_\Delta$. 
Now the proof of \cref*{thm:distto0} follows from the following \cref{lem:topo}.
\begin{lem}\label{lem:topo}
Let $(\cX,d)$ be a~finitely compact metric space, $D\subseteq\cX$ a~closed set and
$(a_n)_{n=1}^\infty$ a~bounded sequence in $\cX$. If every converging subsequence
of $(a_n)_{n=1}^\infty$ has a~limit in $D$ then $\dist(a_n,D)\to 0$, where
$\dist(\cdot,\cdot)$ is the distance function, i.e.
$$\dist(x,D)=\inf_{y\in D} d(x,y).$$
\end{lem}

\begin{proof}
Suppose that $\limsup\dist(a_n,D)>0$. Then there exist a subsequence
$(a_{n_k})_{k=1}^\infty$ and $\eps>0$ such that $\dist(a_{n_k},D)>\eps>0$. This
contradicts the fact that $(a_{n_k})_{k=1}^\infty$ as a bounded sequence in
$\cX$ has a converging subsequence whose limit must belong to the closed set $D$.
\end{proof}

\section*{Proof of \cref*{prop:KMAPlarge}}
\begin{lem}\label{lem:philemma}
If $P$ is a measure on $(\R^d,\cB)$ with bounded support and absolutely continuous with respect to the
Lebesgue measure then for every $\alpha>0$
\begin{equation}\label{eq:phieq}
\Psi(\alpha):=\inf_{ \substack{A\in\cK_r\\ P(A)\geq\alpha} } 
\sup_{\substack{A_1,A_2\in \cB\\ A_1\cup A_2=A\\ A_1\cap A_2=\emptyset}}
P(A_1)\cdot P(A_2)\cdot \norm{E(A_1)-E(A_2)}^2>0
\end{equation}
where $E(B)=\int_B x\d{P(x)}/P(B)$ for $B\in \cB$.
\end{lem}
\begin{proof}
Fix $\alpha>0$. As an easy consequence of \cref*{thm:dPcomp} we obtain that $P(\cdot)$ is a
continuous function in $(\cK_r,\varrho_H)$. Therefore $\cK_r^\alpha:=\{A\in \cK_r\colon
P(A)\geq \alpha\}$ is a closed subspace of compact (by \cref*{thm:dHcomp}) topological space, therefore it
is compact itself. 

Assume that the support of $P$ is contained in the ball $B(\bm{0},r)$ and let
$r>1$. Consider the function
\begin{equation}
\varphi(A)=
\sup_{\substack{A_1,A_2\in \cB\\ A_1\cup A_2=A\\ A_1\cap A_2=\emptyset}}
P(A_1)\cdot P(A_2)\cdot \norm{E(A_1)-E(A_2)}^2\geq 0
\end{equation}
in the compact topological space $(\cK_r,\varrho_H)$. We prove that this function is continuous.

Firstly note that since we operate in a bounded space then if $P(B)\to 0$ then
$\int_B x\d{P(x)}\to 0$. From this it can be easily seen that for every
$\eps>0$ there exist $\delta>0$ such that if $d_P(A,B)<\delta$ then
$\norm{E(A)-E(B)}<\eps$ for $A,B\in \cK^\alpha_r$.

Fix $0<\eps<1$. There exist $\delta_1<\eps$ such that if $d_P(A,A')<\delta_1$ then
$\norm{E(A)-E(A')}<\eps/2$. There exist $\delta_2$ such that if
$\varrho_H(A,A')<\delta_2$ then $d_P(A,A')<\delta_1$ (this is because of
\cref*{thm:dPcomp} and the fact that $(\cK_r,\varrho_H)$ is compact and therefore
the continuity implies the uniform continuity).
Let us take $A,A'\in\cK_r$ such that $\varrho_H(A,A')<\delta_2$. Let $A_1,A_2\in
\cK_r$ be such that $A_1\cap A_2=\emptyset$, $A_1\cup A_2=A$ and
\begin{equation}\label{eq:ppe1}
\varphi(A)-\eps\leq P(A_1)\cdot P(A_2)\cdot \norm{E(A_1)-E(A_2)}^2
\end{equation}
Consider $A_1'=A_1\cup (A'\sm A)\sm (A\sm A')$ and $A_2'=A_2\sm (A\sm A')$. Then
$A_1'\cap A_2'=\emptyset$, $A_1'\cup A_2'=A'$ and 
\begin{equation}
d_P(A_1,A_1'),d_P(A_2, A_2')\leq d_P(A,A')\leq \delta_1.
\end{equation}
Therefore $|P(A_i)-P(A_i')|<\delta_1<\eps$, $\norm{E(A_i)-E(A_i')}<\eps/2$ for $i=1,2$.
Since $|P(A_i)|\leq 1$ and $\norm{E(A_i)}\leq r$ for $i=1,2$ we get
\begin{equation}\label{eq:ppe2}
\big|P(A_1)\cdot P(A_2)\cdot \norm{E(A_1)-E(A_2)}^2-
P(A_1')\cdot P(A_2')\cdot \norm{E(A_1')-E(A_2')}^2\big|
<50r^2\eps. 
\end{equation}
By \eqref{eq:ppe1} and \eqref{eq:ppe2} we obtain
\begin{equation}
\varphi(A)-\eps-50r^2\eps\leq P(A_1')\cdot P(A_2')\cdot \norm{E(A_1')-E(A_2')}^2
\leq \varphi(A').
\end{equation}
By symmetry we get $\varphi(A')-\eps-50r^2\eps\leq \varphi(A)$ which means that
$|\varphi(A)-\varphi(A')|<(1+50r^2)\eps$ for $\varrho_H(A,A')<\delta_2$ which
proofs the continuity of $\varphi$ in the topological space $(\cK_r^\alpha,\varrho_H)$.
Therefore by Weierstrass Theorem we get that
\begin{equation}\label{eq:phieq}
\inf_{ \substack{A\in\cK_r\\ P(A)\geq\alpha} } \varphi(A)=
\varphi(A_0)
\end{equation}
for some $A_0\in \cK_r$ such that $P(A_0)\geq \alpha$. It is easy to see that
$\varphi(A_0)>0$ (it is enough to divide $A_0$ into two subsets of positive
measure by a hyperplane so that the center of masses of two parts do not
coincide) and the Lemma follows.
\end{proof}

\noindent{\bf Proof of \cref*{prop:KMAPlarge}: }
Fix $K>0$. Let $\Psi$ be defined as in the statement of \cref{lem:philemma}. We first prove that for
$\eps=\re{8}e\Psi(\inv{K})$ if $\norm{\Sigma}<\eps$ then every finite maximiser
of the $\Delta$ function is of size larger than $K$. Take any finite partition
$\cG$ of $\R^d$ that consists of at most $K$ convex sets with positive $P$ measure. 
Let $A\in\cG$ be the set of the largest probability in $\cG$; note that
$P(A)\geq \inv{K}$. By definition of $\Psi$ we can divide $A$ into two sets $A_1,A_2$ ($A_1\cup A_2=A$, $A_1\cap
A_2=\emptyset$) such that
\begin{equation}
P(A_1)\cdot P(A_2)\cdot \norm{E(A_1)-E(A_2)}^2>
\Psi(K^{-1})/2.
\end{equation}
Let $\cG'=\cG\cup\{A_1,A_2\}\sm\{A\}$. Then
\begin{equation}
\begin{split}\label{eq:KPMAPlarge1}
\Delta(\cG')-\Delta(\cG)&=
\re{2}\big( P(A_1)\norm{R\cdot E(A_1)}^2+ P(A_2)\norm{R\cdot E(A_2)}^2-
P(A)\norm{R\cdot E(A)}^2 \big)-
\\&-P(A_1)\ln\re{P(A_1)}-P(A_2)\ln\re{P(A_2)}+
P(A)\ln\re{P(A)}.
\end{split}
\end{equation}
It is straightforward to verify that $p\ln \inv{p}\in [0,\re{e}]$ for $p\in
[0,1]$ and, since $P(A_1)E(A_1)+P(A_2)E(A_2)=P(A)E(A)$ we have
\begin{equation}
P(A_1)\norm{R\cdot E(A_1)}^2+ P(A_2)\norm{R\cdot E(A_2)}^2-
P(A)\norm{R\cdot E(A)}^2 =
\frac{P(A_1)P(A_2)}{P(A)}\norm{R\cdot(E(A_1)-E(A_2))}^2.
\end{equation}
Therefore by \eqref{eq:KPMAPlarge1} and \cref{lem:philemma} we get 
\begin{equation}
\begin{split}
\Delta(\cG')-\Delta(\cG)&\geq
\frac{P(A_1)P(A_2)}{P(A)}\norm{R\cdot(E(A_1)-E(A_2))}^2-2\inv{e}\geq\\
&\geq
\frac{P(A_1)P(A_2)}{P(A)}\re{\norm{\inv{R}}}\norm{E(A_1)-E(A_2)}^2-2\inv{e}=\\
&=\frac{P(A_1)P(A_2)}{P(A)}\re{\norm{\Sigma}}\norm{E(A_1)-E(A_2)}^2-2\inv{e}\geq\\
&\geq \inv{\eps}P(A_1)P(A_2)\norm{E(A_1)-E(A_2)}^2-2\inv{e}\geq\\
&\geq \inv{\eps}\Psi(\inv{K})/2-2\inv{e}>2\inv{e}>0.
\end{split}
\end{equation}
Hence $\cG$ is not a maximiser of $\Delta$ function. 

Now let $X_1,X_2,\ldots\iid P$ and $\hat{\cA}_n$ be the family of convex hulls of groups of observations
defined by the sequence of the MAP partitions based on $X_1,\ldots,X_n$ (where the
MAP partitions were computed in the model with the within group covariance matrix of the
norm less than $\eps$). Suppose that there exists a subsequence
$(n_i)_{i=1}^\infty$ such that $|\hat{\cA}_{n_i}|\leq K$ for $i\in \N$. By
the compactness of the space $(F_{\tilde{K}}(\cK_r),\overline{\varrho_H})$
(cf. \cref{thm:famconv}) we get that there is a subsequence
$(\hat{\cA}_{n_{ i_j }})$ that is convergent in this space to a $P$-partition
$\cE$ of $\R^d$ which is a maximiser of $\Delta$ (cf. \cref*{thm:MAPmax}). By our
previous analysis, $|\cE|>K$. On the other hand the probabilities of sets in
$\hat{\cA}_n$ are separated from 0 (this is a consequence of \cref{cor:min}) and
this yields a contradiction.

\section*{Proofs for \Asection{sec:exmp}}

\subsection*{\nameref*{subsec:segment}}
We now find the convex partition that maximises $\Delta$ if $P$ is a~uniform
distribution on $[-1,1]$. Since it is
convex it is defined by the lengths of consecutive
subsegments of $[-1,1]$; let those be $2p_1,\ldots,2p_n$. Let
$s_k=\sum_{i=1}^k p_i$ for $k\geq 1$ and $A_k=[s_{k-1},s_k]$, where $s_0=0$.
Then it follows that the optimal partition maximises
\begin{equation}
F(p_1,\ldots,p_n)=\rho\sum_{i=1}^n p_i(2s_{i-1}-1+p_i)^2+\sum_{i=1}^n p_i\ln p_i,
\end{equation}
where $\rho=R^2/2$, with the constraint $\sum_{i=1}^n p_i=1$. This problem can
be solved using Lagrange
multipliers. We are looking for the local maximum of a~function
\begin{equation}
F_\lambda(p_1,\ldots,p_n)=F(p_1,\ldots,p_n)-\lambda\sum_{i=1}^n p_i.
\end{equation}
We now compute its partial derivatives
\begin{equation}
\begin{array}{rl}
\dd{p_k}\sum_{i=1}^n p_i(2s_{i-1}-1+p_i)^2&= \dd{p_k}
p_k(2s_{k-1}-1+p_k)^2+ \dd{p_k}\sum_{i=k+1}^n
p_i(2s_{i-1}-1+p_i)^2=\\
&=(2s_{k-1}-1+p_k)^2+p_k\cdot 2(2s_{k-1}-1+p_k)+\\
&\hspace{.3cm}+\sum_{i=k+1}^n p_i\cdot 4(2s_{i-1}-1+p_i)=\\
&=(2s_{k-1}-1+p_k)^2-2p_k\cdot (2s_{k-1}-1+p_k)+\\
&\hspace{.3cm}+4\sum_{i=k}^n p_i(2s_{i-1}-1+p_i)=\\
&=( 2s_{k-1}-1 )^2-p_k^2+ 4\sum_{i=k}^n p_i(2s_{i-1}-1+p_i)=\\
&=( 2s_{n}-1 )^2-p_k^2=1-p_k^2\\
\end{array}
\end{equation}
from which
\begin{equation}
\dd{p_k} F_\lambda(p_1,\ldots,p_k)= \ln p_k -\rho p_k^2
+1+\rho-\lambda.
\end{equation}
If all partial derivatives are zero then $\ln p_k-\rho p_k^2=C$ for all
$1\leq k\leq n$ and some $C\in\R$. Therefore we may restrict the search of
maximum of $F$ on the set of probability weights to the subset where $\ln
p_k-\rho p_k^2=C$. On that set function $\tilde{F}$ is equal to
\begin{equation}
\tilde{F}(p_1,\ldots,p_n)=\rho\sum_{i=1}^n p_i(2s_{i-1}-1+p_i)^2+\sum_{i=1}^n
p_i(C+\rho p_i^2)
\end{equation}
and the derivative of this function is equal to
\begin{equation}
\dd{p_k}\tilde{F}(p_1,\ldots,p_n)=\rho(1-p_k^2)+3\rho p_k^2+C.
\end{equation}
If we apply Lagrange multipliers to the function $\tilde{F}$ then we obtain a
condition of the form $p_k^2+\tilde{C}=0$ for all $k\leq n$ and some
$\tilde{C}\in \R$. Since $p_k\in
[0,1]$ and $\sum_{i=1}^n p_k=1$ we get that
$p_1=p_2=\ldots=p_n=1/n$. Here we also have the maximum of
$F$ on the set of probability weights. 

For $p_1=p_2=\ldots=p_n=1/n$ we have
\begin{equation}
\begin{array}{rl}
\E\big(\E(X\cond\cA)\big)^2&= \re{n}\sum_{i=1}^n (2\frac{i-1}{n}-1+\re{n})^2=\\
&=\re{n^3}\sum_{i=1}^n \big(4i^2-4i(n+1)+(n+1)^2\big)=\\
&=\frac{4}{n^3}\frac{n(n+1)(2n+1)}{6}- \frac{4(n+1)}{n^3}\frac{n(n+1)}{2}+
\re{n^3}n(n+1)^2=\\
&=\frac{n+1}{n^2}\big(\frac{2}{3}(2n+1)-2(n+1)+(n+1)\big)=
\frac{(n+1)(n-1)}{3n^2}=\re{3}\big(1-\re{n^2}\big)\\
\end{array}
\end{equation}
and hence
\begin{equation}
f(n):=F(1/n,\ldots,1/n)=\frac{\rho}{3}\big(1-\re{n^2}\big)-\ln
n=\frac{\rho}{3}-\frac{\rho}{3n^2}-\ln n.
\end{equation}
The derivative of the function $\frac{\rho}{3x^2}+\ln x$ is
$-\frac{2\rho}{3x^3}+\re{x}$ so it is increasing for
$x<\sqrt{\frac{2\rho}{3}}=1/\sqrt{3\Sigma}$ and
decreasing for $x>1/\sqrt{3\Sigma}$ and therefore $f(n)$ achieves its maximum for
$\lfloor 1/\sqrt{3\Sigma}\rfloor$ or $\lceil 1/\sqrt{3\Sigma}\rceil$, where $\lfloor x
\rfloor$ and $\lceil x \rceil$ are the largest integer not greater than
$x$ and the smallest integer not less than $x$ respectively. Hence
setting $\Sigma\leq\re{12}$ leads to the overestimation of the number of
clusters.

\subsection*{\nameref*{subsec:exp}}
Let $p_{[ s,t ]}=P([s,t])$ and $e_{[ s,t ]}=\int_s^t x\d{P(x)}/p_{s,t}$. Then
\begin{equation}
p_{[ s,t ]} = e^{-s}-e^{-t},\quad 
e_{[ s,t ]} =
\frac{(1+s)e^{-s}-(1+t)e^{-t}}{e^{-s}-e^{-t}}=s+1-\frac{t-s}{e^{t-s}-1}
\end{equation}
Take any convex partition of $\R$, $\cA=\{[0,t_1),[t_1,t_2),\ldots\}$. Let
$p_i=p_{ [t_{i-1},t_i] }$ and $e_i=e_{ [t_{i-1},t_i] }$. Then
\begin{equation}
\Delta(\cA)=
\frac{R^2}{2} \sum_{i=1}^\infty p_i e_i^2 +
\sum_{i=1}^\infty p_i \ln p_i
\end{equation}
Assume that this sum is finite. Let $s\in [t_{n-1},t_n)$ and 
$\cA'=( \cA\sm\{[t_{n-1},t_n)\} )\cup\{[t_{n-1},s),[s,t_n)\}$ and let
$p_{n,1},e_{n,1},p_{n,2},e_{n,2}$ be defined as $p_i,e_i$, but for the intervals
$[t_{n-1},s)$ and $[s,t_n)$ respectively. Then
\begin{equation}\label{eq:deltadiff}
\Delta(\cA')-\Delta(\cA)=
\frac{R^2}{2}(p_{n,1}e_{n,1}^2+p_{n,2}e_{n,2}^2-p_n e_n^2)+
(p_{n,1}\ln p_{n,1} + p_{n,2}\ln p_{n,2} - p_n\ln p_n).
\end{equation}
Note that $p_{n,1}e_{n,1} + p_{n,2}e_{n,2}=p_n e_n$. Using this we can compute
that
\begin{equation}\label{eq:niceform}
p_{n,1}e_{n,1}^2+p_{n,2}e_{n,2}^2-p_n e_n^2=
\frac{p_{n,1}p_{n,2}}{p_n}|e_{n,1}-e_{n,2}|^2.
\end{equation}
Choose $s$ so that $p_{n,1}=p_{n,2}=\re{2}p_n$. Then it can be computed that $s=t_{n-1}+\ln
2-\ln(1-e^{t_{n-1}-t_n})$. Recall that
\begin{equation}
e_{n,2} = s+1 - \frac{t_n-s}{e^{t_n-s}-1}.
\end{equation}
Since $e_{n,1}\in [t_{n-1},s]$, $t_{n}-s>t_n-t_{n-1}-\ln 2$ and $x\mapsto
x/(e^x-1)$ is a decreasing function, we obtain that
\begin{equation}
\begin{split}
|e_{n,2}-e_{n,1}|&=e_{n,2}-e_{n,1}>e_{n,2}-s=\\
&=1- \frac{t_n-s}{e^{t_n-s}-1} > 1 - \frac{t_n-t_{n-1}-\ln 2}{e^{t_n-t_{n-1}-\ln 2} -1}.
\end{split}
\end{equation}
Hence if $t_n-t_{n-1}>3$ then $|e_{n,2}-e_{n,1}|>\re{2}$. By
\eqref{eq:deltadiff} and \eqref{eq:niceform} we get that for $t_n-t_{n-1}>3$
\begin{equation}
\Delta(\cA')-\Delta(\cA) =
\frac{R^2}{2}\cdot \re{4}p_n|e_{n,2}-e_{n,1}|^2 - p_n\ln 2>
p_n \zz( \frac{R^2}{32}-\ln 2 \zz).
\end{equation}
It means that for $\Sigma<(32\ln 2)^{-1}$ we increase the value of the function
$\Delta$ by dividing every segment of length larger than 3. As a result, no
finite convex partition can be a maximiser of the function $\Delta$ in this
case.

\subsection*{\nameref*{subsec:mnorm}}
Let $g_{a}(x) = \re{2\sqrt{2\pi}}(e^{-(x-a)^2/2} + e^{-(x+a)^2/2})$ be the
density of a mixture of two normal distributions, $\Normal(a,1)$ and
$\Normal(-a,1)$. We prove that for $a>1$ this distribution is bi-modal. 
It is easy to compute its derivatives:
\begin{equation}
\begin{split}
g'_a(x)&= -\re{2\sqrt{2\pi}}\big((x-a)e^{-(x-a)^2/2} + (x+a)e^{-(x+a)^2/2}\big),\\
g''_a(x)&= -\re{2\sqrt{2\pi}}\big(e^{-(x-a)^2/2} + e^{-(x+a)^2/2} -(x-a)^2
e^{-(x-a)^2/2} - (x+a)^2e^{-(x+a)^2/2}\big).
\end{split}
\end{equation}
Hence $g'_a(0)=0$ and $g''_a(0)=-\re{2\sqrt{2\pi}}e^{-a^2/2}(2-2a^2)>0$, which
means that 0 is a local minimum of $g_a$. Moreover the equation $g'_a(x)=0$ is
equivalent to
\begin{equation}\label{eq:bimodal}
U_a(x):=e^{2ax} - \frac{a+x}{a-x}
\end{equation}
Let us look for the solutions of this equation on $x\in (0,\infty)$. It is clear
that there are no solutions for $x\geq a$. 
It is straightforward to verify that $U_a(0)=0$, $U_a(a^-)=-\infty$.
Moreover $U'_a(x)=0$ is for $x\in (0,a)$ equivalent to $V_a(x):=(x-a)^2e^{2ax}=1$.
We have
\begin{equation}
V'_a(x)=2(x-a)e^{2ax}\big(1+a(x-a)\big)
\end{equation}
and hence $V'_a(x)=0$ has exactly one solution for $x\in (0,a)$ (which is
$\frac{a^2-1}{a}$). Since $V_a(0)=a^2>1$ and $V_a(a)=0$ we deduce that
$V_a(x)=1$ has exactly one solution in $(0,a)$, and so the equation
$U'_a(x)=0$. It is straightforward to verify that $U'_a(0)>0$ and therefore
$U_a(x)=0$ has exactly one solution for $x\in (0,a)$.

It follows that
$g'_a$ has exactly one zero on $(0,\infty)$; by symmetry there is also exactly one
zero on $(-\infty,0)$, so there are 3 zeros in total. Since we know that
$x=0$ is the local minimum of $g_a$ and $\lim_{x\to\pm\infty}g_a(x)=0$ it
follows that $g_a$ is bimodal.

\section*{Proof of \cref*{thm:nosmallm}}
Take $d=1$ and
$\alpha=\Tau=\Sigma=1$. Let $y_1,\ldots,y_n\in\R^d$. Take any partition
$\cJ$ of $[n]$. Let $J_n\in \cJ$ be the cluster containing $n$ and assume that
$|J_n|\geq 2$. Let $\cJ_{n,\{n\}}$ be obtained by creating a~singleton out of
$n$, i.e. $\cJ_{n,\{n\}}=\cJ\sm \{J_n\}\cup \{J_n\sm\{n\},\{n\}\}$.
By \eqref{eq:posterior} it is easy to show that the quotient
$\P(\cJ_{n,\{n\}}\cond y_1,\ldots,y_n)/\P(\cJ\cond y_1,\ldots,y_n)$
is equal to
\begin{equation}\label{eq:mcounterex}
h_{J_n}(y_1,\ldots,y_n)=\re{|J_n|-1}\sqrt{\frac{|J_n|+1}{2|J_n|}}
\exp\left\{\frac{y_n^2}{4}+\frac{(\sum y_{
J_n\sm\{n\}})^2}{2|J_n|}-\frac{(\sum
y_{J_n})^2}{2( |J_n|+1 )}\right\}.
\end{equation}
The exponent in the formula above is equal to
\begin{equation}\label{eq:mcounterex2}
y_n^2\frac{|J_n|-1}{4(|J_n|+1)}-y_n\frac{\sum
y_{J_n\sm\{n\}}}{|J_n|+1}+\frac{( \sum y_{J_n\sm\{n\}}
)^2}{2|J_n|(|J_n|+1)},
\end{equation}
which is a~convex quadratic function of $y_n$. Now, since
$|J_n|\geq 2$, it follows that
\begin{equation}\label{eq:mcounterex3}
\frac{|J_n|-1}{4(|J_n|+1)}\geq\re{12}
\quad\textrm{and}\quad
\Big|\frac{\sum y_{J_n\sm\{n\}}}{|J_n|+1}\Big|\leq |\overline{y_{J_n\sm\{n\}}}|.
\end{equation}
Now let $L=2\cdot 18^4$ and $\tilde{x}_m=18^m$.
We show that if
\begin{equation}\label{eq:singlecondini}
n\leq L^{m+1}, \quad 
y_n\geq \tilde{x}_m,\tag{$\star$}
\quad\textrm{and}\quad 
|y_1|,\ldots,|y_{n-1}|\leq\tilde{x}_{m-1} 
\end{equation}
then $h_{J_n}(y_1,\ldots,y_n)>1$ (regardless of $J_n$) and hence in
MAP partition for $[n]$ based on data $(y_i)_{i=1}^n$ singleton $\{n\}$ forms a~separate cluster.
Assume~\eqref{eq:singlecondini}. Note that if $n\leq L^{m+1}$ and 
$|y_1|,\ldots,|y_{n-1}|\leq\tilde{x}_{m-1}$ then by~\eqref{eq:mcounterex},
\eqref{eq:mcounterex2} and~\eqref{eq:mcounterex3} we obtain that
\begin{equation}
h_{J_n}(y_1,\ldots,y_n)\geq\re{L^{m+1}}\sqrt{\re{2}}
\exp\left\{\re{12}y_n^2-\tilde{x}_{m-1}y_n\right\}=:l(y_n).
\end{equation}
Now as we can easily compute zeros of quadratic function, $l(y_n)\geq 1$ is implied
by
\begin{equation}
y_n\geq
6\big(\tilde{x}_{m-1}+\sqrt{\tilde{x}_{m-1}^2+\re{3}[ (m+1)\ln L+(\ln 2)/2 ]}\big).
\end{equation}
It can be easily proved by induction that $3\tilde{x}_{m-1}^2>\re{3}[ (m+1)\ln
L+(\ln 2)/2 ]$ for $m\geq 2$ (note that the left-hand side is geometric
with respect to $m$, while the right-hand side is linear) and
therefore
\begin{equation}
6\big(\tilde{x}_{m-1}+\sqrt{\tilde{x}_{m-1}^2+\re{3}[ (m+1)\ln
L+(\ln 2)/2} ]\big)<18\tilde{x}_{m-1}=\tilde{x}_m
\end{equation}
and as $y_n\geq\tilde{x}_m$ we have that $h_{J_n}(y_1,\ldots,y_n)>1$.

Note that if $(y_n)_{n=1}^\infty$ is a~sequence whose terms belong to
$\{\tilde{x}_m\colon m\in\N\}$ then if for some $m\in\N$
\begin{equation}\label{eq:singlecond}
n\leq L^{m+1}, \quad 
y_n\geq \tilde{x}_m,\tag{$\star'$}
\quad\textrm{and}\quad 
y_1,\ldots,y_{n-1}<y_n
\end{equation}
then condition~\eqref{eq:singlecondini} holds with some $m'\geq m$ (the one that
satisfies $\tilde{x}_{m'}=y_n$). Indeed, if~\eqref{eq:singlecond} is satisfied
and $y_n=\tilde{x}_{m'}$ then as $y_1,\ldots,y_{n-1}<y_n$ we have
$y_1,\ldots,y_{n-1}\leq \tilde{x}_{m'-1}$, moreover $\tilde{x}_{m'}=y_n\geq
\tilde{x}_m$ and hence $m'\geq m$ and $n\leq L^{m+1}\leq L^{m'+1}$ and hence
\eqref{eq:singlecondini} is satisfied. 

\medskip
We now give an example of probability
weights $(p_m)_{m\in\N}$ such that the following probability distribution $P=\sum_{m=1}^\infty
p_m\delta_{\tilde{x}_m}$ has a~finite fourth moment and if
$(X_n)_{n=1}^\infty\iid P$ then~\eqref{eq:singlecond} happens almost surely 
infinitely many times. Let $q=L^{-1}$ and \mbox{$p_m=(1-q)q^{m-1}$}. It is
straightforward to check that in this case $P$ has finite fourth moment, as
\begin{equation}
\sum_{m=1}^\infty p_m\tilde{x}_m^4=(1-L^{-1})\sum_{m=1}^\infty
\frac{(18^m)^4}{(2\cdot 18^4)^{m-1}}=18^4(1-L^{-1})\sum_{m=1}^\infty
\re{2^{m-1}}<\infty.
\end{equation}
Now let $s_m=\sum_{i=1}^m p_i=1-q^{m}$. Then $s_m^{L^{m}}\to e^{-1}$. Let
$$n_m=\sum_{i=0}^m L^i=\frac{L^{m+1}-1}{L-1}<L^{m+1}$$ 
and
$A_m=\{\max_{n_{m-1}\leq i<n_m}X_i\geq \tilde{x}_m \}$. Then the probability of $A_m$ is equal
to $1-s_{m-1}^{L^m}$ which converges to $1-e^{-L}$. By the Borel-Cantelli
Lemma, it follows that almost surely infinitely many of the events
$A_m$ happens. Let $(x_n)_{n=1}^\infty$ be a~realisation of $(X_n)_{n=1}^\infty$
and let $(m_k)_{k=1}^\infty$ be an increasing sequence of all indices $m$ for which $A_{m}$ hold.
Now let $$\hat{n}_m=\min\{n_{m-1}\leq n<n_m\colon x_n=\max_{n_{m-1}\leq
i<n_m}x_i\}.$$ 
Then $x_{\hat{n}_{m_k}}\geq \tilde{x}_{m_k}$ for $k\in\N$. Let
$(k_i)_{i=1}^\infty$ be a
sequence such that $x_{\hat{n}_k}<x_{\hat{n}_{k_i}}$ for $k<k_i$ (such subsequence exists
since $\tilde{x}_m\to\infty$). Note that $x_{\hat{n}_{m_{k_i}}}\geq
\tilde{x}_{m_{k_i}}$, $\hat{n}_{m_{k_i}}<n_{m_{k_i}}<L^{m_{k_i}+1}$ and also
\begin{equation}
\textrm{for $l< \hat{n}_{m_{k_i}}$ we have } \left\{
\begin{array}{cl}
x_l< x_{\hat{n}_{m_{k_i}}} & \textrm{if $m(l)=m_{k_i}$,} \\
x_l\leq x_{\hat{n}_{ m(l) }}< x_{\hat{n}_{m_{k_i}}} & \textrm{if $m(l)=m_{k}$ for some
$k<k_i$,} \\
x_l< \tilde{x}_{m(l)}<\tilde{x}_{m_{k_i}}\leq
x_{\hat{n}_{m_{k_i}}} & \textrm{otherwise,}\\
\end{array}\right.
\end{equation}
where $m(l)=\min\{m\in\N\colon n_m>l\}$. From this it follows that for every $i\in\N$ condition
\eqref{eq:singlecond} is satisfied with $n=\hat{n}_{m_{k_i}}$ and $m=m_{k_i}$.
This proves that almost surely the
MAP partition creates a~new cluster out of a~new observation infinitely many times.

\section*{Proof of \cref*{thm:manyclust}}
Let $X_1,X_2,\ldots \iid P=\Exp(1)$ and $\hat{\cJ}_n$ be the MAP partition
computed on the basis of $X_1,\ldots,X_n$. We can assume that every value of
$X_i$ is unique and hence by \cref*{prop:MAPconv} we obtain that convex hulls of sets in
$\hat{\cJ}_n$ are pairwise disjoint.
Let $M_n=\max\{X_1,\ldots,X_n\}$ and $\cI_n$ be a
partition of $[0,M_n]$ into $|\hat{\cJ}_n|$ segments that `induce'
$\hat{\cJ}_n$, i.e. for every $J\in\hat{\cJ}_n$ there exist $I\in\cI_n$ such
that $\{x_j\colon j\in J\}\subset I$. 

Suppose, contrary to our claim, that
the sequence $|\hat{\cJ}_n|$ is bounded by some $K\in\N$. 
In order to use the results regarding
the behaviour of $\Delta$ function in the exponential case (\cref{subsec:exp}) we need to ensure
that there exist a sequence of segments $\overline{I}_n\in \cI_n$ and subsequence
$(n_k)_{k=1}^\infty$ such that
\begin{enumerate}[(i)]
\setlength{\itemsep}{1pt}
\setlength{\parskip}{0pt}
\setlength{\parsep}{0pt}
\item $\lim_{k\to\infty} |\overline{I}_{n_k}|=\infty$, 
where $|\cdot|$ is segment length,
\item $L:=\limsup_{k\to\infty} \inf \overline{I}_{n_k}<\infty$.
\end{enumerate}
We now construct such a sequence. Let $I^1_n$ be the sequence of the longest
segments in $\cI_n$ (i.e. $\diam I^1_n=\max_{I\in \cI_n}$). Since
almost surely $M_n\to\infty$ and the number of clusters within
the MAP partitions is assumed to be bounded, it follows that $|I^1_n|\to\infty$. If $\limsup_{n\to\infty}
\inf I^1_n<\infty$, set $\overline{I}_n=I^1_n$ and $n_k=k$.
Otherwise proceed inductively; having constructed the sequence $(I^i_n)_{n=1}^\infty$ and
subsequence $(n^i_k)_{k=1}^\infty$
such that $\lim_{k\to\infty} |I^i_{n^i_k}|=\infty$ do as follows: if $\limsup_{k\to\infty} \inf
I^i_{n^i_k}<\infty$ set $\overline{I}_n=I^k_n$ and $n_k=n^i_k$. If not, let
$I^{i+1}_n$ be the sequence of the longest segments to the left of $I^i_n$ in
$\cI_n$ (i.e. $|I^{i+1}_n|=\max\{|I|\colon
I\in \cI_n, \sup I\leq \inf I^i_n\}$). By the assumption about bounded
number of clusters we obtain that $\lim_{k\to\infty} |I^{i+1}_n|=\infty$. Note
that this procedure has to stop after at most $K$ iterations, because by
construction there are at most $K-i$ segments to the left of $I^i_n$.
Therefore requirement (ii) is bound to be finally satisfied.

Note that, because of (i) and (ii) we can deduce that $\liminf_{k\to\infty}
P(\overline{I}_{n_k})\geq P([L,\infty))=:p>0$. 

Let $\overline{J}_n$ be the cluster in $\hat{\cJ}_n$ induced by $\overline{I}_n$, i.e.
$\overline{J}_n=\{i\leq n\colon X_i\in \overline{I}_n\}$. Let
$\{\overline{I}^1_n,\overline{I}^2_n\}$ be a partition of $\overline{I}_n$ into
two equally probable segments, which induces partition
$\{\overline{J}^1_n,\overline{J}^2_n\}$ of $\overline{J}_n$.
Let $\tilde{\cJ}_n$ be obtained from
$\hat{\cJ}_n$ by replacing $\overline{J}_n$ by two sets $\overline{J}^{1}_n$
and $\overline{J}^{2}_n$. Then
\begin{equation}\label{eq:expxxx1}
\frac{\P(\tilde{\cJ}_n\cond\textbf{x})}{\P(\hat{\cJ}_n\cond\textbf{x})}
= C\frac{a_{n,1}!a_{n,2}!}{a_n!}\left( \frac{a_n}{a_{n,1} a_{n,2}} \right)^{3/2}\frac{
R_{a_n}}{R_{a_{n,1}}R_{a_{n,2}}}\cdot
\exp\left\{D_n\right\}^{1/2},
\end{equation}
where $a_n, a_{n,1},a_{n,2}$ are the sizes of
$\overline{J}_n,\overline{J}^1_n,\overline{J}^2_n$ respectively,
$R_n=\sqrt{R^2+U^2/n}$ and
\begin{equation}
D_n=a_{n,1}\znorm{R_{a_{n,1}}^{-1}R^2\overline{\bx_{\overline{J}_n^{1}}}}^2+
a_{n,2}\znorm{R_{a_{n,2}}^{-1}R^2
\overline{\bx_{\overline{J}_n^{2}}}}^2-a_n\znorm{R_{a_n}^{-1}R^2\overline{\bx_{\overline{J}_n^M}}}^2.
\end{equation}
By \cref*{lem:gliwcan} we have 
\begin{equation}
a_n/n-P(\overline{I}_n)\to 0,\quad
a_{n,1}/n-P(\overline{I}^1_n)\to 0,\quad
a_{n,2}/n-P(\overline{I}^2_n)\to 0
\end{equation}
Since $\liminf_{k\to\infty}
P(\overline{I}_{n_k})\geq p>0$ it follows that $a_{n_k}\to\infty$ as
$k\to\infty$. By Stirling formula and the Strong Law of Large Numbers
\begin{equation}
\begin{split}
\sqrt[n_k]{ \frac{\P(\tilde{\cJ}_{n_k}\cond\textbf{x})}{\P(\hat{\cJ}_{n_k}\cond\textbf{x})} }
&\approx
\sqrt[n_k]{\frac{a_{n_k,1}^{a_{n_k,1}}a_{n_k,2}^{a_{n_k,2}}}{a_{n_k}^{a_{n_k}}} }
\exp\left\{\frac{R^2}{2}( \frac{a_{{n_k},1}}{{n_k}}\overline{\bx_{\overline{J}_{n_k}^1}}^2+
\frac{a_{{n_k},2}}{{n_k}}\overline{\bx_{\overline{J}_{n_k}^2}}^2-\frac{a_{n_k}}{{n_k}}\overline{\bx_{\overline{J}_{n_k}}}^2
)\right\}=\\
&= \left(\frac{a_{n_k,1}}{a_{n_k}}\right)^{a_{n_k,1}/n_k}
\left(\frac{a_{n_k,2}}{a_{n_k}}\right)^{a_{n_k,2}/n_k}
\exp\left\{\frac{R^2}{2}( \frac{a_{{n_k},1}}{{n_k}}\overline{\bx_{\overline{J}_{n_k}^1}}^2+
\frac{a_{{n_k},2}}{{n_k}}\overline{\bx_{\overline{J}_{n_k}^2}}^2-\frac{a_{n_k}}{{n_k}}\overline{\bx_{\overline{J}_{n_k}}}^2
)\right\}\approx\\
&\approx 2^{-a_{n_k}/{n_k}} 
\exp\left\{\frac{R^2}{2}( \frac{a_{{n_k},1}}{{n_k}}\overline{\bx_{\overline{J}_{n_k}^1}}^2+
\frac{a_{{n_k},2}}{{n_k}}\overline{\bx_{\overline{J}_{n_k}^2}}^2-\frac{a_{n_k}}{{n_k}}\overline{\bx_{\overline{J}_{n_k}}}^2
)\right\}\approx\\
&\approx\exp\left\{\frac{R^2}{8}\frac{a_{n_k}}{{n_k}}( \overline{\bx_{\overline{J}_{n_k}^1}}-
\overline{\bx_{\overline{J}_{n_k}^2}})^2-\frac{a_{n_k}}{{n_k}}\ln 2\right\},
\end{split}
\end{equation}
By \cref*{cor:EapproxWn} (with $\delta=p/2$) we know that
$\overline{\bx_{\overline{J}_{n_k}^{1}}}$ and
$\overline{\bx_{\overline{J}_{n_k}^{2}}}$ approximate $\E(X\cond X\in
\overline{I}_n^1)$ and $\E(X\cond X\in \overline{I}_n^2)$. Since $\overline{I}_n$ becomes
arbitrarily large, its length finally exceeds 3 and previous considerations
(cf. \cref*{subsec:exp}) lead to the conclusion, that $\P(\hat{\cJ}_{n_k}\cond
\bx_{1:{n_k}})<\P(\tilde{\cJ}_{n_k}\cond
\bx_{1:{n_k}})$ for large enough $k$. This is a contradiction that proves our
assertion.

\section*{Proof of \cref*{lem:induced}}
Note that we may assume that $P(A)>0$ for $A\in\cA$. Indeed, if $P(A)=0$ then
$\Delta(\{A\})=0$ (by a natural convention that $0\ln 0=0$) and on the other
hand if $X_1,X_2,\ldots\iid P$ then almost surely $X_i\notin A$ for $i\in\N$.

We abuse the notation slightly and denote $p_{J^A_n}=|J^A_n|/n$ for $A\in\cA$. By the law of
large numbers the sequence $(X_n)_{n=1}^\infty$ almost surely satisfies
$p_{J^A_n}\to p_A>0$. By Stirling formula
\begin{equation}\label{eq:stirling1}
\prod_{J\in\cJ^A_n}(np_J)!\approx \prod_{J\in\cJ^A_n}\left(\frac{np_J}{e}\right)^{np_J}\sqrt{2\pi
np_J}
=\sqrt{2\pi
n}^{|\cJ^A_n|}\sqrt{\prod_{J\in\cJ^A_n}p_J}\cdot\big(\frac{n}{e}\prod_{J\in\cJ^A_n}p_J^{p_J}\big)^n
\end{equation}
from which it follows by the Strong Law of Large Numbers that
$ \sqrt[n]{\prod_{J\in\cJ^A_n}(np_J)!}\approx \frac{n}{e}\prod_{J\in\cJ^A_n}p_J^{p_J}
\approx \frac{n}{e}\prod_{J\in\cA}p_A^{p_A}$.
Note that since $J^{\cA}_n$ has at most $|\cA|$ elements,
\begin{equation}\label{eq:convto1}
\lim_{n\to\infty}\sqrt[n]{C^{|\cJ^{\cA}_n|}}=1
\quad\textrm{and}\quad
\lim_{n\to\infty}\sqrt[n]{\prod_{J\in\cJ^{\cA}_n}|J|^{(d+2)/2}\det R_{|J|}}=1.
\end{equation}
It follows from the law of large numbers that
$\overline{\bX_{J^A_n}}\to\E (X\cond X\in A)$ for $A\in\cA$ almost surely.
It follows that
\begin{equation}\label{eq:convVar}
\re{n}\sum_{J\in\cJ^{\cA}_n}|J|\znorm{ R_{|J|}^{-1}R^2 \overline{\bX_J}}^2
\approx \sum_{J\in\cJ^{\cA}_n}p_J\znorm{ R~\overline{\bX_J}}^2
\approx\sum_{A\in\cA}p_A\znorm{R\E (X\cond X\in A)}^2.
\end{equation}
Applying~\eqref{eq:stirling1},~\eqref{eq:convto1} and~\eqref{eq:convVar}
together with~\eqref{eq:DeltaDef} to
the formula~\eqref{eq:posterior} for $\cJ^\cA$ 
completes the proof of the Lemma.

\section*{Proof of \cref*{lem:approxAn}}
From \cref*{cor:max} (a) we know $\min\{np_J\colon J\in \hat{\cJ}_n\}\to\infty$.
By applying Stirling formula to each factor $(np_J)!$ and taking into account
that by \cref*{cor:max} (b) the number of factors is bounded, we obtain that
\begin{equation}
\prod_{J\in\hat{\cJ}_n}(np_J)!\approx
\prod_{J\in\hat{\cJ}_n}\left(\frac{np_J}{e}\right)^{np_J}\sqrt{2\pi
np_J}=
\big(\frac{n}{e}\big)^n\sqrt{2\pi
n}^{|\hat{\cJ}_n|-1}\sqrt{\prod_{J\in\hat{\cJ}_n}p_J}\cdot\left(\prod_{J\in\hat{\cJ}_n}p_J^{p_J}\right)^{n}.
\end{equation}
By definition the elements of $\hat{\cA}_n$ are convex and hence by
\cref*{lem:gliwcan} the frequencies $p_J$ for $J\in\cJ_n$ 
approximate the respective probabilities of sets in $\hat{\cA}_n$ uniformly. Hence, as
$(|\hat{\cJ}_n|)_{n=1}^\infty$ is bounded almost surely, it follows that
$
\sqrt[n]{\prod_{J\in\hat{\cJ}_n}(np_J)!}\approx
\frac{n}{e}\prod_{J\in\hat{\cJ}_n}p_J^{p_J} \approx \frac{n}{e}\prod_{A\in\hat{\cA}_n}p_A^{p_A}.
$
By applying a~similar argument to the remaining part of formula
\eqref{eq:posterior}, the result follows by \cref*{cor:EapproxWn}, which is an
easy consequence of \cref{lem:Eapprox} (here we also
use \cref*{cor:max} (a)).

\begin{lem}\label{lem:Eapprox}
If $P$ satisfies~\eqref{eq:gcc} and for $X\sim P$ we have $\E \norm{X}^2<\infty$ then
$P$ satisfies
\begin{equation}\label{eq:gccE}
\lim_{n\to\infty}\sup_{C\in\cK} \znorm{\E_n X\1_{X\in C}-\E X\1_{X\in C}}=0
\quad\textrm{almost surely.}\tag{$\ast\ast$}
\end{equation}
where $\E_n f(X)=\int_{\cX}f(X)\d{P_n}=\re{n}\sum_{i=1}^n f(X_i)$.
\end{lem}
\begin{proof}
Let $x^{(i)}$ ($i\leq d$) be the $i$-th coordinate of vector $x$. We now prove
that for every $r>0$
\begin{equation}\label{eq:Eapp}
\lim_{n\to\infty}\sup_{C\in\cK\cap [-r,r]^d} \big|\E_n X^{(1)}\1_{X\in C} -\E
X^{(1)}\1_{X\in C}\big|=0.  
\end{equation}
Fix $r>0$ and $C\in\cK\cap [-r,r]^d$. For $m\in\N$ and
$-m\leq k\leq m-1$ let $C^m_k=C\cap [rk/m,r(k+1)/m)\times\R^{d-1}$. Then
\begin{equation}\label{eq:Eapp1}
\Big|\E X^{(1)}\1_{X\in C}-\sum_{k=-m}^{m-1} r\frac{k}{m}
P(C^m_k)\Big|\leq\frac{r}{m}P(C)\leq \frac{r}{m}.
\end{equation}
It follows from the same reasoning
\begin{equation}\label{eq:Eapp2}
\Big|\E_n X^{(1)}\1_{X\in C}-\sum_{k=-m}^{m-1} r\frac{k}{m}
P_n(C^m_k)\Big|\leq\frac{r}{m}\quad\textrm{for every $n\in\N$}.
\end{equation}
Now choose $\eps>0$ and $m>r/\eps$. Note that $C^m_k$ are convex sets (as
intersections of two convex sets) and hence by~\eqref{eq:gcc} we may choose
$N$ so that for $n>N$ and any convex $C'$ we have that
$|P_n(C')-P(C')|<\eps/(2m)$ and
hence
\begin{equation}\label{eq:Eapp3}
\Big|\sum_{k=-m}^{m-1} r\frac{k}{m} P(C^m_k)-
\sum_{k=-m}^{m-1} r\frac{k}{m} P_n(C^m_k)\Big|\leq
\sum_{k=-m}^{m-1} r\frac{|k|}{m} |P(C^m_k)-P_n(C^m_k)|< r\eps.
\end{equation}
By combining~\eqref{eq:Eapp1},~\eqref{eq:Eapp2} and~\eqref{eq:Eapp3} we obtain
that $|\E_n X^{(1)}\1_{X\in C}-\E X^{(1)}\1_{X\in C}|< (2+r)\eps$
for $n>N$ and since the choice of $N$ does not depend on $C$, 
\eqref{eq:Eapp} follows.

We now prove that almost surely
\begin{equation}\label{eq:EappNotRestricted}
\lim_{n\to\infty}\sup_{C\in\cK} 
\big|\E_n X^{(1)}\1_{X\in C} -\E X^{(1)}\1_{X\in C}\big|=0.  
\end{equation}
The same result for the remaining coordinates of~\eqref{eq:gccE} follows in the
same way, from which follows the statement of the Lemma.  
Note that the function $r\mapsto \E |X^{(1)}|\1_{X\notin [-r,r]^d}$
is decreasing to 0 as $r$ goes to infinity. By the Strong Law of Large Numbers almost
surely $\lim_{n\to\infty} \E_n |X^{(1)}|\1_{X\notin [-K,K]^d}=\E
|X^{(1)}|\1_{X\notin [-K,K]^d}$ for every $K\in \N$.

Fix $C\in\cK$ and $\eps>0$. Since $\lim_{K\to\infty} \E |X^{(1)}|\1_{X\notin [-K,K]^d}=0$ 
it follows that there exist $K\in \N$ such that $\E
|X^{(1)}|\1_{X\notin [-K,K]^d}<\eps$ and $\lim_{n\to\infty} \E_n
|X^{(1)}|\1_{X\notin [-K,K]^d}<\eps$. The latter means that there exist
$n_1$ such that $\E_n |X^{(1)}|\1_{X\notin [-K,K]^d}<\eps$ for every $n>n_1$. By
\eqref{eq:Eapp} there exist $n_2\in\N$ such that for every $n>n_2$ 
\begin{equation}
\big|\E_n X^{(1)}\1_{X\in C\cap [-K,K]^d} -\E X^{(1)}\1_{X\in C\cap[-K,K]^d}\big|<\eps.  
\end{equation}
Therefore for $n>\max\{n_1,n_2\}$ we get
\begin{equation}
\begin{split}
\big|\E_n X^{(1)}\1_{X\in C} -\E X^{(1)}\1_{X\in C}\big|&<  
\big|\E_n X^{(1)}\1_{X\in C\cap [-K,K]^d} -\E X^{(1)}\1_{X\in
C\cap[-K,K]^d}\big|+\\
&\quad+\E_n |X^{(1)}|\1_{X\notin [-K,K]^d}+ \E |X^{(1)}|\1_{X\notin [-K,K]^d}<3\eps
\end{split}
\end{equation}
Because $n_1,n_2$ do not depend on $C$, \eqref{eq:EappNotRestricted} follows,
which finishes the proof of the Lemma.
\end{proof}

\begin{cor}\label{cor:EapproxWn}
If $P$ satisfies~\eqref{eq:gcc} and for $X\sim P$ we have $\E \norm{X}^2<\infty$ then for every $\delta>0$ we have
\begin{equation}
\lim_{n\to\infty}\sup_{\substack{C\in\cK\\ P(C)>\delta}} \znorm{\E_n ( X\cond X\in C )-\E ( X\cond X\in C )}=0
\quad\textrm{almost surely.}
\end{equation}
\end{cor}
\begin{proof}
This is a straightforward consequence of \cref{lem:Eapprox} and the definition
$\E_n ( X\cond X\in C )=\E_n X\1_{X\in C} /P_n(C)$. \end{proof}

\section*{Proof of \cref*{thm:famconv}}
Assume that $(\cF,d)$ is a~(pseudo)metric space. We prove that $(F_K(\cF),\bar{d})$ is also a
(pseudo)metric space. Take any $\cA=\{A^{(1)},\ldots,A^{(k)}\}\in F_K(\cF)$ and
$\cB=\{B^{(1)},\ldots,B^{(l)}\}\in F_K(\cF)$. By definition
\begin{equation}
\bar{d}(\cA,\cB)=\min_{\sigma\in\Sigma_K} \max_{i\leq K}
d(A^{(i)},B^{(\sigma(i))})\geq 0,
\end{equation}
since $d(A^{(i)}, B^{(j)})\geq 0$ for any $i,j\leq K$ (as in the definition we
assume that $A^{(i)}=\emptyset$ and $B^{(j)}=\emptyset$ for $i>k$ or $j>l$
respectively).
Let $\cC=\{C^{(1)},\ldots,C^{(l)}\}\in F_K(\cF)$ and let $\sigma_1$, $\sigma_2$
and $\sigma_3$ be permutations of $[K]$ that satisfy
\begin{equation}
\overline{d}(\cA,\cB)=\max_{i\leq K}
d(A^{(i)},B^{(\sigma_1(i))})\quad\textrm{and}\quad
\overline{d}(\cB,\cC)=\max_{i\leq K} d(B^{(i)},C^{(\sigma_2(i))})
\end{equation}
Note that 
$d(A^{(i)},B^{(\sigma_1(i))})+d(B^{(\sigma_1(i))},C^{(\sigma_2(\sigma_1(i)))})
\geq d(A^{(i)},C^{(\sigma_2(\sigma_1(i)))})$ and hence
\begin{equation}
\begin{split}
\overline{d}(\cA,\cB) + \overline{d}(\cB,\cC) &=
\max_{i\leq K} d(A^{(i)},B^{(\sigma_1(i))})+
\max_{i\leq K} d(B^{(\sigma_1(i))},C^{(\sigma_2(\sigma_1(i)))})\geq\\
&\geq\max_{i\leq K}
\big(d(A^{(i)},B^{(\sigma_1(i))})+d(B^{(\sigma_1(i))},C^{(\sigma_2(\sigma_1(i)))})\big)\geq\\
&\geq\max_{i\leq K}d(A^{(i)},C^{(\sigma_2\circ\sigma_1(i))})\geq
\overline{d}(\cA,\cC)
\end{split}
\end{equation}
and the triangle inequality follows. This means that $\overline{d}$ is a pseudometric on $\cF_K$.

Now assume that $(\cF,d)$ is finitely compact. Let $(\cA_n)_{n=1}^\infty$ be a~sequence in $F_K(\cF)$ and let 
$\cA_n=\{A_n^{(1)},A_n^{(2)},\ldots,A_n^{(k_n)}\}$. As the sequence
$(k_n)_{n=1}^\infty$ is bounded by $K$ we may choose a~subsequence $\cA_{n_k}$
and $\tilde{K}\in\N$ such that $|\cA_{n_k}|=\tilde{K}$ for every $k\in\N$. Consider the
sequence $(A_{n_k}^{(1)})_{k=1}^\infty$. This sequence is bounded (as
$(\cA_n)_{n=1}^\infty$ is bounded). Therefore it has a~subsequence
$(A_{n_{k_l}}^{ (1) })_{l=1}^\infty$ converging in $d$ to $A^{(1)}\in\cF$. 
Now we consider $(A_{n_{k_l}}^{(2)})_{l=1}^\infty$ and again we choose
a subsequence $(A_{n_{k_{l_m}}})_{m=1}^\infty$ converging in $d$ to
$A^{(2)}\in\cF$. By
iterating this procedure $\tilde{K}$ times we obtain a~family
$\hat{A}=\{A^{(1)},\ldots,A^{(\tilde{K})}\}$ of `limiting' sets. It is easy to
verify that the final subsequence of $(\cA_n)_{n=1}^\infty$ converges in
$\bar{d}$ to $\hat{A}$, which finishes the proof.

\section*{Proof of \cref*{thm:MAPmax}}
Take any $\cE=\{E^{(1)},\ldots,E^{(\tilde{K})}\}\in\mathbf{E}$ and assume that it
is a~limit of $( \hat{\cA}_{n_k} )_{k=1}^\infty$ in~$\overline{d_H}$. By \cref*{thm:dPcomp} the sequence $(
\hat{\cA}_{n_k} )_{k=1}^\infty$ converges to $\cE$ also in $\overline{d_P}$.
Since for every $k\in\N$ every pair of sets in the family $\hat{\cA}_{n_k}$ has at most one point in
common (\cref*{prop:MAPconv}) then by the continuity of $P$ with respect
to the Lebesgue measure every pair of sets within $\hat{\cA}_{n_k}$ has an
intersection of $P$ measure 0. Therefore by the continuity of the intersection with respect to $d_P$ (\cite{bib:doob},
Chapter III, Formula (13.3)) we get that $P(E^{(i)}\cap E^{(j)})=0$ for $1\leq i<j\leq
\tilde{K}$. 

To prove that $\cE$ is a~$P$-partition it is left to
show that $P(\bigcup \cE)=1$ (we denote $\bigcup \cE=\bigcup_{E\in\cE} E$). Suppose this is not the case. It means that $E_0=\R^d\sm
\bigcup\cE$ is an open set with positive probability. Therefore it includes
a~ball $B'$ of positive probability. Since $B'$ is a convex set, we get
$p_{J^{B'}_n}\to p_{B'}>0$ and therefore there exist $n'\in\N$ such that
$X_{n'}\in B'$. This is not possible, since $X_{n'}\in \bigcup \hat{\cA}_{n}$
for every $n\geq n'$ and therefore $X_{n'}\in \bigcup \cE$, which is a~contradiction. 

By \cref*{lem:approxAn} and the continuity of $\Delta$
with respect to the metric $\overline{d_P}$ we obtain:
\begin{equation}\label{eq:Qx1}
\sqrt[n]{Q_{\bX_{1:n}}(\hat{\cJ}_n)}\approx \exp(\Delta(\hat{\cA}_n))\approx
\exp(\Delta(\cE)).
\end{equation}

Now take any finite $P$-partition $\cA$. We can assume that each $X_n$ belongs to exactly one of the sets in $\cA$,
$p_{\cJ^A_n}\to p_A$ and $\overline{\bX_{\cJ^A_n}}\to \E(X\cond X\in A)$ for $A\in\cA$ (it just requires adding a
countable number of conditions on the
infinite iid sequence with distribution $P$, each of which is satisfied almost surely). By definition of $\hat{\cJ}_n$ and
\cref{lem:induced} we get
\begin{equation}\label{eq:Qx2}
\sqrt[n]{Q_{\bX_{1:n}}(\hat{\cJ}_n)}\geq \sqrt[n]{Q_{\bX_{1:n}}(\cJ^{\cA}_n)}
\approx \exp(\Delta(\cA)).
\end{equation}

Equations~\eqref{eq:Qx1} and~\eqref{eq:Qx2} together give us $\Delta(\cE)\geq
\Delta(\cA)$ which proves that $\cE$ is a~finite partition that maximises
$\Delta$.

\section*{Proofs for \Asection{sec:discuss}}

\begin{lem}\label{lem:consistinput}
Let $\alpha>0$, $G_0=\Normal(0,\Tau)$, $F_\theta=\Normal(\theta,\Sigma)$ for $\theta\in\R^d$.
Let $X_1,X_2,\ldots$ be an infinite sample from the Gaussian DPMM, i.e. a
sequence of random variables, obtained by the following construction defined by
\cref{eq:DPMgen}.
Then almost surely $\limsup_{n\to\infty} \re{n}\sum_{i=1}^n \norm{X_i}^2<\infty$.
\end{lem}
\begin{proof}
It is a well known fact
(\cite{bib:sethuraman1994constructive}) that an infinite sample from the DPMM may be
performed by the following procedure
\begin{equation}
\begin{split}
\bm{p}=(p_1,p_2,\ldots)&\sim SB(\alpha)\\
\bm{\theta}=(\theta_1,\theta_2,\ldots) &\iid \Normal(0,\Tau)\\
X_1,X_2,\ldots\cond \bm{p},\bm{\theta}&\iid \sum_{i=1}^\infty
p_i\Normal(\theta_i,\Sigma)
\end{split}
\end{equation}
where $SB(\alpha)$ is the so called \emph{stick-breaking construction}
(i.e. $p_n=V_n\prod_{i=1}^{n-1}(1-V_i)$, where $V_1,V_2,\ldots\iid
\Beta(1,\alpha)$). Therefore 
\begin{equation}\label{eq:DPMMexp}
\E(\norm{X_1}^4\cond\bm{p},\bm{\theta})) = 
\sum_{n=1}^\infty p_n Q(\theta_n)
\end{equation}
where by definition $Q(\theta)=\E\norm{X}^4$ when
$X\sim\Normal(\theta,\Sigma)$. Note that $Q$ is a multivariate polynomial in the
coefficients of $\theta$; it is given by the formula
\begin{equation}
Q(\theta)=\sum_{k=1}^d (\theta_{(k)}^4+6\sigma_{k,k}\theta_{(k)}^2
+3\sigma_{k,k}^2) +
\sum_{k\neq l}(
2\sigma_{k,l} + 4\sigma_{k,l}\theta_{(k)} \theta_{(l)}+
\theta_{(k)}^2\theta_{(l)}^2+
\sigma_{k,k}\sigma_{l,l} + \sigma_{k,k}\theta_{(l)}^2 +
\sigma_{l,l}\theta_{(k)}^2
)
\end{equation}
where $\theta_{(k)}$ is the $k$-th coefficient of $\theta$ and
$[\sigma_{k,l}]_{k,l\leq d}=\Sigma$. We show that for $\sum_{i=1}^n p_i Q(\theta_i)$
is almost surely finite. Note that
\begin{equation}\label{eq:DPMMinttheta}
\P(\sum_{n=1}^\infty p_i Q(\theta_i)<\infty)=
\int_{\overline{\bm{p}}\in\R^\infty} \P(\sum_{n=1}^\infty p_i Q(\theta_i)<\infty\cond
\bm{p}=\overline{\bm{p}}) \d{SB(\overline{\bm{p}})}
\end{equation}

Note that given $\bm{p}=\overline{\bm{p}}=(\overline{p_n})_{n=1}^\infty$, where
$\sum_{n=1}^\infty \overline{p_n}=1$ and $\overline{p_n}\in[0,1]$, the series $\sum_{n=1}^\infty p_n
Q(\theta_n)$ is a sequence of independent random variables, which is almost
surely bounded by the Kolmogorov's Two Series
Theorem (\cite{bib:durrett2010probability}, Theorem 2.5.3). Indeed, when
$\theta\sim \Normal(0,\Tau)$ then all mixed moments of $\theta$ are finite
(i.e. $\E \prod_{i=1}^d \theta_{(i)}^{w_i}<\infty$ for every
$w_1,\ldots,w_d\in\N$) and therefore $\E Q(\theta)<\infty$ and $\Var
Q(\theta)<\infty$, so
\begin{equation}
\begin{split}
\sum_{n=1}^\infty \E_{\bm{p}=\overline{\bm{p}}} \overline{p_n} Q(\theta_n) &=
\E_{\bm{p}=\overline{\bm{p}}} Q(\theta) \sum_{n=1}^\infty \overline{p_n}=
\E Q(\theta)<\infty,\\
\sum_{n=1}^\infty \Var_{\bm{p}=\overline{\bm{p}}} \overline{p_n} Q(\theta_n) &=
\Var_{\bm{p}=\overline{\bm{p}}} Q(\theta) \sum_{n=1}^\infty \overline{p_n}^2<
\Var Q(\theta)<\infty
\end{split}
\end{equation}
Therefore $\P(\sum_{n=1}^\infty p_i Q(\theta_i)<\infty\cond
\bm{p}=\overline{\bm{p}})=1$ and \eqref{eq:DPMMinttheta} implies 
$\P(\sum_{n=1}^\infty p_i Q(\theta_i)<\infty)=1$.
Finally, conditioned on $\bm{p}=\overline{\bm{p}}$ and
$\bm{\theta}=\overline{\bm{\theta}}$, the sequence $X_1,X_2,\ldots$ is a
sequence of independent random variables and therefore if $\sum_{n=1}^\infty p_i
Q(\theta_i)<\infty$ then by the Strong Law of
Large Numbers and \eqref{eq:DPMMexp}, we have 
\begin{equation}
\lim_{n\to\infty} \re{n}\sum_{i=1}^n \norm{X_i}^2
\xlongequal{\P_{\overline{\bm{p}},\overline{\bm{\theta}}}\ a.s.}
\E(\norm{X_1}^2\cond\bm{p}=\overline{\bm{p}},\bm{\theta}=\overline{\bm{\theta}}))<
\E(\norm{X_1}^4\cond\bm{p}=\overline{\bm{p}},\bm{\theta}=\overline{\bm{\theta}}))^{\re{2}}<\infty.
\end{equation}
This means that
\begin{equation}
\begin{split}
\P(\limsup_{n\to\infty} \re{n}\sum_{i=1}^n \norm{X_i}^2 <\infty)&=\\
&\hspace*{-3cm}=\int_{\overline{\bm{p}}\in\R^\infty,\overline{\bm{\theta}}\in\R^\infty}
\P(\limsup_{n\to\infty} \re{n}\sum_{i=1}^n \norm{X_i}^2 <\infty \cond
\bm{p}=\overline{\bm{p}},\bm{\theta}=\overline{\bm{\theta}})
\d{SB(\overline{\bm{p}})}\d{G_0^\infty(\overline{\bm{\theta}})}\geq\\
&\hspace*{-3cm}\geq\int_{\overline{\bm{p}}\in\R^\infty,\overline{\bm{\theta}}\in\R^\infty}
\P(\sum_{n=1}^\infty \overline{p_n}Q(\overline{\theta_n})<\infty)
\d{SB(\overline{\bm{p}})}\d{G_0^\infty(\overline{\bm{\theta}})}=
\P(\sum_{n=1}^\infty p_nQ(\theta_n)<\infty)=1
\end{split}
\end{equation}
and the proof follows.
\end{proof}
\section*{A comment on extending the results to the finite Dirichlet mixture
models}
The \emph{finite Dirichlet zero-mean Gaussian mixture model} is a model of the form
\begin{equation}
\begin{split}
\bm{p}=(p_1,\ldots,p_K)&\sim\Dir(\alpha,\alpha,\ldots,\alpha)\\
\bm{\theta}=(\theta_1,\ldots,\theta_K) &\iid \Normal(0,T)\\
\bm{x}=(x_1,\ldots,x_n)\cond \bm{\theta},\bm{p}&\iid \sum_{k=1}^K p_k\Normal(\theta_k,\Sigma)
\end{split}
\end{equation}
Given the vector $\bm{p}$ of probabilities of belonging to respective clusters,
the marginal probability on partition $\cJ$ of indices is
\begin{equation}
\P(\cJ\cond\bm{p})=\sum_{\tau\colon\cJ\stackrel{1-1}{\to}[K]}
\prod_{J\in\cJ} p_{\tau(J)}^{|J|},
\end{equation}
where the sum is over all injective functions $\tau$ from $\cJ$ to $[K]$. By properties
of the Dirichlet distribution we have that for every
$\tau\colon\cJ\stackrel{1-1}{\to}[K]$
\begin{equation}
\E\prod_{J\in\cJ} p_{\tau(J)}^{|J|} = \frac{\Gamma(K\alpha)}{\Gamma(K\alpha+n)}
\prod_{J\in\cJ} \frac{\Gamma(\alpha+|J|)}{\Gamma(\alpha)}=
\re{(K\alpha)^{(n)}}\prod_{J\in\cJ} \alpha^{(|J|)},
\end{equation}
where $a^{(b)}=\frac{\Gamma(a+b)}{\Gamma(a)}=a(a+1)\ldots (a+b-1)$. There are
$K^{[|\cJ|]}=K(K-1)\ldots (K-|\cJ|+1)$ injections from $\cJ$ to $[K]$, therefore
\begin{equation}
\P(\cJ)=\E\P(\cJ\cond\bm{p})
=\frac{K^{[|\cJ|]}}{(K\alpha)^{(n)}}\prod_{J\in\cJ} \alpha^{(|J|)}
=\frac{K^{[|\cJ|]}\alpha^{|\cJ|}}{(K\alpha)^{(n)}}\prod_{J\in\cJ} (1+ \alpha )^{(|J|-1)}.
\end{equation}
We will call the resulting distribution the \emph{finite Chinese Restaurant
Process} and write $\cJ\sim \CRP^K(\alpha)_n$. It is easy to see that
$\CRP^K(\frac{\alpha}{K})_n$ converges to $\CRP(\alpha)_n$, which is a
well-known fact (see \cite{bib:neal}).
We can now consider a \emph{finite CRP based zero-mean Gaussian model}, defined by the following
scheme
\begin{equation}
\begin{array}{rcll}
\cJ&\sim&\CRP^K(\alpha)_n&\\
\bm{\theta}=(\theta_J)_{J\in \cJ}\cond \cJ&\iid&\Normal(0,T)&\\
\bx_J=(x_j)_{j\in J}\cond \cJ,\bm{\theta}&\iid&\Normal(\theta_J,\Sigma)&\textrm{for $J\in \cJ$}.
\end{array}
\end{equation}
This leads to the following formula for the posterior in the exactly same way,
as in the proof of \cref*{prop:normalForm}.
\begin{rem}
The
conditional probability of partition $\cJ$ in the finite zero-mean Gaussian model,
given the observation vector $\mathbf{x}=(x_j)_{j=1}^n$, is proportional to
\begin{equation}
K^{[|\cJ|]}C^{|\cJ|}\prod_{J\in \cJ}\frac{(1+\alpha)^{|J|-1}}{|J|^{d/2}\det R_{|J|}}\cdot
\exp\Big\{\frac{1}{2}\sum_{J\in \cJ}|J|\cdot \znorm{ R_{|J|}^{-1}R^2 \overline{\bx_J}}^2\Big\}
=:\mQ{\cJ}_K,
\end{equation}
where the notation is the same as in \cref*{prop:normalForm}.
\end{rem}
We denote the MAP partition in this model by $\hat{\cJ}_n^K(x_1,\ldots,x_n)$. We
now discuss the applicability of the results from the article to this partition.

\subsection*{\cref*{prop:MAPconv} for finite mixture model}

In this case the proof from \cref*{prop:MAPconv} remains unchanged. It used the
fact that the Chinese Restaurant Process prior on partitions depends only on the number of observations
within each cluster and hence it is not changing when we replace clusters
$I,J$ with $\tilde{I},\tilde{J}$ such that $|\tilde{I}|=|I|$,
$|\tilde{J}|=|J|$ and $\tilde{I}\cup\tilde{J}=I\cup J$. This is also the case
with the finite Chinese Restaurant Process. 

\subsection*{\cref*{prop:xmininter} for finite mixture model}

In our proof, \cref{lem:facprod} and \cref{lem:sqrtn} were used only for
\cref{prop:xmax}, but the latter is obvious for the finite mixture model since
when the number of clusters is bounded by $K$ then the number of observations in
the largest cluster is at least $n/K$. Therefore \cref{cor:centerMax} also holds
for the finite mixture model.

\smallskip
As for the proof of \cref{lem:xminbound}, in the \cref{eq:xxx1} the factor
$\frac{m!M!}{(m+M)!}$ should be replaced by 
$$\frac{K^{[m]}K^{[M]}}{K^{[m+M]}}
\cdot \frac{\alpha+m+M}{(\alpha+m)(\alpha+M)}\frac{(1+\alpha)^{(m)}(1+\alpha)^{(M)}}{(1+\alpha)^{(m+M)}}<
2\frac{(1+\alpha)^{(m)}(1+\alpha)^{(M)}}{(1+\alpha)^{(m+M)}}
$$ 
Note that
\begin{equation}
\begin{split}
\frac{(1+\alpha)^{(m)}(1+\alpha)^{(M)}}{(1+\alpha)^{(m+M)}}&=
\frac{(1+\alpha)(2+\alpha)\ldots(M+\alpha)}
{(m+1+\alpha)(m+2+\alpha)\ldots(m+M+\alpha)}=\\
&=\frac{1+\alpha}{m+1+\alpha}\cdot \frac{2+\alpha}{m+2+\alpha}\cdot
\ldots\cdot \frac{M+\alpha}{m+M+\alpha}<\\
&<\frac{1+\ceil{\alpha}}{m+1+\ceil{\alpha}}\cdot \frac{2+\ceil{\alpha}}{m+2+\ceil{\alpha}}\cdot
\ldots\cdot \frac{M+\ceil{\alpha}}{m+M+\ceil{\alpha}}=\\
&=\frac{( M+\ceil{\alpha}
)!(m+\ceil{\alpha})!}{\ceil{\alpha}!(m+M+\ceil{\alpha})!}<
\frac{( M+\ceil{\alpha} )!(m+\ceil{\alpha})!}{(m+M+\ceil{\alpha})!}.
\end{split}
\end{equation}
By Stirling inequality \eqref{eq:stirling} we get that
\begin{equation}
\begin{split}
\frac{( M+\ceil{\alpha} )!(m+\ceil{\alpha})!}{(m+M+\ceil{\alpha})!}&<
\frac{2\pi e^2 \sqrt{(M+\ca)(m+\ca)}\big((M+\ca)/e\big)^{M+\ca}\big((m+\ca)/e\big)^{m+\ca}}
{\sqrt{2\pi}\sqrt{m+M+\ca}\big(( m+M+\ca )/e\big)^{m+M+\ca}}=\\
&=\sqrt{2\pi}e^{2-\ca}\sqrt{\frac{(M+\ca)(m+\ca)}{m+M+\ca}}
\frac{(M+\ca)^{M+\ca}(m+\ca)^{m+\ca}}{( m+M+\ca )^{m+M+\ca}}<\\
&<\sqrt{2\pi}e^{2-\ca}\sqrt{m+\ca}
\frac{(m+\ca)^{m+\ca}}{( m+M+\ca )^{m}}=\\
&<\sqrt{2\pi}e^{2-\ca}( m+\ca )^{\ca+1/2}
\Big(\frac{m+\ca}{ m+M+\ca }\Big)^m.
\end{split}
\end{equation}
Therefore we can transform \eqref{eq:finalmM} into
\begin{equation}
\liminf_{n\to\infty}\frac{\P(\hat{\cJ}_n\cond\textbf{x})}{\P(\tilde{\cJ}_n\cond\textbf{x})}\leq
\liminf_{n\to\infty}\tilde{C}'\left(\frac{m+M}{mM}\right)^{d/2}( m+\ca )^{\ca+1/2}
\Big(C''\frac{m+\ca}{ m+M+\ca }\Big)^m= 0.
\end{equation}
Indeed, note that as $m/M\to 0$, we have $C''\frac{m+\ca}{ m+M+\ca }\to 0$, so
even if $m\to\infty$ so that $( m+\ca )^{\ca+1/2}\to \infty$ we have $( m+\ca
)^{\ca+1/2} \Big(C''\frac{m+\ca}{ m+M+\ca
}\Big)^m\to 0$. Therefore the proof of \cref{prop:xmininter} does not require major
changes.

\smallskip
Finally in the proof of \cref*{prop:xmininter} the only place where the prior is
important is \cref{eq:interto0}, which now become
\begin{equation}
\frac{(1+\alpha)^{ ( a+b-1 ) }(1+\alpha)^{(M-1)}}
{(1+\alpha)^{(b-1)}(1+\alpha)^{ (a+M-1) }}=\frac{(b+\alpha)^{(a)}}{(M+\alpha)^{(a)}}<
\frac{b+\alpha}{M+\alpha}\stackrel{k\to\infty}{\longrightarrow} 0.
\end{equation}
The rest of the proof of \cref*{prop:xmininter} remains unchanged.

\subsection*{\cref*{thm:distto0} for finite mixture model}

\cref*{thm:distto0} was an easy consequence of \cref*{thm:MAPmax}. It still
holds for finite mixture model if we restrict our attention to the
$P$-partitions of the observation space with at most $K$ sets (which is the number of clusters
assumed by the model).

The place of the prior distribution on the space of partitions in the proof of \cref*{thm:MAPmax} was in \cref*{lem:induced} and
\cref*{lem:approxAn}. Proof of \cref*{lem:induced} remains unchanged, provided that
\begin{equation}
\sqrt[n]{\prod_{J\in\cJ^\cA_n} (1+\alpha)^{(|J|-1)}} \approx \frac{n}{e}
\prod_{A\in\cA} p_A^{p_A},
\end{equation}
where $\cA$ is a partition of the observation space with at most $K$ sets
(otherwise clearly $Q_{ \bm{x} }(\hat{\cJ}_n^\cA)_K=0$).
Note that 
\begin{equation}
\begin{split}
\prod_{J\in\cJ^\cA_n} (1+\alpha)^{(|J|-1)}&\leq
\prod_{J\in\cJ^\cA_n} (1+\ceil{\alpha})^{(|J|-1)}=
\prod_{J\in\cJ^\cA_n} \frac{(|J|+\ceil{\alpha}-1)!}{\ceil{\alpha}!}=\\
&=(\ceil{\alpha})^{-K}\prod_{J\in\cJ^\cA_n} (|J|+\ceil{\alpha}-1)!\approx\\
&\approx (\ceil{\alpha})^{-K}\prod_{J\in\cJ^\cA_n}
\sqrt{2\pi(|J|+\ceil{\alpha}-1)}\Big(\frac{|J|+\ceil{\alpha}-1}{e}\Big)^{|J|+\ceil{\alpha}-1}=\\
&= e^{-n-|\cA|(\ceil{\alpha}-1)}(\ceil{\alpha})^{-K}\prod_{A\in \cA}
\sqrt{2\pi(|J^A_n|+\ceil{\alpha}-1)}\Big(\frac{|J^A_n|+\ceil{\alpha}-1}{e}\Big)^{|J^A_n|+\ceil{\alpha}-1}.
\end{split}
\end{equation}
Since $|J^A_n|/n\to p_A$ it is easy now to deduce that 
$\sqrt[n]{\prod_{J\in\cJ^\cA_n} (1+\alpha)^{(|J|-1)}} \approx \frac{n}{e}
\prod_{A\in\cA} p_A^{p_A}$. The proof of \cref*{lem:approxAn} can be easily
modified in the same way.

\newpage
\renewcommand{\thesection}{B}
\begin{center}
\Large\textbf{Supplement B}\\
\small Simulations
\end{center}

\vspace*{1cm}
This supplementary material provides computer simulations for considerations
presented in
\cref*{subsec:segment} of the paper together with a short simulation study of
the clustering properties of the MAP for the finite mixture input. The experiments were performed on a 64-bit
Linux machine with \texttt{R} version 3.2.3 (\cite{bib:R}). Sampling from the
posterior was performed by using MCMC methods, i.e. running a Markov Chain for
which the posterior probability is a stationary distribution. Our Markov Chain
was the one used in Algorithm 2 from \cite{bib:neal}. The choice of the MAP was
always performed by the following procedure:
firstly 100 MCMC steps were recorded after 100 burn-in
period (the initial partition is a single cluster). Then the posterior
probability of every resulting partition (up to the norming constant) was
computed and the one with the highest output was chosen as the MAP. 

Sampling from
multivariate normal and computing values of its distribution function were
performed using \texttt{mvtnorm} package (\cite{bib:mvtnorm}).

\section*{\nameref*{subsec:segment}} 
The experiment involved creating a sample of size 200
from $\Unif[-1,1]$ distribution and constructing the MAP partition. This
procedure was performed for all possible 
combinations of parameters $( \Sigma,\Tau,\alpha )\in\{1,.1,.01,.001\}\times
\{1,.1,.01\}^2$.
The results are presented in \cref{tab:difpar}.

The analysis of \cref{tab:difpar} yields to the following conclusions.
Firstly, there is an apparent connection between the observed number of clusters
in the MAP and the value predicted in \cref*{subsec:segment}, which is
$1/\sqrt{3\Sigma}$. Secondly, decreasing values of $\alpha$ and
$\Tau$ leads to the smaller number of clusters in the estimated MAP, but their
impact is significantly smaller than the impact of $\Sigma$. In case of $\alpha$
it is easily justified by the formula
\eqref{eq:posterior}, where we have the factor $\alpha^{\#\textrm{clusters}}$ in the
prior weight. The role of $\Tau$ is more difficult to explain as it occurs in
two factors: as $U^{\#\textrm{clusters}}$ and in
$R_{m}=\sqrt{\Sigma^{-1}+\Tau^{-1}/m}$. In the later it is divided by the
cluster size, so the intuition is that there is indeed positive correlation
between the value of $T$ and the number of clusters.

\section*{\nameref*{subsec:exp}} 
\begin{figure}
\centering
\includegraphics[width=.8\textwidth]{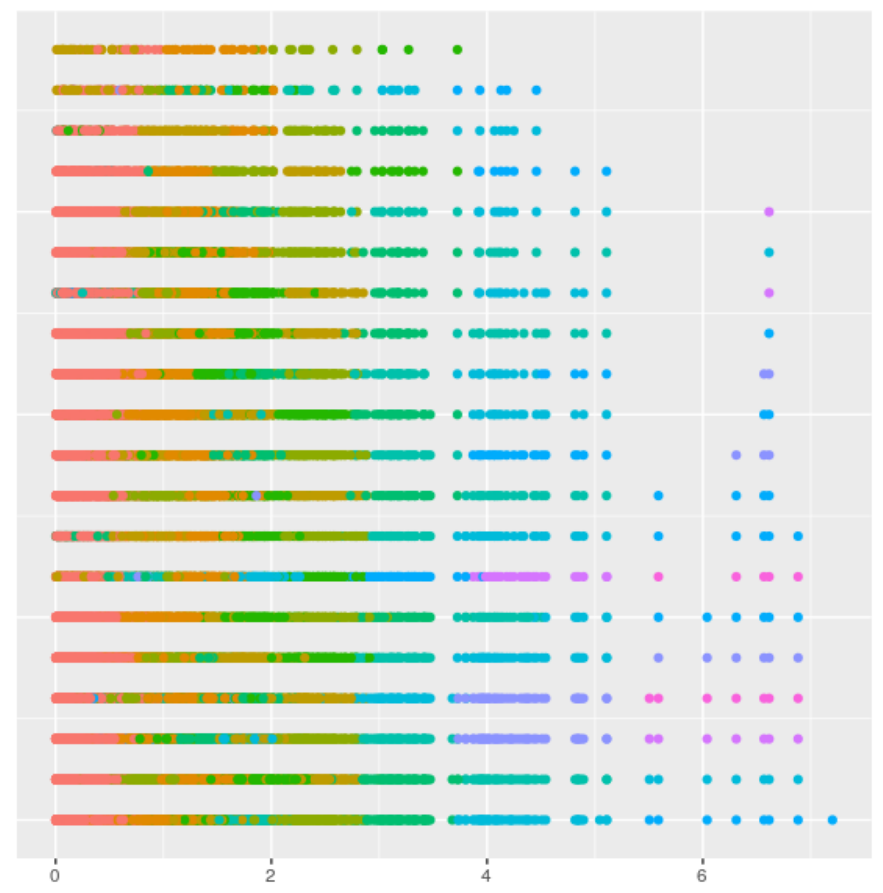}
\caption{Clustering in the MAP partition of the first $k=100,200,300,\ldots,2000$
observations in the iid sample from $\Exp(1)$. Each row corresponds to a
different value of $k$; different clusters are denoted by different colors.}
\label{fig:exponential}
\end{figure}

We sample an iid sequence from $\Exp(1)$ of size 2000. Then we construct
the MAP division of first $k$ observations for $k=100, 200, 300, \ldots, 2000$.
The parameters of the model were $\alpha=\Tau=1$, $\Sigma=(32\log 2)^{-1}$. The
results are presented in \cref{fig:exponential}, where each row corresponds to
different value of $k$ and within a row the partition is indicated by colors.

\cref{fig:exponential} is consistent with the considerations presented in
\cref{subsec:exp} regarding the number of clusters in the MAP. It suggests that at some
stage a group of extremal observations will create a new cluster and therefore
the number of clusters in the MAP tends to infinity.

\section*{\nameref*{subsec:mnorm}} 
\begin{figure}
\centering
\includegraphics[width=\textwidth]{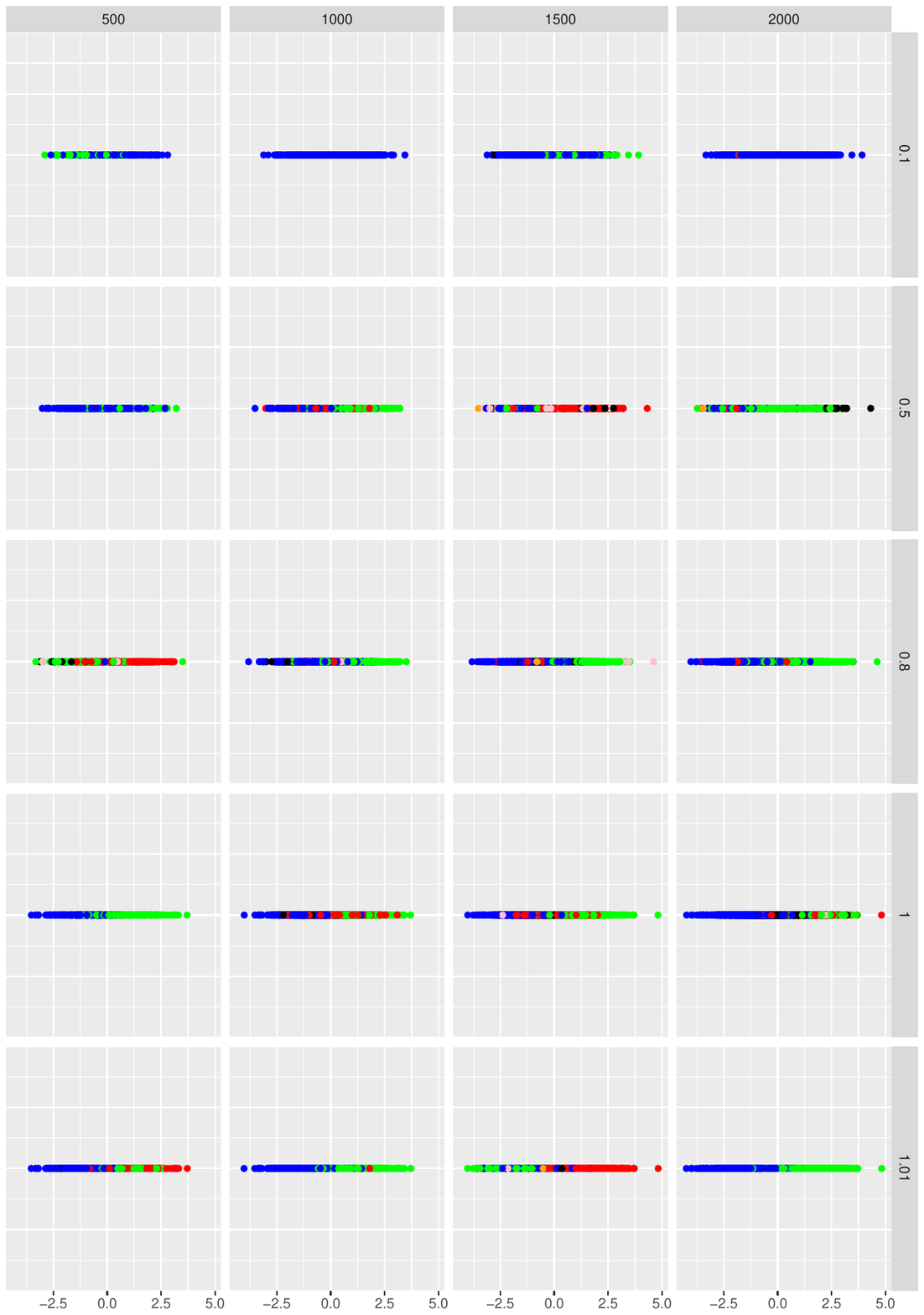}
\caption{Clustering in the MAP partition of the first $k=100,500,1000,1500,2000$
observations (in columns) in the iid sample from the mixture of two normal
distributions $\Normal(a,1)$ and $\Normal(-a,1)$ (in rows). Different clusters are denoted by different colors.}
\label{fig:normal}
\end{figure}

In this numerical experiment firstly we sample an iid sequence $X_1,\ldots,
X_{2000}\sim\Normal(0,1)$ and an iid sequence of Rademacher random variables
$R_1,\ldots,R_{2000}\sim \re{2}(\delta_1+\delta_{-1})$. 
Note that for every $a\in \R$ the distribution of the random variable
$Z_{a,i}=X_i+aY_i$ is
the mixture of two normals with means $a$ and $-a$.  
The MAP partition of the vector
$Z_{a,1},\ldots,Z_{a,k}$ is computed for $a\in\{.1, .5, .8, 1\}$ and $k\in \{500, 1000,
1500, 2000\}$. The parameters of the mode were $\alpha=\Sigma=\Tau=1$. Note that
the choice of $\Sigma$ parameter is consistent with the variance within
mixtures. The results are shown on \cref{fig:normal}.

\begin{figure}
\centering
\includegraphics[width=\textwidth]{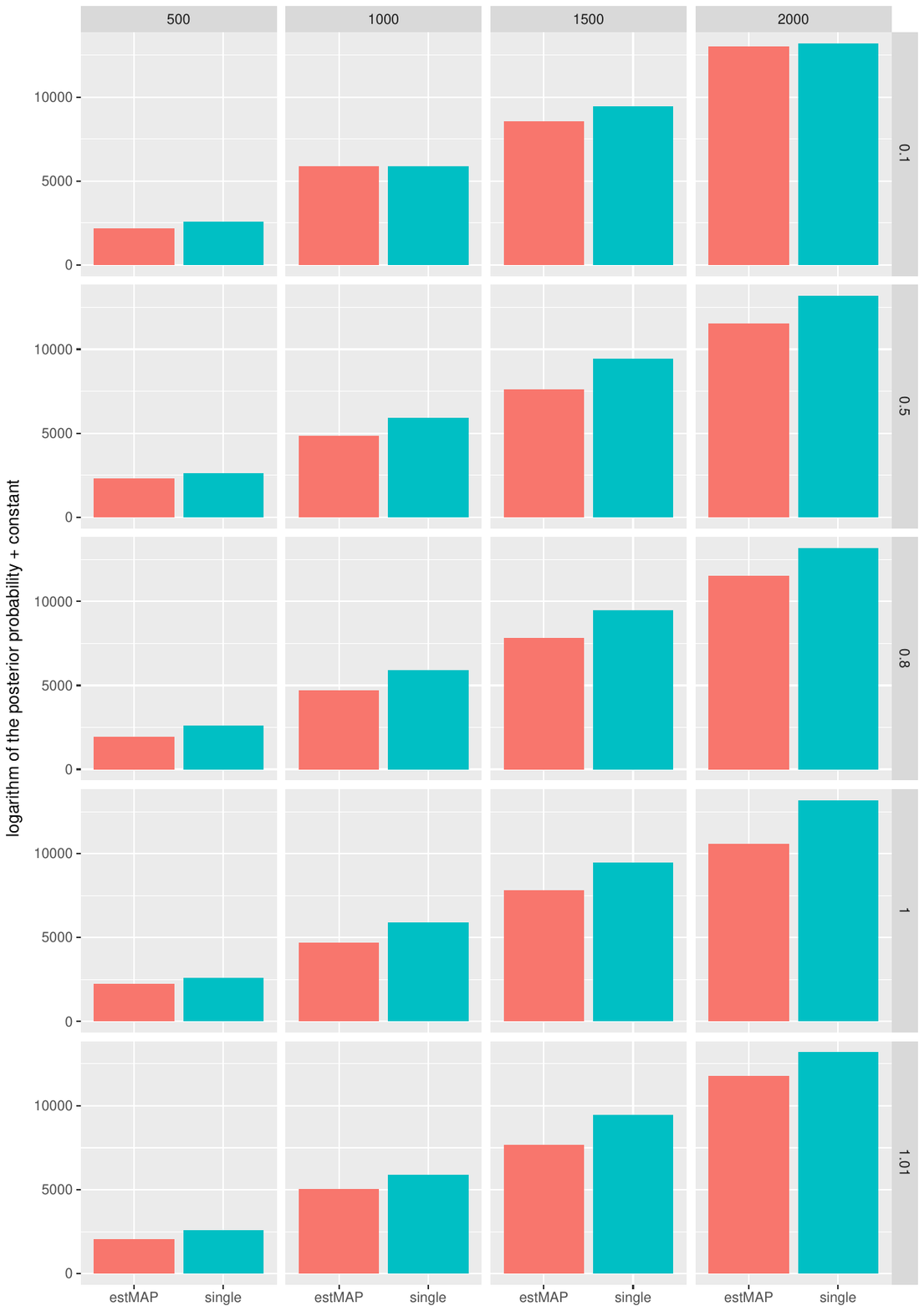}
\caption{Comparison of posterior probabilities of the approximation of the MAP
(as in \cref{fig:normal}) and the single cluster. In all cases single cluster
has higher posterior probability.}
\label{fig:normalProbs}
\end{figure}

\cref*{subsec:mnorm} predicts that with this choice of $a$ parameters, a single cluster partition has
larger posterior probability than
splitting the observations into two clusters of equal size. \cref{fig:normal}
does not fully support these predictions. However this is due to the
imperfection of the MAP MCMC approximation. Indeed, \cref{fig:normalProbs} shows
the comparison of the logarithms of posterior probabilities of the sample
estimators of MAP and single cluster partition. In all cases the single cluster
partition has higher posterior probability.

\begin{table}
\centering\small
\begin{tabular}{cccclc}
$\Sigma$& $\Tau$& $\alpha$& sizes of clusters &  \# clusters & est \\
\hline\\
1 & 1 & 1 & 200 & 1 & 0.58\\
1 & 1 & 0.1 & 200 & 1 & 0.58\\
1 & 1 & 0.01 & 200 & 1 & 0.58\\
1 & 0.1 & 1 & 180, 15, 3, 1, 1 & 5 & 0.58\\
1 & 0.1 & 0.1 & 200 & 1 & 0.58\\
1 & 0.1 & 0.01 & 200 & 1 & 0.58\\
1 & 0.01 & 1 & 200 & 1 & 0.58\\
1 & 0.01 & 0.1 & 200 & 1 & 0.58\\
1 & 0.01 & 0.01 & 200 & 1 & 0.58\\
0.1 & 1 & 1 & 110, 90 & 2 & 1.83\\
0.1 & 1 & 0.1 & 106, 94 & 2 & 1.83\\
0.1 & 1 & 0.01 & 105, 95 & 2 & 1.83\\
0.1 & 0.1 & 1 & 105, 95 & 2 & 1.83\\
0.1 & 0.1 & 0.1 & 121, 79 & 2 & 1.83\\
0.1 & 0.1 & 0.01 & 110, 90 & 2 & 1.83\\
0.1 & 0.01 & 1 & 116, 84 & 2 & 1.83\\
0.1 & 0.01 & 0.1 & 103, 97 & 2 & 1.83\\
0.1 & 0.01 & 0.01 & 200 & 1 & 1.83\\
0.01 & 1 & 1 & 52, 40, 30, 30, 28, 18, 2 & 7 & 5.77\\
0.01 & 1 & 0.1 & 35, 34, 31, 31, 30, 20, 19 & 7 & 5.77\\
0.01 & 1 & 0.01 & 59, 52, 51, 38 & 4 & 5.77\\
0.01 & 0.1 & 1 & 46, 42, 37, 36, 22, 12, 4, 1 & 8 & 5.77\\
0.01 & 0.1 & 0.1 & 49, 48, 43, 37, 23 & 5 & 5.77\\
0.01 & 0.1 & 0.01 & 51, 45, 43, 31, 30 & 5 & 5.77\\
0.01 & 0.01 & 1 & 56, 44, 37, 34, 27, 1, 1 & 7 & 5.77\\
0.01 & 0.01 & 0.1 & 55, 44, 38, 35, 28 & 5 & 5.77\\
0.01 & 0.01 & 0.01 & 51, 48, 45, 43, 13 & 5 & 5.77\\
0.001 & 1 & 1 & 20, 17, 16, 14, 14, 13, 13, 12, 12, 12, 12, 11, 11, 10, 8, 5 & 16 & 18.26\\
0.001 & 1 & 0.1 & 20, 19, 18, 16, 16, 16, 15, 14, 14, 13, 12, 11, 9, 7 & 14 & 18.26\\
0.001 & 1 & 0.01 & 21, 20, 17, 17, 17, 16, 14, 14, 13, 12, 11, 10, 9, 9 & 14 & 18.26\\
0.001 & 0.1 & 1 & 20, 18, 17, 16, 16, 16, 15, 14, 13, 13, 12, 9, 7, 6, 5, 3 & 16 & 18.26\\
0.001 & 0.1 & 0.1 & 21, 20, 18, 15, 15, 13, 13, 12, 11, 10, 10, 9, 9, 9, 8, 7 & 16 & 18.26\\
0.001 & 0.1 & 0.01 & 25, 21, 17, 17, 16, 16, 15, 14, 13, 13, 12, 11, 10 & 13 & 18.26\\
0.001 & 0.01 & 1 & 27, 25, 24, 23, 23, 21, 19, 14, 13, 10, 1 & 11 & 18.26\\
0.001 & 0.01 & 0.1 & 21, 20, 20, 19, 18, 18, 17, 16, 13, 12, 10, 9, 7 & 13 & 18.26\\
0.001 & 0.01 & 0.01 & 33, 28, 25, 22, 20, 20, 16, 14, 11, 11 & 10 & 18.26\\
\end{tabular}
\caption{\footnotesize Results of numerical experiments with all combination of
$\Sigma,\Tau,\alpha\in\{1,.1,.01,.001\}\times\{1,.1,.01\}^2$. Column `sizes of
clusters` presents the sizes of the clusters
created in the MAP, sorted in decreasing order. Column `est' is equal to
$1/(\sqrt{3\Sigma})$ and should approximate the number of clusters, which is
given in column `\# clusters`.}
\label{tab:difpar}
\end{table}

\end{document}